\documentclass[12pt]{amsart}
\usepackage[utf8]{inputenc}
\usepackage{MnSymbol}
\usepackage{tikz-cd}

\usepackage{bm}
\usepackage[textsize=footnotesize,shadow]{todonotes}
\usepackage{fancyhdr}
\usepackage[colorlinks=true,unicode=true]{hyperref}
\usepackage{lpic}

\usepackage{environ}
\usepackage{pgffor}
\newcommand{\skipthis}[2][]{\relax}

\def\mkmathletter#1#2{%
    \expandafter\gdef\csname#1#2\endcsname%
           {\ensuremath{\csname math#2\endcsname{#1}}}%
}
\def\mkmathletters#1{%
\mkmathletter{A}{#1}%
\mkmathletter{B}{#1}%
\mkmathletter{C}{#1}%
\mkmathletter{D}{#1}
\mkmathletter{E}{#1}
\mkmathletter{F}{#1}
\mkmathletter{G}{#1}
\mkmathletter{H}{#1}
\mkmathletter{I}{#1}
\mkmathletter{J}{#1}
\mkmathletter{K}{#1}
\mkmathletter{L}{#1}
\mkmathletter{M}{#1}
\mkmathletter{N}{#1}
\mkmathletter{O}{#1}
\mkmathletter{P}{#1}
\mkmathletter{Q}{#1}
\mkmathletter{R}{#1}
\mkmathletter{S}{#1}
\mkmathletter{T}{#1}
\mkmathletter{U}{#1}
\mkmathletter{V}{#1}
\mkmathletter{W}{#1}
\mkmathletter{X}{#1}
\mkmathletter{Y}{#1}
\mkmathletter{Z}{#1}
\mkmathletter{a}{#1}
\mkmathletter{b}{#1}
\mkmathletter{c}{#1}
\mkmathletter{d}{#1}
\mkmathletter{e}{#1}
\mkmathletter{f}{#1}
\mkmathletter{g}{#1}
\mkmathletter{h}{#1}
\mkmathletter{i}{#1}
\mkmathletter{j}{#1}
\mkmathletter{k}{#1}
\mkmathletter{l}{#1}
\mkmathletter{m}{#1}
\mkmathletter{n}{#1}
\mkmathletter{o}{#1}
\mkmathletter{p}{#1}
\mkmathletter{q}{#1}
\mkmathletter{r}{#1}
\mkmathletter{s}{#1}
\mkmathletter{t}{#1}
\mkmathletter{u}{#1}
\mkmathletter{v}{#1}
\mkmathletter{w}{#1}
\mkmathletter{x}{#1}
\mkmathletter{y}{#1}
\mkmathletter{z}{#1}
}
\mkmathletters{frak}
\mkmathletters{bb}
\mkmathletters{cal}
\mkmathletters{bf}
\mkmathletters{rm}
\mkmathletters{sc}
\mkmathletters{sf}

\def\Gammabf{\ensuremath{\bm\Gamma}}

\def\Lambdabf{\ensuremath{\bm\Lambda}}

\newcounter{remarkcounter}

\def\st{\,\middle\vert\mkern-3mu\middle|\,}
\def\set#1{\left\{#1\right\}}
\def\choose#1{
\begin{pmatrix}
        #1
\end{pmatrix}
}
\def\<#1>{\left\langle#1\right\rangle}
\def\|#1|{\left|\left|#1\right|\right|}

\def\pre#1#2#3#4{}

\let\termtexthook=\relax
\let\termmarginhook=\relax
\let\seehook=\relax
\def\termmargin{no}

\def\termtexthook{\em}
\def\termmarginhook#1{\color{blue}\fbox{\parbox{\marginparwidth}{\raggedright#1}}}
\def\seehook#1{\color{gray}\underline{#1}}

\def\allterms{}

\newcommand{\term}[2][EMPTY]{%
  \ifthenelse{\equal{#1}{EMPTY}}%
    {\hypertarget{t:#2}{{\termtexthook{#2}}}%
      \xdef\allterms{\allterms\ #2 \thepage,}}%
    {\hypertarget{t:#1}{{\termtexthook{#2}}}%
      \xdef\allterms{\allterms\ #1 \thepage,}}%
  \ifthenelse{\equal{\termmargin}{yes}}
    {\marginpar{\termmarginhook{#2}}}%
    {\relax}%
}

\renewcommand{\see}[2][EMPTY]{%
  \hypersetup{linkcolor=blue}      
  \ifthenelse{\equal{#1}{EMPTY}}%
    {\hyperlink{t:#2}{{\seehook{#2}}}}%
    {\hyperlink{t:#1}{{\seehook{#2}}}}%
  \hypersetup{linkcolor=red}  
}  

\makeatletter\@addtoreset{equation}{section}\makeatother
\def\defaulteqtag{\thesection.\arabic{equation}}

\newcommand{\tageq}[2][empty]{%
   \ifthenelse{\equal{#1}{empty}}
     {%
      \stepcounter{equation}%
      \tag{\defaulteqtag}%
     }
     {%
     \tag{#1}%
     }%
     \label{eq:#2}%
}

\def\empty{}
\xdef\listofclaims{}
\def\globallet#1#2{\global\let#1#2}
\def\epropsymbol{$\boxtimes$}
\def\eprop{ \rule{0.00mm}{3mm}\nolinebreak\hfill\epropsymbol}
\def\printlabel#1{\makebox[0mm][l]{%
    \hspace{-3mm}\raisebox{1em}{\tiny\color{red}#1}}}
\def\Label#1{\printlabel{#1}\label{#1}}
\newcommand{\NewTheorem}[3]{
  \newtheorem{#1nosave}[#3]{#2}
  \NewEnviron{#1}[1]
  {
    \xdef\lastclaim{##1}%
    \expandafter\xdef\csname##1countername\endcsname{the#3}
    \expandafter\xdef\csname##1envname\endcsname{#1nosave}
    \expandafter\globallet\csname##1content\endcsname=\BODY
    \xdef\listofclaims{\listofclaims,##1}
    \begin{#1nosave}\label{##1}%
      \expandafter\xdef\csname##1number\endcsname{\csname
        the#3\endcsname}
        \expandafter\ifx\csname r@##1.rep\endcsname\relax\else
          \hyperref[##1.rep]{\boldmath$\downarrow$}
        \fi          
        \BODY\eprop
    \end{#1nosave}
  }
}

\newif\ifclaimismissing
\claimismissingfalse

\newcommand{\repeatclaim}[1]{%
  \xdef\lastclaim{#1}%
   \expandafter\ifx\csname#1content\endcsname\relax
      {\mbox{\relax}\par\noindent\color{red}\bf Missing claim
      ``#1''.\par%
      \global\claimismissingtrue%
      }
   \else
     \global\claimismissingfalse%
     \begingroup%
     \expandafter\def\csname\csname#1countername\endcsname\endcsname{%
       \csname#1number\endcsname}
     \begin{\csname#1envname\endcsname}\label{#1.rep}
         \hyperref[#1]{\boldmath$\uparrow$} 
       \let\label=\skipthis
       \let\Label=\skipthis
       \let\marginpar=\skipthis
       \let\todo=\skipthis
       \csname#1content\endcsname
       \eprop
     \end{\csname#1envname\endcsname}
   \endgroup%
   \fi
}

\def\repeatallclaimsaux,#1;EndOfTheArgs{
  \foreach \claim in {#1}{
    \repeatclaim{\claim}
  }
}

\def\repeatallclaims{%
  \ifx\listofclaims\empty%
    \relax%
  \else%
    \expandafter\repeatallclaimsaux\listofclaims;EndOfTheArgs%
  \fi%
}

\NewTheorem{definition}{Definition}{withinsection}
\NewTheorem{proposition}{Proposition}{withinsection}
\NewTheorem{theorem}{Theorem}{withinsection}
\NewTheorem{corollary}{Corollary}{withinsection}
\NewTheorem{lemma}{Lemma}{withinsection}
\NewTheorem{tlemma}{Lemma}{withinsection1}

\def\eproofsymbol{$\blacklozenge$}
\def\eproofsymbol{$\boxtimes$}

\NewEnviron{Proof}[1][\lastclaim]{%
  \expandafter\ifx\csname#1content\endcsname\relax
     {\mbox{\relax}\par\color{red}\bf Skipping proof of ``#1''.%
      \nolinebreak\hfill\eproofsymbol\par}%
   \else
     \par\noindent\textit{Proof:} \BODY%
     \rule{0.00mm}{3mm}\nolinebreak\hfill\eproofsymbol\par
   \fi  
}
\def\eproof{\rule{0.00mm}{3mm}\nolinebreak\hfill\eproofsymbol\par}

\def\skippar#1\par{\relax\par}
\def\justpar#1\par{#1\par}

\def\specialpar#1#2#3#4\par{%
  \rule{0mm}{3mm}\par%
  #3%
  \textcolor{#1}{%
    \makebox[0mm][r]{#2}\unskip%
    {#4}%
  }%
  \par%
}

\def\csletsecond#1#2{%
  \expandafter\let\expandafter#1\csname#2\endcsname}

\def\cslet#1#2{%
  \expandafter\csletsecond\csname#1\endcsname{#2}}

\newcommand{\getoption}[4]{%
  \def\get@option##1#2=##2,##3ENDOFTHEARGS{%
    \expandafter\gdef\csname#1#2\endcsname{##2}}%
  \get@option#4,#2=#3,ENDOFTHEARGS}

\newcommand{\setcolors}[2][]{%
  \getoption{#2}{rmcolor}{green}{#1}%
  \getoption{#2}{inscolor}{blue}{#1}%
  \getoption{#2}{commentcolor}{red}{#1}%
  \getoption{#2}{symbol}{$bullet\,\,$}{#1}%
  \getoption{#2}{commentsymbol}{off}{#1}%
  \getoption{#2}{rmsymbol}{off}{#1}%
  \getoption{#2}{inssymbol}{off}{#1}%
  \expandafter\gdef\csname#2@par\endcsname{%
    \specialpar{\csname#2commentcolor\endcsname}%
               {}%
               {\noindent}%
    \ifthenelse{\equal{\csname#2commentsymbol\endcsname}{on}}{%
    \marginpar{\center\textcolor{\csname#2commentcolor\endcsname}{%
               \csname#2symbol\endcsname}}}{\relax}%
  }
  \expandafter\gdef\csname#2@inspar\endcsname{%
    \specialpar{\csname#2inscolor\endcsname}%
               {}%
               {}%
    \ifthenelse{\equal{\csname#2inssymbol\endcsname}{on}}{%
    \marginpar{\center\textcolor{\csname#2inscolor\endcsname}{%
               \csname#2symbol\endcsname}}}{\relax}%
  }
  \expandafter\gdef\csname#2@rmpar\endcsname{%
    \specialpar{\csname#2rmcolor\endcsname}%
               {}%
               {}%
    \ifthenelse{\equal{\csname#2rmsymbol\endcsname}{on}}{%
    \marginpar{\center\textcolor{\csname#2rmcolor\endcsname}{%
               \csname#2symbol\endcsname}}}{\relax}%
  }
  \expandafter\def\csname#2@rm\endcsname{\relax}
  \expandafter\renewcommand\csname#2@rm\endcsname[2][]{%
    \ifthenelse{\equal{\csname#2rmsymbol\endcsname}{on}}{%
    \marginpar{\center\textcolor{\csname#2rmcolor\endcsname}{%
               \csname#2symbol\endcsname}}}{\relax}%
    \textcolor{\csname#2inscolor\endcsname}{##1}%
    \ifthenelse{\equal{##1}{}}{\relax}{/}%
    \textcolor{\csname#2rmcolor\endcsname}{##2}%
  }%
  \expandafter\def\csname#2@ins\endcsname{\relax}
  \expandafter\renewcommand\csname#2@ins\endcsname[2][]{%
    \ifthenelse{\equal{\csname#2inssymbol\endcsname}{on}}{%
    \marginpar{\center\textcolor{\csname#2inscolor\endcsname}{%
               \csname#2symbol\endcsname}}}{\relax}%
    \textcolor{\csname#2rmcolor\endcsname}{##1}%
    \ifthenelse{\equal{##1}{}}{\relax}{/}%
    \textcolor{\csname#2inscolor\endcsname}{##2}%
  }%
  \expandafter\gdef\csname#2on\endcsname{%
    \cslet{#2par}{#2@par}%
    \cslet{#2inspar}{#2@inspar}%
    \cslet{#2rmpar}{#2@rmpar}%
    \cslet{#2ins}{#2@ins}%
    \cslet{#2rm}{#2@rm}%
  }%
  \expandafter\gdef\csname#2off\endcsname{%
    \cslet{#2par}{skippar}%
    \cslet{#2inspar}{justpar}%
    \cslet{#2rmpar}{skippar}%
    \cslet{#2ins}{justins}%
    \cslet{#2rm}{justrm}%
  }%
  \csname#2on\endcsname
} 

\input{defs.tex}

\xdef\SectionName{*}
\xdef\SubsectionName{}
\let\Section=\section
\renewcommand{\section}[2][empty]{%
	\xdef\SubsectionName{}%
	\ifthenelse{\equal{#1}{empty}}%
	{\xdef\SectionName{#2}%
		\Section{#2}}%
	{\xdef\SectionName{#1}%
		\Section[#1]{#2}}%
}

\let\Subsection=\subsection
\renewcommand{\subsection}[2][empty]{%
	\ifthenelse{\equal{#1}{empty}}%
	{\xdef\SubsectionName{#2}%
		\Subsection{#2}}%
	{\xdef\SubsectionName{#1}%
		\Subsection[#1]{#2}}%
}

\def\draft{\fbox{\tiny draft: \today}}
\let\draft=\relax
\setcolors[rmcolor=cyan,
           commentcolor=purple,
           inscolor=blue,
           symbol=$\boxed{\mathbb{S}}\,$,
           commentsymbol=on,
           rmsymbol=on,
           inssymbol=on]
          {s}
\setcolors[rmcolor=green,
           commentcolor=red,
           inscolor=brown,
           symbol=$\boxed{\mathbb{J}}\,$,
           commentsymbol=on,
           rmsymbol=on,
           inssymbol=on]
          {j}

\title[Tropical Limits]{Tropical Limits of Probability Spaces, 
       Part I
       \\[1ex] 
       \footnotesize\mdseries 
       The Intrinsic Kolmogorov-Sinai Distance and the 
       \\ 
       Asymptotic Equipartition Property for Configurations} 

\author[RM]{R. Matveev}
\author[JWP]{J. W. Portegies}

\begin{document}
	\thispagestyle{fancy} 
	\begin{abstract}
The entropy of a finite probability space $X$ measures the
observable cardinality of large independent products $X^{\otimes n}$
of the probability space. If two probability spaces $X$ and $Y$ have
the same entropy, there is an almost measure-preserving bijection
between large parts of $X^{\otimes n}$ and $Y^{\otimes n}$. In this
way, $X$ and $Y$ are asymptotically equivalent.

It turns out to be challenging to generalize this notion of asymptotic
equivalence to configurations of probability spaces, which are
collections of probability spaces with measure-preserving maps between
some of them.

In this article we introduce the intrinsic Kolmogorov-Sinai distance
on the space of configurations of probability spaces. Concentrating on
the large-scale geometry we pass to the asymptotic Kolmogorov-Sinai
distance.  It induces an asymptotic equivalence relation on sequences
of configurations of probability spaces. We will call the equivalence
classes \emph{tropical probability spaces}.

In this context we prove an Asymptotic Equipartition Property for
configurations. It states that tropical configurations can always be
approximated by homogeneous configurations.  In addition, we show that
the solutions to certain Information-Optimization problems are
Lipschitz-con\-tinuous with respect to the asymptotic Kolmogorov-Sinai
distance.  It follows from these two statements that in order to solve
an Information-Optimization problem, it suffices to consider
homogeneous configurations. 

Finally, we show that spaces of trajectories of length $n$ of certain stochastic
processes, in particular stationary Markov chains, have a tropical limit.
\end{abstract}

	\maketitle
	
	
	\setcounter{section}{-1}
	\section{Introduction}\label{s:introduction}
	  The aim of the present article is to develop a theory of tropical
  probability spaces, which are asymptotic classes of finite
    probability spaces.  Together with the accompanying techniques,
  we expect them to be relevant to problems arising in information
  theory, causal inference, artificial intelligence and neuroscience.

  As a matter of introduction and motivation of the research presented
  in the article, we start by considering a few simple examples.

\subsection{Single probability spaces}\label{s:intro-single}
  We consider a {finite probability space} $X = (S, p)$, where $S$
  is a finite set, and $p$ is a probability measure on $S$. For
  simplicity, assume for now that the measure $p$ has full
  support. Next, we consider the, so-called, Bernoulli sequence of
  {probability spaces}
  \[
  X^{\otimes n} = ( S^n, p^{\otimes n})
  \] 
  where $S^n$ denotes the $n$-fold Cartesian product of $S$, and
  $p^{\otimes n}$ is the $n$-fold product measure. 
  
  This situation arises in several contexts.  For example, in physics,
  $X^{\otimes n}$ would encode the state of the system comprised of many
  identical non-interacting (weakly interacting) subsystems with state
  space $X$.  In information theory, $X^{\otimes n}$ would describe the
  output of an i.i.d.~random source.  In dynamical systems or stochastic
  processes the setting corresponds to Bernoulli shifts and Bernoulli
  processes.
  
  The \emph{entropy} of $X$ is the exponential growth rate of the
  \emph{observable} cardinality of tensor powers of $X$.  The
  observable cardinality, loosely speaking, is the cardinality of the
  set $X^{\otimes n}$ after (the biggest possible) set of small measure
  of elements, each with negligible measure, has been removed.  It turns
  out that the observable cardinality of $X^{\otimes n}$ might be much
  smaller than $|S|^{n}$, the cardinality of the whole of $X^{\otimes
    n}$, in the following sense.

  The \emph{Asymptotic Equipartition Property} states that for every
  $\epsilon>0$ and sufficiently large $n$ one can find a, so-called,
  \emph{typical subset} $A^{\c n}_{\epsilon}\subset S^{n}$, such that
  it takes up almost all of the mass of $X^{\otimes n}$ and the
  probability distribution on $A^{\c n}_{\epsilon}$ is almost uniform
  on the normalized logarithmic scale,
  
  \begin{enumerate}
  \item\label{intro-aep1}
    $p^{\otimes n}(A^{\c n}_{\epsilon})\geq 1-\epsilon$
  \item\label{intro-aep2}
    For any $a,a'\in A^{\c n}_{\epsilon}$ holds 
    $\left|\frac1n \ln p(a)-\frac1n\ln p(a')\right|\leq\epsilon$
  \end{enumerate}

  The cardinality $|A^{\c n}_{\epsilon}|$ may be much smaller than
  $|S|^{n}$, but it will still grow exponentially with $n$.  Even
  though there are many choices for such a set $A^{\c n}_\epsilon$,
  the exponential growth rate with respect to $n$ is well-defined upto
  $2\epsilon$. In fact, there exists a number $h_{X}$ such that for
  any choice of the typical subset $A^{(n)}_{\epsilon}$ holds
\[
\ebf^{n \cdot h_X - \epsilon} \leq |A_\epsilon^{(n)}| \leq \ebf^{n \cdot h_X + \epsilon}
\]

The limit of the
growth rate as $\epsilon\to0+$ is called the entropy of $X$, as
explained in more detail in Section~\ref{s:category-entropy}
\[
\Ent(X) 
:= 
\lim_{\epsilon \downarrow 0} \lim_{n \to \infty} \frac{1}{n} \ln |A_\epsilon^{(n)}|
\]
By the law of large numbers
\[
\ent(X) 
= 
h_X 
= 
- \sum_{x \in S} p(x) \ln p(x)
\]
which is the formula by which the Shannon entropy is usually introduced. 

Entropy is especially easy to evaluate if the space is uniform, since
for any {finite probability space} with the uniform distribution holds
\[\tageq{entropy-uniform}
  \ent(X)=\ln|X|      
\]  

This point of view on entropy goes back to the original idea of Boltzmann,
according to which entropy is the logarithm of the number of
equiprobable states, that a system, comprised of many identical weakly
interacting subsystems, may take on.

\subsubsection{Asymptotic equivalence}\label{s:intro-single-ae}
The Asymptotic Equipartion Property implies that the sequence
$X^{\otimes n}$ is asymptotically equivalent to a sequence of uniform spaces in
the following sense. Let us denote by $p_U$ the probability
distribution that is supported on $A_\epsilon^{(n)}$, and is uniform
on its support. Then a sequence from independent samples according to $p_U$
is very hard to discriminate from a sequence of independent samples
from $p_X^{\otimes n}$.

Similarly, when $X$ and $Y$ are probability spaces with the same entropy, 
the sets $X$ and $Y$ are asymptotically equivalent in the sense that there is 
a bijection between the typical sets, which can be seen as a change of code. 
This is essentially the content of Shannon's source coding theorem.

In \cite{Gromov-Search-2012}, Gromov proposed this existence of an 
``almost-bijection'' as a basis of an
asymptotic equivalence relation on sequences of probability spaces.
Even though we were greatly influenced by ideas in
\cite{Gromov-Search-2012}, we found that Gromov's definition does not
extend easily to \textbf{configurations} of probability spaces. 

By a configuration of probability spaces we mean a collection of probability spaces with 
measure-preserving maps between some of them. We will give a precise definition in Section \ref{s:category-config}, but will consider some particular examples below.
Formalizing
and studying a notion of asymptotic equivalence for configurations of
probability spaces is the main topic of the present article.

\subsection{Configurations of probability spaces}
\label{s:intro-2fan}
Suppose that now instead of a single probability space, we consider
a pair of probability spaces $X=(\un X,p_{X})$ and $Y=(\un Y,p_{Y})$ with a
joint distribution, that is a probability measure on $\un X\times \un Y$ that
pushes forward to $p_{X}$ and $p_{Y}$ under coordinate projections.
In other words, we consider a triple of probability spaces $X$, $Y$
and $U$ with a pair of measure-preserving maps $U\to X$ and $U\to
Y$. This is what we later call a minimal two-fan of probability spaces
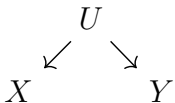
\[  
\begin{tikzcd}[row sep=small,column sep=tiny,ampersand replacement=\&]
  \mbox{}
  \&
  U
  \arrow{dl}{}
  \arrow{dr}{}
  \&
  \\
  X
  \&
  \&
  Y
\end{tikzcd}
\]  
and is a particular instance of a configuration of probability spaces.
 
\begin{figure}
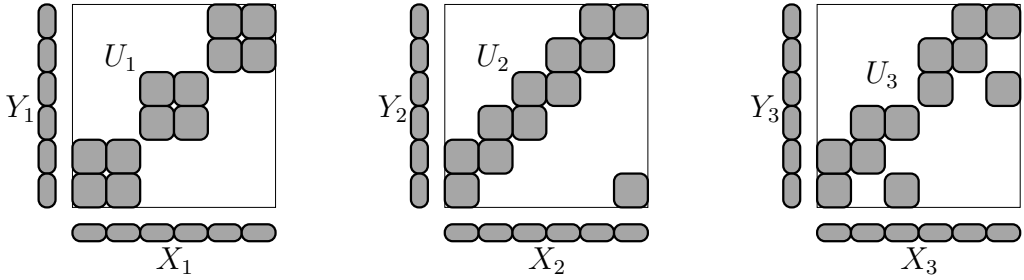

  \begin{lpic}[l(6mm),b(5mm),draft,clean]{squares-snakes(0.45)}
    \lbl[r]{0,40; $Y_{1}$}
    \lbl[t]{41,-1; $X_{1}$}
    \lbl{25,55;$U_{1}$}
    \lbl[r]{110,40; $Y_{2}$}
    \lbl[t]{151,-1; $X_{2}$}
    \lbl{135,55;$U_{2}$}
    \lbl[r]{220,40; $Y_{3}$}
    \lbl[t]{261,-1; $X_{3}$}
    \lbl{250,50;$U_{3}$}
  \end{lpic}
  \caption{Examples of pairs of probability spaces together with a
    joint distribution.}
  \label{f:squares-snakes}
\end{figure}

\subsubsection{Three examples}
Three examples of such an object are shown on Figure
\ref{f:squares-snakes}, which is to be interpreted in the following
way. Each of the spaces $X_{i}$ and $Y_{i}$, $i=1,2,3$, have
cardinality six and a uniform distribution, where the weight of each
atom is $\frac16$.  The spaces $U_{i}$, $i=1,2,3$, have cardinality
$12$ and the distribution is also uniform with all weights being
$\frac{1}{12}$. The support of the measure on $U_{i}$'s is colored
grey on the pictures.  The maps from $U_{i}$ to $X_{i}$ and $Y_{i}$
are coordinate projections.

In view of equation (\ref{eq:entropy-uniform}) we have for each
$i=1,2,3$,
\begin{align*}
  \ent(X_{i})&=\ln 6\\
  \ent(Y_{i})&=\ln 6\\
  \ent(U_{i})&=\ln 12
\end{align*}

Now we would like to ask the following.\\
\textbf{Question}. 
\begin{quote}\em
  Is it possible to find an almost-bijection between sufficiently high
  powers of $(X_{i}^{\otimes n}\ot Z_{i}^{\otimes n}\to Y_{i}^{\otimes
    n})$ by $(X_{j}^{\otimes n}\ot Z_{j}^{\otimes n}\to Y_{j}^{\otimes
    n})$ for $i\neq j$ with an arbitrary given precision as in
  Shannon's coding theorem as described at the end of the previous
  subsection~\ref{s:intro-single}?  More generally, what is the proper
  generalization of an asymptotic equivalence relation as discussed in
  the previous subsection to sequences of tensor powers of two-fans?
\end{quote}

  We would like to argue that even though the entropies of the
  constituent spaces are all (pairwise) the same, all three examples
  above should be pairwise asymptotically different.

  To establish that the examples in Figure~\ref{f:squares-snakes} are
  different, that is, not isomorphic (see also Section
  \ref{s:category-config}) is relatively easy, since they have
  non-isomorphic symmetry groups.  However, we present a different
  argument, that lends itself for generalization to prove that the
  examples at hand are not \emph{asymptotically} equivalent and
  that also gives a quantitative difference between them.

  To distinguish Example 1 from both 2 and 3, one could argue along
  the following lines. We could try to add a third space $Z=(\un
  Z,p_{Z})$ to the pair $X$ and $Y$ and provide a \emph{joint
    distribution} $p_{Q}$ on
  \[
  Q=(\un X\times\un Y\times\un Z,p_{Q})
  \]
  such that the projection of $p_{Q}$ on the first two factors is
  $p_{U}$ and on the third factor is $p_{Z}$.
  
  Once we do that, we could evaluate entropies of various push-forwards
  of $p_{Q}$. Denote by $V=(\un X\times \un Z,p_{V})$ and $W=(\un Y\times
  \un Z,p_{W})$, where $p_{V}$ and $p_{W}$ are push-forwards of $p_{Q}$
  under corresponding coordinate projections. All the probability spaces
  now fit into a commutative diagram
  
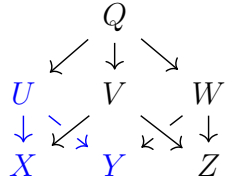
\[
  \begin{tikzcd}[row sep=small,column sep=small,ampersand replacement=\&]
    \&
    Q
    \arrow{dl}{}
    \arrow{d}{}
    \arrow{dr}{}
    \&
    \\
    \blue{U}
    \arrow[color=blue]{d}{}
    \arrow[color=blue]{dr}{}
    \&
    V
    \&
    W
    \arrow{dl}{}
    \arrow{d}{}
    \\
    \blue{X}
    \arrow[leftarrow, crossing over]{ur}{}
    \&
    \blue{Y}
    \&
    Z
    \arrow[leftarrow, crossing over]{ul}{}
    \\
  \end{tikzcd}
  \]
  where each arrow is a reduction, which is simply a measure-preserving
  map between probability spaces.

We consider the set of all possible extensions of the above form
and denote it by $\Ext(X,Y,U)$.
For any extension $\Ebf=(X,Y,Z,U,V,W,Q)$ in $\Ext(X,Y,U)$ we have
four ``new'' entropies
\[\tageq{extension-entropies}
\ent(Q),\quad\ent(V),\quad\ent(W),\quad\ent(Z)
\]
in addition to the ``known'' entropies of $X$, $Y$ and $U$.
The vector 
\[
\ent_{*}(\Ebf):=\big(\ent(X),\ent(Y),\ent(Z),\ent(U),\ent(V),\ent(W),\ent(Q)\big)
\]
is the entropy vector of the extension $\Ebf$.

The set of all possible values of the entropy vector for all extensions of
$(X,Y,U)$
\begin{align*}
&\res(X,Y,U)
:=
\\
&\set{\ent_{*}(\Ebf)\in\Rbb^{7}\st\text{$\Ebf$ is an extension of $(X,Y,U)$}}
\subset
\Rbb^{7}
\end{align*}
is what we call the unstabilized relative entropic set of the
two-fan $(X\ot U\to Y)$. 

\subsubsection{The unstabilized relative entropic sets for the examples}
\label{s:unstabilizedexamples}
It turns out that these \emph{unstabilized} relative entropic sets of $(X_1 \ot U_1 \to Y_1)$ and $(X_2 \ot U_2 \to Y_2)$ are different
\[
\res(X_{1},Y_{1},U_{1})\neq\res(X_{2},Y_{2},U_{2})
\]

To see this, let us calculate some particular points in the unstabilized relative entropic sets of the Examples 1--3. We consider the constrained Information-Optimization problem, of finding an extension $\Ebf = (X, Y, Z, U, V, W, Q)$ of $(X, Y, U)$ such that
\begin{enumerate}
\item \label{i:low}
the space $Z$ is a reduction of $U$, that is
\[
\Ent(Q) = \Ent(U)
\]
\item \label{i:indep}
the spaces $X$ and $Y$ are independent conditioned on $Z$, 
\[
\Ent(Q) + \Ent(Z) = \Ent(V) + \Ent(W)
\]
\item\label{i:max}
the sum
\[
\Ent(X \rel Z) + \Ent(Y \rel Z)
\]
is maximal, subject to conditions (\ref{i:low}) and (\ref{i:indep}).
\end{enumerate}

  It is very easy to read the solutions $\hat\Ebf_{1}$, $\hat\Ebf_{2}$
  and $\hat\Ebf_{3}$ of this optimization problem for Examples 1, 2
  and 3 right from the pictures in Figure~\ref{f:squares-snakes}.
  Indeed, condition~(\ref{i:low}) says that $Z_i$ must be a partition
  of $U_i$. Condition~(\ref{i:indep}) says that each set in the
  partition must be ``rectangular'', that is it must be a Cartesian
  product of a subset of $\un X_{i}$ and a subset of $\un Y_{i}$. The
  quantity to be maximized is the average $\log$-area of the sets in
  the partition.

  Optima are very easy to find ``by hand''. For Example 1 it is a
  partition of $U_{1}$ into three $2\times2$ squares. In
  Examples 2 and 3, one of the solutions is the partition of $U_{2}$,
  (resp. $U_{3}$) into $1\times2$ rectangles. Thus, the optimal values
  are $(2\ln2)$ for Example 1, and $(\ln2)$ for Examples
  2 and 3.

\subsubsection{The stabilized relative entropic set}
  We have just seen that Examples 1 and 2 can be told apart by
  determining the unstabilized relative entropic set.  However, this
  is not really what we are interested in.  Rather, we wonder whether
  high tensor powers can be distinguished this way.

  This relates to why we used the adjective ``unstabilized'': because
  the relative entropic set usually grows (is not stable) under taking
  tensor powers. That is, for every $n,k \in \Nbb$ it holds that
  \[\tageq{inclusionrelset}
    k\cdot \res (X^{\otimes n}, Y^{\otimes n}, U^{\otimes n})
    \subset 
    \res (X^{\otimes (k\cdot n)}, Y^{\otimes (k\cdot n)},
    U^{\otimes (k\cdot n)})
  \]
  but in general the set on the right-hand side can be strictly larger
  than the set on the left.

  In view of the inclusion (\ref{eq:inclusionrelset}) we may define the
  \emph{stabilized} relative entropic set
  \[
  \sres(X, Y, U) 
  =
  \closure\left( \lim_{n \to \infty}
  \frac{1}{n}\res(X^{\otimes n}, Y^{\otimes n}, U^{\otimes
    n})\right)
  \]
This set turns out to be convex. 

\subsubsection{The stabilized relative entropic set for the examples}

  In fact, the \emph{stabilized} relative entropic set also
  differentiates between Examples 1 and 2
  \[
    \sres(X_{1},Y_{1},U_{1})\neq\sres(X_{2},Y_{2},U_{2})
  \]
  The proof of this fact follows the same lines as in Section
  \ref{s:unstabilizedexamples}, but the stabilization makes the argument
  much more technical.

We expect that the stabilized relative entropic set cannot
differentiate between Examples 2 and 3. However, there are other types of relative
entropic sets, and other Information-Optimization problems that
\emph{can} differentiate between Examples 2 and 3.

The relative entropic sets are discussed in Section \ref{s:extensions}.

\subsection{Information-Optimization problems and relative entropic sets}

In Section \ref{s:unstabilizedexamples} we used an
Information-Optimization problem to find particular points in the
(unstabilized) relative entropic set.  This is no coincidence, and the
link between stabilized Information-Optimization problems and the
stable relative entropic set can be made very explicit.  Because the
stable relative entropic set is convex, it can be completely
characterized by Information-Optimization problems and vice versa.

Such Information-Optimization problems play a very important role in information theory \cite{Yeung-First-2012}, causal inference \cite{Steudel-Information-2015}, artificial intelligence
\cite{Dijk-Informational-2013}, information decomposition \cite{Bertschinger-Quantifying-2014}, robotics \cite{Ay-Predictive-2008},
and neuroscience \cite{Friston-Free-2009}.
The techniques developed in the article allow one to address this type of problems easily and efficiently.

\subsection{The intrinsic Kolmogorov-Sinai distance}

\skippar In \cite{Gromov-Search-2012}, Gromov proposed this as a basis of an asymptotic equivalence relation on 
sequences of probability spaces. 
Even though we were greatly influenced by ideas in \cite{Gromov-Search-2012}, we found that Gromov's definition does not extend easily to configurations of probability spaces. 
Formalizing and studying a notion of asymptotic equivalence for configurations of probability spaces is the main topic of the present article.

As we mentioned at the end of Section \ref{s:intro-single}, one is
tempted to define asymptotically equivalent configurations along the
lines of Shannon's source coding theorem following
\cite{Gromov-Search-2012}.  Two configurations would be asymptotically
equivalent if there is an almost measure-preserving bijection between
subspaces of almost full measure in their high tensor powers.

However, we found this approach inconvenient.  Instead of finding an
almost measure-preserving bijection between large parts of the two spaces,
 we consider a stochastic coupling
(transportation plan, joint distribution) between a pair of spaces and
measure its deviation from being an isomorphism of probability spaces,
that is a measure-preserving bijection.  Such a measure of deviation
from being an isomorphism then leads to the notion of intrinsic
Kolmogorov-Sinai distance, and its stable version -- the asymptotic
Kolmogorov-Sinai distance, as explained in Section \ref{s:kolmogorov}.

In the case of single probability spaces we define the \emph{intrinsic
  Kolmogorov-Sinai distance} between two probability spaces $X=(\un
X,p_{X})$ and $Y=(\un Y,p_{Y})$ by
\[
\ikd(X, Y) := \inf \set{\big[\ent(Z)-\ent(X)\big] +
\big[\ent(Z)-\ent(Y)\big]}
\]
where the infimum is taken over all choices of the joint distribution
$Z=(\un X\times\un Y,p_{Z})$. Note that each of the summands is
nonnegative and vanishes if and only if the corresponding
marginalization $Z\to X$ or $Z\to Y$ is an isomorphism of probability
spaces.  In this sense the distance measures the deviation from the
existence of a measure-preserving bijection between $X$ and $Y$.

Furthermore, we define the \emph{asymptotic Kolmogorov-Sinai distance}
between two probability spaces $X$ and $Y$ by
\[
\aikd(X, Y) = \lim_{n \to \infty} \frac{1}{n} \ikd (X^{\otimes n }, Y^{\otimes n}).
\]

This definition could be generalized to configurations of probability
spaces and we will say that two configurations are asymptotically
equivalent if the asymptotic Kolmogorov-Sinai distance between them
vanishes.

\subsection{Asymptotic Equipartition Property}\label{s:intro-aep}  
Examples 1, 2, and 3 above have the property that the symmetry group
acts transitively on the support of the measure on $U_{i}$ and they
are particular instances of what we call homogeneous configurations.

In Section \ref{s:ac-aep-conf}, we show an \emph{Asymptotic
  Equipartion Property for configurations}: Theorem
\ref{p:aep-complete} states that every sequence of tensor powers of a
configuration can be approximated in the asymptotic Kolmogorov-Sinai
distance by a sequence of homogeneous configurations.

This Asymptotic Equipartition Property allows one to substitute
configurations of probability spaces by homogeneous approximations.
Homogeneous probability spaces are just uniform probability spaces,
and as a first simple consequence of the Asymptotic Equipartition
Property, the asymptotic Kolmogorov-Sinai distance between probability
spaces $X$ and $Y$ can be computed and equals
\[
\aikd(X , Y) = |\ent(X) - \ent(Y)|.
\]

Homogeneous \emph{configurations} are, unlike homogeneous probability
spaces, rather complex objects.  Nonetheless, they seem to be simpler
than arbitrary configurations of probability spaces for the types of
problems that we would like to address.

More specifically, we show in Section \ref{s:extensions} that the
optimal values in (stabilized) Information-Optimization problems only
depend on the asymptotic class of a configuration and that they are
continuous with respect to the asymptotic Kolmogorov-Sinai distance;
in many cases, the optimizers are continuous as well.  The Asymptotic
Equipartition Property implies that for the purposes of calculating
optimal values and approximate optimizers, one only needs to consider
homogeneous configurations and this can greatly simplify computations.

Summarizing, the Asymptotic Equipartition Property and the continuity
of Information-Optimization problems are important justifications for
the choice of asymptotic equivalence relation and the introduction of
the intrinsic and asymptotic Kolmogorov-Sinai distances.

\subsection{The article}
The article has the following structure. Section~\ref{s:category} is
devoted to the basic setup used throughout the text.  In Section
\ref{s:config} we explain what we mean by configurations of
probability spaces, give examples, describe simple properties and
operations. Further, in Section~\ref{s:disttypes} we generalize the
notion of probability distribution to that on configurations and
discuss the theory of types for configurations.  In
Section~\ref{s:kolmogorov} the intrinsic Kolmogorov-Sinai distance and
the asymptotic Kolmogorov-Sinai distance are introduced and some
technical tools for the estimation of Kolmogorov distance are
developed.  Section~\ref{s:lagging} contains estimates on the
distances between types.  We use these estimates in the proof of the
Asymptotic Equipartition Property for configurations in
Section~\ref{s:ac-aep-conf}.  Section \ref{s:extensions} deals with
extensions of configurations. We prove there the Extension Lemma,
which is used to show continuity of extensions and implies, in
particular, that solutions of the constrained optimization problem for
the entropies of extensions are Lipschitz-continuous with respect to
the asymptotic Kolmogorov-Sinai distance, thus they only depend on the
asymptotic classes of configurations.  In Section~\ref{s:mixtures} we
briefly discuss a special type of configurations called mixtures,
which will play an important role in the construction of tropical
probability spaces.  Finally, in Section~\ref{s:tropical} we introduce
the notion of tropical probability spaces and configurations thereof,
and list some of their properties. We will continue our study of
tropical probability spaces and configurations in subsequent articles.

Some technical and not very illuminating proofs are deferred to
Section~\ref{s:technical} Technical Proofs.  In the electronic version
one can move between the proof in the technical section and the
statement in the main text by following the link (arrow up or down).

	\bigskip
	
	\paragraph{\bf Acknowledgments.} 
	We would like to thank Tobias Fritz, Franti{\v s}ek Mat{\' u}{\v s}, Misha Movshev and Johannes Rauh for inspiring discussions. We are grateful to the participants of the Wednesday Morning Session at the CASA group at the Eindhoven University of Technology for valuable feedback on the introduction of the article.
Finally, we thank the Max Planck Institute for Mathematics in the Sciences, Leipzig, for its hospitality.
	
	\section[Category]{Category of probability spaces and configurations}
	\label{s:category}
	  This section is devoted to the basic setup used throughout the
  present article.  We introduce a category of probability spaces and
  reductions, similar to categories introduced in \cite{Baez-Characterization-2011} and
  \cite{Gromov-Search-2012}, and define configurations of probability
  spaces and the corresponding category.  The last subsection recalls
  the notion of entropy and its elementary properties.

\subsection{Probability spaces and reductions}\label{s:category-prob}
  Below we will consider probability spaces such that the support of
  the probability measure is finite. Any such space contains a
  full-measure subspace isomorphic to a finite space, thus we call
  such objects \term[finite probability space]{finite probability
  spaces}.  For a probability space $X=(S,p_{X})$ denote by $\un
  X=\supp p_{X}$ the support of the measure and by $|X|$ its
  cardinality. Slightly abusing the language, we call this quantity
  the \term[cardinality of probability space]{cardinality} of $X$.

  For a pair of probability spaces a \term{reduction} $X\to Y$ is a
  class of measure-preserving maps, with two maps being equivalent if
  they coincide on a set of full measure. The composition of two
  reductions is itself a reduction.  Two probability spaces are
  isomorphic if there is a measure-preserving bijection between
  the supports of the probability measures. Such a bijection defines an
  invertible reduction from one space into another. Clearly the
  cardinality $|X|$ is an isomorphism invariant.
  The automorphism group $\Aut(X)$ is the group of all
  self-isomorphisms of $X$.

  A probability space $X$ is called \term[homogeneous probability
  space]{homogeneous} if the automorphism group $\Aut(X)$ acts
  transitively on the support $\un X$ of the measure.  The property of
  being homogeneous is an isomorphism invariant.  In the isomorphism
  class of a homogeneous space there is a representative with uniform
  measure.

  The finite probability spaces and reductions form a category, that
  we denote by $\prob$. The subcategory of homogeneous spaces will be
  denoted by $\probhom$.  The isomorphism in the category coincides
  with the notion of isomorphism above.

  The category $\prob$ is not a small category. However it has a small
  full subcategory, that contains an object for every isomorphism class in
  $\prob$ and for every pair of objects in it, it contains all the
  available morphisms between them. From now on we imagine that such a
  subcategory was chosen and fixed and replaces $\prob$ in all
  considerations below.

  There is a product in $\prob$ given by the Cartesian product of
  probability spaces, that we will denote by $X\otimes Y:=(\un
  X\times\un Y,p_{X}\otimes p_{Y})$. There are canonical reductions
  $X\otimes Y\to X$ and $X\otimes Y\to Y$ given by projections to
  factors. For a pair of reductions $f_{i}:X_{i}\to Y_{i}$, $i=1,2$
  their tensor product is the reduction $f_{1}\otimes
  f_{2}:X_{1}\otimes X_{2}\to Y_{1}\otimes Y_{2}$, which is equal to
  the class of the Cartesian product of maps representing $f_{i}$'s.
   The tensor product is not however a categorical product. The product
  leaves the subcategory of homogeneous spaces invariant.

  The probability measure on $X$ will be usually denoted by $p_{X}$ or
  simply $p$, when the risk of confusion is low.

\subsection{Configurations of probability spaces}
\label{s:category-config}
  Essentially, a configuration $\Xcal=\set{X_{i};f_{ij}}$ is a
  commutative diagram consisting of a finite number of probability spaces and
  reductions between some of them, that is transitively closed, while
  a morphism $\rho:\Xcal\to\Ycal$ between two configurations
  $\Xcal=\set{X_{i};f_{ij}}$ and $\Ycal=\set{Y_{i};g_{ij}}$ of the
  same combinatorial type is a collection of reductions between
  corresponding individual objects $\rho_{i}:X_{i}\to Y_{i}$, that
  commute with the reductions within each configuration,
  $\rho_{j}\circ f_{ij}=g_{ij}\circ\rho_{i}$.
  
  We need to keep track of the combinatorial structure of the
  collection of reductions within a configuration. There are several
  possibilities for doing so: 
  \begin{itemize}
  \item 
    the reductions form a directed transitively closed graph without
    loops;
  \item 
    the spaces in the configuration form a poset; 
  \item
    the underlying combinatorial structure could be recorded as a finite
    category. 
  \end{itemize}

  The last option seems to be most convenient since it has many
  operations necessary for our analysis already built-in.

  A \term{diagram category} $\Gbf$ is a finite category such that
  for each pair of objects $O_{1}$, $O_{2}$ in $\Gbf$ the morphism
  space between them
  \[
    \Hom_{\Gbf}(O_{1},O_{2})\cup\Hom_{\Gbf}(O_{2},O_{1})
  \] 
  contains at most one element.

  For a diagram category $\Gbf$ a \term{configuration of probability
    spaces} modeled on $\Gbf$ is a functor $\Xcal: \Gbf\to \prob$.
  The collection of all configurations of probability spaces modeled
  on a fixed diagram category $\Gbf$ forms the category of functors
  $\prob\<\Gbf>:=[\Gbf,\prob]$. The objects of $\prob\<\Gbf>$ are
  configurations, that is functors from $\Gbf$ to $\prob$, while
  morphisms in $\prob\<\Gbf>$ are natural transformations between
  them.  For a configuration $\Xcal\in\prob\<\Gbf>$, the diagram
  category $\Gbf$ will be called the \term[combinatorial type of a
    configuration]{combinatorial type} of $\Xcal$.

  For a diagram category $\Gbf$ or a configuration
  $\Xcal\in\prob\<\Gbf>$ we denote by $\size{\Gbf}=\size{\Xcal}$ the
  number of objects in the category $\Gbf$.

  An object $O$ in a diagram category $\Gbf$ will be called
  \term[initial object/space]{initial}, if it is not a target of any
  morphism except for the identity. Likewise a \term[terminal
    object/space]{terminal} object is not a source of any morphism,
  except for the identity morphism. Note that this terminology is
  somewhat unconventional from the point of view of category theory.
  
   A diagram category is called \term[complete diagram
     category]{complete} if it has a unique initial object. Thus a
   configuration modeled on a complete category includes a space that
   reduces to all other spaces in the configuration.

  The above terminology transfers to configurations modeled on
  $\Gbf$: An initial space in $\Xcal\in\prob\<\Gbf>$ is one that
  is not a target space of any reduction within the configuration, a
  terminal space is not a source of any non-trivial reduction and
  $\Xcal$ is complete if $\Gbf$ is, that is there is a unique initial
  space.

  The tensor product of probability spaces extends to a tensor product
  of configurations. For $\Xcal,\Ycal\in\prob\<\Gbf>$, such that
  $\Xcal=\set{X_{i};f_{ij}}$ and $\Ycal=\set{Y_{i};g_{ij}}$ define
  \[
    \Xcal\otimes\Ycal
    :=
    \set{X_{i}\otimes Y_{i};f_{ij}\otimes g_{ij}}
  \] 

  Occasionally we will also talk about configuration of sets. Denote by
  $\Set$ the category of finite sets and surjective maps. Then all of the
  above constructions could be repeated for sets instead of probability
  spaces. Thus we could talk about the category of configurations of sets
  $\Set\<\Gbf>$.

  Given a reduction $f:X\to Y$ between two probability spaces, the
  restriction $\un f:\un X\to\un Y$ is a well-defined surjective map.
  Given a configuration $\Xcal=\set{X_{i};f_{ij}}$ of probability
  spaces, there is an underlying configuration of sets, obtained by
  taking the supports of measures on each level and restricting
  reductions on these supports. We will denote it by
  $\un\Xcal=\set{\un X_{i};f_{ij}}$, where $\un X_{i}:=\supp
  p_{X_{i}}$. Thus we have a forgetful functor
  \[
  \underline{\mkern5mu\cdot\mkern5mu}:\prob\<\Gbf>\to\Set\<\Gbf>
  \]
  
  For now we will consider two important examples of diagram
  categories and configurations modeled on them. We give further
  examples in Section \ref{s:config-examples}.

\subsubsection{Two-fans:}\label{s:category-config-2fan}
  A two-fan is a configuration modeled on the category $\Lambdabf$
  with three objects, one initial and two terminal.
  \[
    \Lambdabf=(O_{1}\ot O_{12}\to O_{2})
  \]
  
  There is a special significance to two-fans, since these are the
  simplest non-trivial configurations, as we will see later.

  Essentially, a two-fan $X\ot Z\to Y$ is a triple of probability spaces and a pair of
  reductions  between them. 

  A reduction of a two-fan $X\ot Z\to Y$ to another two-fan $X'\ot Z'\to Y'$
  is a triple of reductions $Z\to Z'$, $Y\to Y'$ and $X\to X'$ that
  commute with the reductions within each fan, that is, the following diagram is commutative
  \[
  \begin{tikzcd}
  X 
  \arrow{d}
  & Z
  \arrow{l}
  \arrow{d}
  \arrow{r}
  & Y 
  \arrow{d}\\
  X' 
  & 
  Z' 
  \arrow{l}
  \arrow{r}
  & Y'  
  \end{tikzcd}
  \] 
  
  Isomorphisms and the automorphism group $\Aut(\cdot)$
  are defined accordingly.  Note that terminal spaces in a two-fan are
  labeled and reductions preserve the labeling.
  
  A two-fan $X\ot Z\to Y$ is called \term[minimal 2-fan]{minimal} if
  for a.e.~$x\in X$ and $y\in Y$ there is a unique $z\in Z$, that
  reduces to $y$ and to $x$.  Given a two-fan $X\ot Z\to Y$, there is
  always a reduction to a minimal two-fan $X\ot Z'\to Y$. Such minimal
  reduction is unique up to isomorphism. Explicitly, take $Z':=\un
  X\times \un Y$ as a set and consider a probability distribution on
  $Z'$ induced by a map $Z\to Z'$ which is the Cartesian product of
  the reductions $X\ot Z\to Y$ in the original two-fan.

  The notion of being minimal is in fact a categorical notion. It could
  be equivalently defined by saying that a two-fan $\Xcal=(X_{1}\ot X_{12}\to
  X_{2})$ is minimal if for any reduction $\lambda:\Xcal\to\Xcal'$
  holds: if both $\lambda_{1}$ and $\lambda_{2}$ are isomorphisms, then
  $\lambda_{12}$ is also an isomorphism. Consequently, if one specifies
  reductions from the terminal spaces of a minimal two-fan to another
  two-fan, then there exists at most one extension to the reduction of the
  whole fan.

  The inclusion of a pair of probability
  spaces $X$ and $Y$ as terminal vertices in a minimal two-fan is
  equivalent to specifying a joint distribution on $\un X\times \un Y$.

  An arbitrary configuration $\Xcal$ will be called \term[minimal
    configuration]{minimal} if with every two-fan, it also contains a
  minimal two-fan with the same terminal spaces.
  We will denote the space of minimal configurations modelled on a
    diagram category $\Gbf$ by $\prob\<\Gbf>_{\mbf}$.

  Given a two-fan 
  \[
    \Fcal=(X\ot Z\to Y)
  \]
  with terminal spaces $X$ and $Y$, and a point $x \in X$ with $p_X(x) > 0$, one may construct a \term{conditional probability distribution}
  $p_Y(\;\cdot\;\rel x)$ on $\un Y$. We denote the corresponding space
  $Y\rel x:=\big(\un Y,p_Y(\;\cdot\;\rel x)\big)$. The usual bar
  ``$|$'', that is normally used for conditioning, interferes with our
  notations for cardinality of spaces. 
  We will give more details in Section \ref{s:config-conditioning}.

\subsubsection{A diamond configuration}
\label{s:category-config-diamond}
  A ``diamond'' configuration is modeled on a diamond category $\dia$,
  that consists of a two-fan and a ``co-fan'':
   \[
   \dia =
  \left(  \begin{tikzcd}[row sep=small, column sep=small]
  \mbox{}
  & 
  O_{12}
  \arrow{dl}{}
  \arrow{dr}{}
  &
  \\
  O_{1}
  \arrow{dr}{}
  & 
  &
  O_{2}
  \arrow{dl}{}
  \\
  \mbox{}
  &
  O_{\bullet}
  &
  \end{tikzcd} \right)
  \]
  
   Of course, there is also a morphism $O_{12}\to O_{\bullet}$, which
  lies in the transitive closure of the given four morphisms. As a
  rule, we will skip writing morphisms, that are implied by the
  transitive closure.

  A diamond configuration is minimal if the top two-fan in it is
  minimal.

\subsection{Entropy}\label{s:category-entropy}
  Our working definition of \term{entropy} will be based on the
  following version of the asymptotic equipartition theorem for
  Bernoulli process, see \cite{Cover-Elements-1991}.
  \begin{theorem}{p:aepbernoullisingle}
    Suppose $X$ is a finite probability space, then for any $\epsilon>0$ and
    any $n>\!>0$ there exists a subset $A^{\c n}_{\epsilon}\subset
    X^{\otimes n}$ such that
    \begin{enumerate}
    \item 
      $p(A^{\c n}_{\epsilon})\geq1-\epsilon$
    \item
      For any $a,a'\in A^{\c n}_{\epsilon}$ holds
      \[
        \left|\frac{\ln p(a)}{n}-\frac{\ln p(a')}{n}\right|\leq\epsilon
      \]
    \end{enumerate}
    Moreover, if $A^{\c n}_{\epsilon}$ and $B^{\c n}_{\epsilon}$ are
    two subsets of $X^{\otimes n}$ satisfying two conditions above, then their
    cardinalities satisfy
    \[\tageq{cardinalityrange}
      \left|
        \frac{\ln |A^{\c n}_{\epsilon}|}{n}
        -
        \frac{\ln |B^{\c n}_{\epsilon}|}{n}
      \right|
      \leq 2\epsilon
    \]
  \end{theorem}
 
  Then we define
  \[
    \ent(X):=\lim_{\epsilon\to0}\lim_{n\to\infty}\frac1n\ln|A^{\c n}_{\epsilon}|
  \]

  Clearly, in view of the property (\ref{eq:cardinalityrange}) in the
  Theorem \ref{p:aepbernoullisingle} the limit above is well-defined
  and is independent of the choice of the typical subsets $A_\epsilon^{(n)}$.

  Entropy satisfies the so-called Shannon inequality,
  see for example \cite{Cover-Elements-1991}, namely for any minimal
  diamond configuration
  
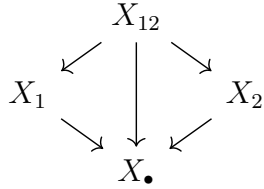
\[
  \begin{tikzcd}[row sep=small, column sep=small]
    \mbox{}
    & 
    X_{12}
    \arrow{dl}{}
    \arrow{dr}{}
    \arrow{dd}{}
    &
    \\
    X_{1}
    \arrow{dr}{}
    & 
    &
    X_{2}
    \arrow{dl}{}
    \\
    \mbox{}
    &
    X_{\bullet}
    &
  \end{tikzcd}
  \]
  the following inequality holds, 
  \[\tageq{shannonineq}
    \ent(X_{1})+\ent(X_{2})
    \geq
    \ent(X_{12})+\ent(X_{\bullet})
  \]

  Furthermore, entropy is additive with respect to the tensor product,
  that is, for a pair of probability spaces $X,Y\in\prob$ holds
  \[\tageq{entropyadditive}
  \ent(X\otimes Y)=\ent(X)+\ent(Y)
  \]

  Further, for a pair $X$, $Y$ of probability spaces included in a
  minimal two-fan $(X\ot Z\to Y)$ we define the conditional entropy
  \[
    \ent(X\rel Y):=\ent(Z)-\ent(X)
  \]

  The above quantity is always non-negative in view of Shannon
  inequality~(\ref{eq:shannonineq}).  Moreover, the following
  identity holds, see \cite{Cover-Elements-1991}
  \[\tageq{relentinteg}
    \ent(X\rel Y)
    =
    \int_{Y}\ent(X\rel y)\d p_{Y}(y)
  \]

	\section{Configurations}\label{s:config}
	  In this section we will look at configurations in more detail. We
  start by considering some important examples.

\subsection[Examples]{Examples of configurations}\label{s:config-examples}

\subsubsection{Singleton}\label{s:config-examples-single}
  We denote by $\bullet$ a diagram category with a single object.
  Clearly configurations modeled on $\bullet$ are just probability
  spaces and we have $\prob\equiv\prob\<\bullet>$.

\subsubsection{Chains}\label{s:config-examples-chain}
  The chain $\Cbf_{n}$ of length $n\in\Nbb$ is a category with $n$ objects
  $\set{O_{i}}_{i=1}^{n}$ and morphisms from $O_{i}$ to $O_{j}$
  whenever $i\leq j$. A configuration $\Xcal\in\prob\<\Cbf_{n}>$ is a
  chain of reductions
  \[
    \Xcal=(X_{1}\to X_{2}\to\cdots\to X_{n})
  \]

\subsubsection{Two-fan}\label{s:config-examples-fan}
  The two-fan $\Lambdabf$ is a category with three objects
  $\set{O_{1},O_{12},O_{2}}$ and two non-identity morphisms $O_{12}\to
  O_{1}$ and $O_{12}\to O_{2}$.  See also
  Section~\ref{s:category-config-2fan}.  Two-fans are the simplest
  configurations for which asymptotic equivalence classes contain more
  information than just entropies of the entries.
  
  Recall that a fan $(X\ot Z\to Y)$ is called minimal, if for any pair
  of points $x\in X$ and $y\in Y$ with positive weights there exists at
  most one $z\in Z$, that reduces to $x$ and to $y$. Equivalently, for
  any super-configuration
  
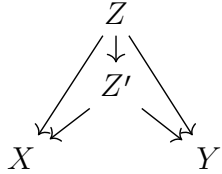
\[
  \begin{tikzcd}[column sep=small,row sep=small]
    \mbox{}
    &
    Z
    \arrow{ddl}{}
    \arrow{ddr}{}
    \arrow{d}{}
    &
    \mbox{}
    \\
    \mbox{}
    &
    Z'
    \arrow{dl}{}
    \arrow{dr}{}
    &
    \mbox{}
    \\
    X
    &
    &
    Y
  \end{tikzcd}
  \]
  the reduction $Z\to Z'$ must be an isomorphism.

\subsubsection{Full configuration}\label{s:config-examples-full}
  The full category $\Lambdabf_{n}$ on $n$ objects is a category with
  objects $\set{O_{I}}_{I\in
    2^{\set{1,\ldots,n}}\setminus\set{\emptyset}}$ indexed by all
  non-empty subsets $I\in 2^{\set{1,\ldots,n}}$ and a morphism from
  $O_{I}$ to $O_{J}$, whenever $J\subseteq I$.

  For a collection of random variables $X_{1},\ldots,X_{n}$ one may
  construct a minimal full configuration $\Xcal\in\prob\<\Lambdabf_{n}>$ by
  considering all joint distributions and ``marginalization''
  reductions. We denote such a configuration by
  $\<X_{1},\ldots,X_{n}>$. 
  On the other hand, the terminal
  vertices of a full configuration can be viewed as random variables
  on the domain of definition given by the (unique) initial space.
  
  Suppose $\Xcal\in\prob\<\Lambdabf_{n}>$ is a minimal full
  configuration with terminal vertices $X_{1},\ldots,X_{n}$. It is
  convenient to view $\Xcal$ as a distribution on the Cartesian
  product of the underlying sets of the terminal vertices:
  \[
    p_{\Xcal}\in\Delta(\un X_{1}\times\cdots\times\un X_{n})
  \]
  Once the underlying sets of the terminal spaces are fixed, there is
  a one-to-one correspondence between the full minimal configurations
  and distributions as above.

\subsubsection{``Two-tents'' configuration}
\label{s:config-examples-twotents}
  The ``two-tents'' category $\Mbf_{2}$ consists of five objects, of
  which two are initial and three are terminal, and morphisms are as
  follows
  \[
  \Mbf_{2}=
  \left(
  \begin{tikzcd}[column sep=small,ampersand replacement=\&]
  O_{12}
  \arrow[xshift=-0.3em]{d}{}
  \arrow{dr}{} 
  \&
  \&
  O_{23} 
  \arrow[xshift=-0.3em]{dl}{}
  \arrow{d}{}
  \\
  O_1 
  \&
  O_2
  \&
  O_3
  \end{tikzcd}
  \right)
  \]
  
  Thus, a typical two-tents configuration consists of five probability
  spaces and reduction as in
  \[
    \Xcal=(X\ot U\to Y\ot V\to Z)
  \]
  The probability spaces $U$ and $V$ are initial and $X$, $Y$ and $Z$
  are terminal.
 
\subsubsection{``Many-tents'' configuration}
\label{s:config-examples-manytents} 
  The previous example could be generalized to a ``many-tents'' category
  \[
    \Mbf_{n}
    =
    (O_{1}\ot O_{12}\to O_{2}\ot\cdots\to O_{n-1}\ot O_{n-1,n}\to O_{n})
  \]

\subsubsection{``Fence'' configuration}
\label{s:config-examples-fence} 
  The ``fence'' category $\Wbf_{3}$ consists of six objects and the
  morphisms are 
  \[
    \Wbf_{3}
    =
    \left(
    \begin{tikzcd}[ampersand replacement=\&,row sep=small, column sep=small]
      O_{12}
      \arrow{d}
      \&
      O_{13}
      \arrow{dr}
      \arrow{dl}
      \&
      O_{23}
      \arrow{d}
      \\
      O_{1}
      \&
      O_{2}
      \arrow[leftarrow, crossing over]{ur}
      \arrow[leftarrow, crossing over]{ul}
      \&
      O_{3}      
    \end{tikzcd}
    \right)
  \]

\subsubsection{Co-fan}\label{s:config-examples-cofan}
A co-fan $\Vbf$ is a category with three objects and morphisms as in the diagram
\[
\Vbf =
\left(  \begin{tikzcd}[row sep=small, column sep=small]
\mbox{}
O_{1}
\arrow{dr}{}
& 
&
O_{2}
\arrow{dl}{}
\\
\mbox{}
&
O_{\bullet}
&
\end{tikzcd} \right)
\]

\subsubsection{``Diamond'' configurations}
\label{s:config-examples-diamond} 
  A ``diamond'' configuration $\dia$ is modeled on a diamond category
  that consists of a fan and a co-fan
  \[
  \dia =
  \left(  
  \begin{tikzcd}[ampersand replacement=\&,row sep=small, column sep=small]
    \mbox{}
    \& 
    O_{12}
    \arrow{dl}{}
    \arrow{dr}{}
    \&
    \\
    O_{1}
    \arrow{dr}{}
    \& 
    \&
    O_{2}
    \arrow{dl}{}
    \\
    \mbox{}
    \&
    O_{\bullet}
    \&
  \end{tikzcd} \right)
  \]
  See also Section~\ref{s:category-config-diamond}.

  \bigskip

  Examples
  \ref{s:config-examples-single}, 
  \ref{s:config-examples-chain}, 
  \ref{s:config-examples-fan},
  \ref{s:config-examples-full} and 
  \ref{s:config-examples-diamond}
  are complete. Examples \ref{s:config-examples-single},
  \ref{s:config-examples-chain} and \ref{s:config-examples-cofan} do
  not contain a two-fan. Tropical limits of such configurations are
  very simple. Essentially, such tropical limits correspond to the
  tuple of numbers corresponding to the entropies of the constituent
  spaces. Therefore, we call configurations not containing a two-fan
  \term[simple configuration]{simple}.
  
\subsection{Constant configurations}
\label{s:config-constantconfig}
  Suppose $X$ is a probability space and $\Gbf$ is a diagram category.
  One may form a \term[constant configuration]{constant
    $\Gbf$-configuration} by considering a functor that maps all
  objects in $\Gbf$ to $X$ and all the morphisms to the identity
  morphism $X\to X$. We denote such a constant configuration by
  $X^{\Gbf}$ or simply by $X$, when $\Gbf$ is clear from the context.
  Any constant configuration is automatically minimal.

  If $\Ycal=\set{Y_{i};f_{ij}}$ is another $\Gbf$-configuration, then
  a reduction $\rho:\Ycal\to X^{\Gbf}$ (which we write sometimes
  simply as $\rho:\Ycal\to X$) is a collection of reductions
  $\rho_{i}:Y_{i}\to X$, such that
  \[
    f_{ij}\circ\rho_{i}=\rho_{j}
  \]
  
\subsection{Configurations of configurations}
\label{s:config-config}
  Of course, the operation of ``configuration'' could be iterated, so
  given a pair $\Gbf_{1}$, $\Gbf_{2}$ of diagram categories we could
  form a $\Gbf_{2}$-configuration of $\Gbf_{1}$-configurations, so we
  could speak, for example, about a two-fan of configurations of the same
  type.
  \[
    \Prob\<\Gbf_{1},\Gbf_{2}>
    :=
    \Prob\<\Gbf_{1}>\<\Gbf_2>
    =
    \prob\<\Gbf_{1}\square\Gbf_2>
  \]
  where $\Gbf_{1}\square\Gbf_2$ is the ``Cartesian product of
  graphs'' (as every diagram category could be considered as a
  transitively closed directed graph).  This operation is commutative,
  thus, for example, a two-fan of $\Gbf$-configurations is a
  $\Gbf$-configuration of two-fans.
  
  We will rarely need anything beyond a two-fan of configurations.

  A two-fan $\Fcal=(\Xcal\ot\Zcal\to\Ycal)$ of
  $\Gbf$-configurations is called minimal if in any extension of
  $\Fcal$ of the form
  
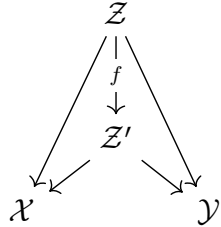
\[
  \begin{tikzcd}[ row sep=small, column sep=small,ampersand replacement=\&]
    \mbox{}
    \&
    \Zcal
    \arrow{dd}[description]{f}
    \arrow{dddl}{}
    \arrow{dddr}{}
    \&
    \mbox{}
    \\
    \mbox{}
    \&
    \&
    \\
    \mbox{}
    \&
    \Zcal'
    \arrow{dr}{}
    \arrow{dl}{}
    \&
    \mbox{}
    \\
    \Xcal
    \&
    \&
    \Ycal
  \end{tikzcd}
  \]
  the reduction $f:\Zcal\to\Zcal'$ must be an isomorphism of
  $\Gbf$-configurations.

  Recall that a two-fan of $\Gbf$-configurations could also be viewed
  as a $\Gbf$-configuration of two-fans of probability spaces. In the
  following lemma we show that in order to verify the minimality of a
  two-fan of configurations it is sufficient to check the minimality
  of all the constituent two-fans.
  \begin{lemma}{p:minimalfansconfig}
    Let $\Gbf$ be a diagram category. Then
    \begin{enumerate}
    \item\label{p:minimalfansconfig1} A two-fan
      $\Fcal=(\Xcal\ot\Zcal\to\Ycal)$ of $\Gbf$-configurations is
      minimal, if and only if the constituent two-fans of probability
      spaces $\Fcal_{i}=(X_{i}\ot Z_{i}\to Y_{i})$ are all minimal.
    \item \label{p:minimalfansconfig2}
      For any two-fan $\Fcal=(\Xcal\ot\Zcal\to\Ycal)$ of
      $\Gbf$-configurations its minimal reduction exists,
      that is, there exists a minimal two-fan
      $\Fcal'=(\Xcal\ot\Zcal'\to\Ycal)$ included in the following
      diagram
      
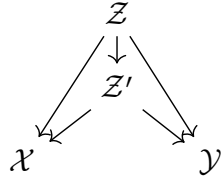
\[
      \begin{tikzcd}[row sep=small, column sep=small, ampersand replacement=\&]
        \mbox{}
        \&
        \Zcal
        \arrow{d}{}
        \arrow{ddl}{}
        \arrow{ddr}{}
        \&
        \mbox{}
        \\
        \mbox{}
        \&
        \Zcal'
        \arrow{dr}{}
        \arrow{dl}{}
        \&
        \mbox{}
        \\
        \Xcal
        \&
        \&
        \Ycal
      \end{tikzcd}
      \]      
    \end{enumerate}
  \end{lemma}

  Even though this lemma is rather elementary, there are many similar
  statements that are not true. Thus we are compelled to provide a
  proof, which can be found in Section~\ref{s:technical} on
  page~\pageref{p:minimalfansconfig.rep}.

  Similarly, a full configuration $\Fcal$ of $\Gbf$-configurations is
  called minimal if for every two-fan of $\Gbf$-configurations in
  $\Fcal$ there is a minimal two-fan in $\Fcal$ of
  $\Gbf$-configurations with the same terminal configurations.

  Lemma \ref{p:minimalfansconfig} has the following corollary and
  counterpart for full configurations of $\Gbf$-configurations.
  \begin{corollary}{p:minimalfullconfig}
    Let $\Gbf$ be a diagram category. Then
    \begin{enumerate}
    \item\label{p:minimalfullconfig1} 
      A full configuration $\Fcal$ of $\Gbf$-configurations is minimal,
      if and only if the constituent full configurations of probability
      spaces $\Fcal_{i}$ are all minimal.
    \item \label{p:minimalfullconfig2}   
      For any full configuration $\Fcal$ of $\Gbf$-configurations its
      minimal reduction exists.
    \end{enumerate}
  \end{corollary}
  
\subsection{Restrictions and extensions}
  \label{s:config-restrictions}
  Suppose $R:\Gbf_{1}\to \Gbf_{2}$ is a functor between two diagram
  categories. For a configuration $\Xcal:\Gbf_{2}\to \prob$, the pull-back
  configuration $\Ycal=R^{*}\Xcal\in\prob\<\Gbf_{1}>$ defined as the
  composition
  \[
    \Ycal
    :=
    \Xcal\circ R
  \]
  is called an \term[restriction functor]{$R$-restriction} of $\Xcal$
  to $\Gbf_1$ and $\Xcal$ is the \term[extension of a
    configuration]{extension} of $\Ycal$. If the functor $R$ is
  injective then we call $\Ycal=R^{*}\Xcal$ a \term{sub-configuration}
  of $\Xcal$ and write $\Ycal\subset\Xcal$, likewise $\Xcal$ will be
  called a \term{super-configuration} of $\Ycal$.

  The restriction operation is functorial in the sense that given
  two configurations 
  $\Xcal,\Xcal'\in\prob\<\Gbf_{1}>$ and a reduction
  $f:\Xcal\to\Xcal'$, there is a canonical reduction
  $R^{*}f:R^{*}\Xcal\to R^{*}\Xcal'$. Thus $R^{*}$ can be considered
  as a functor
  \[
    R^{*}:\prob\<\Gbf_{2}>\to\prob\<\Gbf_{1}>
  \]
  
  Some important examples of restrictions and extensions are below.

\subsubsection{Restriction of a full configuration to a smaller full
  configuration} 
\label{s:config-restrictions-fullfull}
  Recall that, as explained in Section
  \ref{s:config-examples-full}, the terminal vertices of a full
  configuration could be considered as random variables and any
  collection of random variables ``generates'' a full configuration.

  For a full configuration $\Xcal=\langle X_{i}\rangle_{i=1}^{n}$ and
  a subset $I\subset\set{1,\ldots,n}$ we denote by $R_{I}^*\Xcal=\langle
  X_{i}\rangle_{i\in I}$ the restriction of $\Xcal$ to a full
  configuration generated by $X_{i}$, $i\in I$.  We will also make use
  of the notation $R_{k,l}^*$ for the restriction operator $R_{\set{1,\ldots,k}}^*: \prob\<\Lambdabf_l> \to \prob\<\Lambdabf_k>$.

\subsubsection{Restriction of a $\Lambdabf_{3}$-configuration to
  an $\Mbf_{2}$-configuration}
\label{s:config-restrictions-full2tents}

  Given a full configuration $\Xcal\in\prob\<\Lambdabf_{3}>$ we may
  ``forget'' part of the data. If we, for example, forget the top
  space and the relation between a pair out of three terminal spaces
  we end up with the two-tents configuration. This operation
  corresponds to the inclusion functor
  \[
  M:
  \left(
  \begin{tikzcd}[column sep=small,row sep=normal,ampersand replacement=\&]
    O_{12}
    \arrow{d}{}
    \arrow{dr}{}
    \&
    \&
    O_{23}
    \arrow{d}{}
    \arrow{dl}{}
    \\
    O_{1}
    \&
    O_{2}
    \&
    O_{3}
  \end{tikzcd}
  \right)
  \longrightarrow
  \left(
  \begin{tikzcd}[column sep=small, row sep=normal,ampersand replacement=\&]
    \mbox{}
    \&
    Q_{123}
    \arrow{dl}{}
    \arrow{d}{}
    \arrow{dr}{}
    \\
    Q_{12}
    \arrow{d}{}
    \&
    Q_{31}
    \arrow{dl}{}
    \arrow{dr}{} 
    \&
    Q_{23}
    \arrow{d}
    \\
    Q_{1}
    \&
    Q_{2}
    \arrow[leftarrow, crossing over]{ul}
    \arrow[leftarrow, crossing over]{ur}
    \&
    Q_{3}
  \end{tikzcd}
  \right)
  \]
  that preserves the sub-indices.

We show in Section \ref{s:config-adhesion} below that the corresponding 
restriction operator
\[
M^{*}:\prob\<\Lambdabf_{3}>\to\prob\<\Mbf_{3}> 
\]
is surjective, both on objects and all morphisms. Thus, as a map of
collections of objects it has a right inverse. However, no natural right
inverse exists. 

\subsubsection{Restriction of a $\Lambdabf_{3}$-configuration to
  a $\Wbf_{3}$-configuration} 
\label{s:config-restrictions-fullfence}

Starting with a full configuration we might choose to forget the
initial space (and reductions, for which it was the domain). The
remaining configuration has the combinatorial type of a fence. This
operation corresponds to the functor
\[
W:
\left(
\begin{tikzcd}[column sep=small, row sep=normal,ampersand replacement=\&]
  O_{12}
  \arrow{d}{}
  \&
  O_{31}
  \arrow{dl}{}
  \arrow{dr}{} 
  \&
  O_{23}
  \arrow{d}{}
  \\
  O_{1}
  \&
  O_{2}
  \arrow[leftarrow, crossing over]{ul}
  \arrow[leftarrow, crossing over]{ur}
  \&
  O_{3}
\end{tikzcd}
\right)
\longrightarrow
\left(
\begin{tikzcd}[column sep=small, row sep=normal, ampersand replacement=\&]
  \&
  Q_{123}
  \arrow{dl}{}
  \arrow{d}{}
  \arrow{dr}{}
  \\
  Q_{12}
  \arrow{d}{}
  \&
  Q_{31}
  \arrow{dl}{}
  \arrow{dr}{} 
  \&
  Q_{23}
  \arrow{d}
  \\
  Q_{1}
  \&
  Q_{2}
  \arrow[leftarrow, crossing over]{ul}
  \arrow[leftarrow, crossing over]{ur}
  \&
  Q_{3}
\end{tikzcd}
\right)
\]
The corresponding operator 
\[
W^{*}:\prob\<\Lambdabf_{3}>\to\prob\<\Wbf_{3}>
\]
is not surjective. To find out when a $\Wbf_{3}$-configuration is extendable
to a $\Lambdabf_{3}$-configuration is an interesting problem, see for
example, \cite{Abramsky-Contextuality-2015} and references
therein.
It is our hope that the methods developed in this article might be
useful to address these questions.

\subsubsection{Doubling}
\label{s:config-restrictions-doubling}
 This will be the first example of an interesting
    functor between diagram categories, which is not injective. In
    this situation the term ``restriction'' does not really reflect
    the operation of pull-back well, however we did not come up with a
    better terminology.
 
  The doubling operation is the restriction of a two-fan to a two-tents
  configuration.
  Consider the two-fan category $\Lambdabf=(O_{1}\ot
  O_{12}\to O_{2})$ and a two-tents category $\Mbf_{2}=(Q_{1}\ot
  Q_{12}\to Q_{2}\ot Q_{23}\to Q_{3})$. Define the functor
  $D:\Mbf_{2}\to\Lambdabf$ by setting
  \[
  D:
  \begin{cases}
    Q_{1}\mapsto O_{1}\\
    Q_{2}\mapsto O_{2}\\
    Q_{3}\mapsto O_{1}\\
    Q_{12}\mapsto O_{12}\\
    Q_{23}\mapsto O_{12}\\
  \end{cases}
  \]

  Note that $D$ extends uniquely to the spaces of morphisms, since each
  morphism space is either empty or a one-point set.
  
  Thus $D^{*}(X\ot Z\to Y)=(X\ot Z\to Y\ot Z\to X)$, where the ``left''
  and ``right'' two-fans are isomorphic.
  
  This operation along with a particular Information-Optimization
    problem is related to the so-called copy operation, that was
  used to find many non-Shannon information inequalities, as
  described, for example, in \cite{Dougherty-Non-Shannon-2011}.

\subsection{Adhesion}\label{s:config-adhesion}
  Given a minimal two-tents configuration $\Xcal=(X\ot U\to Y\ot V\to
  Z)\in\prob\<\Mbf_{2}>_{\mbf}$ one could always construct an extension of
  $\Xcal$ to a full configuration $\ad(\Xcal)\in\prob\<\Lambdabf_{3}>_{\mbf}$ in the
  following way: As explained in Section \ref{s:config-examples-full},
  to construct a minimal full configuration with terminal vertices
  $X$, $Y$ and $Z$ it is sufficient to provide a distribution on $\un
  X\times\un Y\times\un Z$ with the correct marginals. We do this
  by setting
  \[
  p(x,y,z):=\frac{p_{U}(x,y)\cdot p_{V}(y,z)}{p_{Y}(y)}
  \]
  It is straightforward to check that the appropriate restriction of
  the full configuration defined in the above manner is indeed the
  original two-tents configuration.  Essentially, to extend we need to
  provide a relationship (coupling) between spaces $X$ and $Z$ and we
  do it by declaring $X$ and $Z$ independent relative to $Y$.  This is
  an instance of operation called \term{adhesion}, see
  \cite{Matus-Infinitely-2007}.

  If we call the top vertex in the full configuration $W$, the
  entropies achieve equality in the Shannon inequality, that is
  \[
    \Ent( U ) + \Ent( V ) - \Ent( W ) - \Ent(Y) = 0.
  \]
 
  Adhesion provides a right inverse $\ad$ to the restriction
  functor $M^{*}$ described in
  Section~\ref{s:config-restrictions-full2tents}
  \[
  \begin{tikzcd}[ampersand replacement=\&]
    \prob\<\Lambdabf_{3}>_{\mbf}
    \arrow{r}{M^{*}}
    \&
    \prob\<\Mbf_{2}>_{\mbf}
    \arrow[dashed,bend left]{l}{\ad}
  \end{tikzcd}
  \]

  It is important to note though, that the map $\ad$ is not functorial
  and, in fact, no functorial inverse of $M^{*}$ exists.

\subsection{Homogeneous configurations}\label{s:config-homo}
\begin{figure}
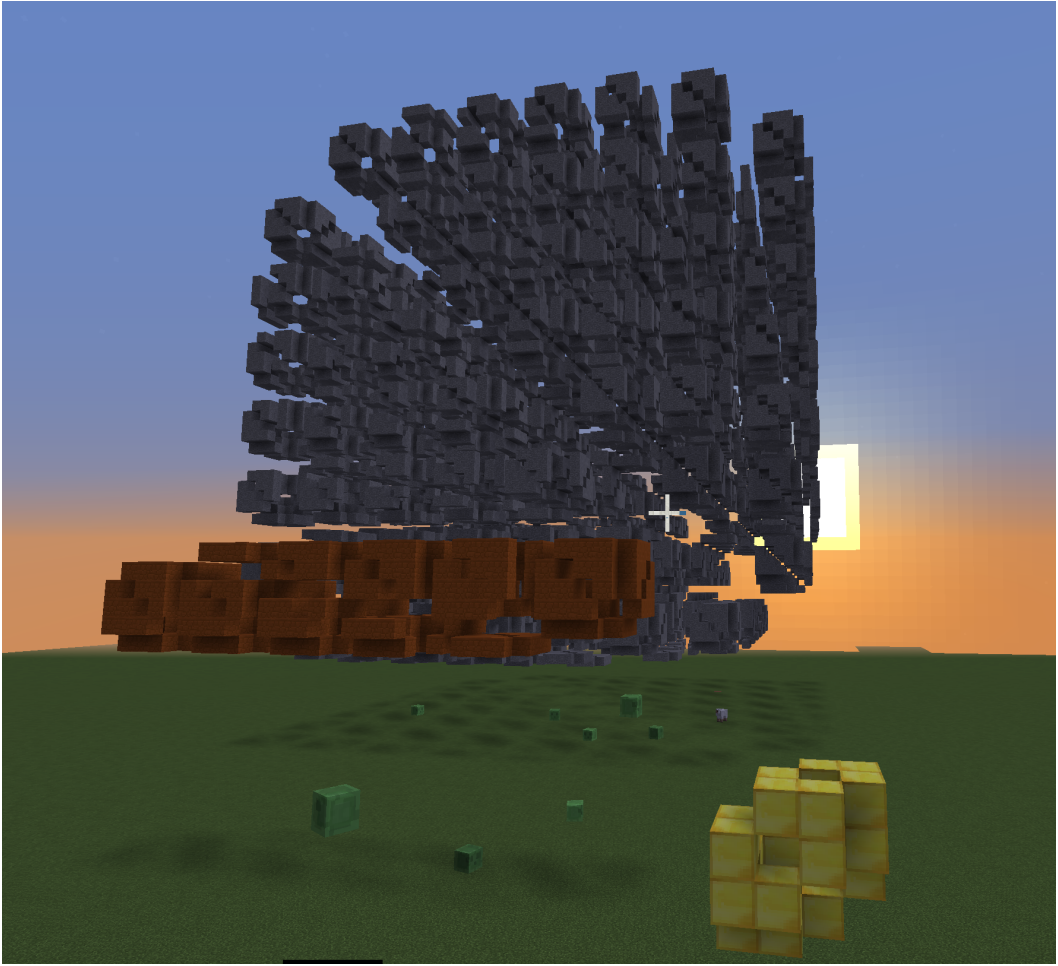

	\begin{lpic}{MonstersDusk(140mm,)}
	\end{lpic}
	\caption{Examples of homogeneous configurations}\label{f:MonstersDusk}
  \end{figure}

  A configuration $\Xcal\in\prob\<\Gbf>$ modeled on some diagram
  category $\Gbf$ is called \term[homogeneous
    configuration]{homogeneous} if its automorphism group
  $\Aut(\Xcal)$ acts transitively on every probability space in
  $\Xcal$.  Three examples of homogeneous configurations were given in
  the introduction. Other examples of a homogeneous configurations (of
  combinatorial type $\Lambdabf_{3}$) are shown in
  Figure~\ref{f:MonstersDusk}.  The
  subcategory of all homogeneous configurations modeled on $\Gbf$ will
  be denoted $\Prob\<\Gbf>_{\hbf}$.
  
  In fact, for $\Xcal$ to be homogeneous it is sufficient that the
  $Aut(\Xcal)$ acts transitively on every initial space in $\Xcal$.
  Thus, if $\Xcal$ is complete with initial space $X_{0}$, to
  check homogeneity it is sufficient to check the transitivity of the
  action of the symmetries of $\Xcal$ on $X_{0}$.

  By functoriality of the restriction operator, any restriction of a
  homogeneous configuration is also homogeneous. In other words, if
  $R:\Gbf\to\Gbf'$ is a functor and
  \[
  R^{*}:\prob\<\Gbf'>\to\prob\<\Gbf>
  \] 
  is the associated restriction operator, then
  \[
  R^{*}\big(\prob\<\Gbf'>_{\hbf}\big)\subset\prob\<\Gbf>_{\hbf}
  \]

  In particular, all the individual spaces of a homogeneous
  configuration are homogeneous 
  \[
    \prob\<\Gbf>_{\hbf}\subset\probhom\<\Gbf>
  \]
  However homogeneity of the whole of the configuration is a stronger
  property than homogeneity of the individual spaces in the
  configuration, thus in general
  \[
    \prob\<\Gbf>_{\hbf}
    \subsetneq
    \probhom\<\Gbf>
  \]

  A single probability space is homogeneous if and only if there is a
  representative in its isomorphism class with uniform measure and the
  same holds true for chain configurations, for the co-fan or any
  other configuration that does not contain a two-fan.  However, for
  more complex configurations, for example for two-fans, no such
  simple description is available.
 
 \subsubsection{Universal construction of homogeneous configurations} 
  Examples of homogeneous configurations could be constructed in the
  following manner.  Suppose $\Gamma$ is a finite group and
  $\set{H_{i}}$ is a collection of subgroups. Consider a collection of
  sets $\un X_{i}:=\Gamma/H_{i}$ and consider a natural surjection
  $f_{ij}:\un X_{i}\to \un X_{j}$ whenever $H_{i}$ is a subgroup of
  $H_{j}$. Equipping each $\un X_{i}$ with the uniform distribution
  one can turn the configuration of sets $\set{\un X_{i};f_{ij}}$ into
  a homogeneous configuration. It will be complete if there is a
  smallest subgroup (under inclusion) among $H_{i}$'s.  
  
  Such a configuration will be complete and minimal, if
  together with any pair of groups $H_{i}$ and $H_{j}$ in the
  collection, their intersection $H_{i}\cap H_{j}$ also belongs to the
  collection $\set{H_{i}}$.

  In fact, any homogeneous configuration arises this way.  Suppose
  configuration $\Xcal = \set{ X_i ; f_{ij} }$ is homogeneous, then
  we set $\Gamma = \Aut(\Xcal)$ and choose a collection of points $x_i
  \in X_i$ such that $f_{ij} (x_i) = x_j$ and denote by $H_i :=
  \stab(x_i) \subset \Gamma$.  Then, if one applies the construction
  of the previous paragraph to $\Gamma$, with the collection of
  subgroups $\set{H_i}$, one recovers the original configuration
  $\Xcal$.

\subsection{Conditioning}\label{s:config-conditioning}
  Suppose a configuration $\Xcal$ contains a fan 
  \[
  \Fcal = \left(X \oot[f] Z \too[g] Y\right)
  \] 
  Given a point $x\in X$ with a non-zero weight one may consider
  \term{conditional probability distributions} $p_Z(\;\cdot\;\rel x)$
  on $\un Z$, and $p_Y( \; \cdot \;\rel x)$ on $\un Y$.  The
  distribution $p_Z(\; \cdot \; \rel x)$ is supported on $f^{-1}(x)$
  and is given by
  \[
  p_Z( z \rel x) 
  = 
  \frac{p_Z(z)}{p_X(x)}
  \]
  The distribution $p_Y(\;\cdot\;\rel x)$ is the pushforward of $p_Z(
  \; \cdot \; \rel x)$ under $g$
  \[
   p_Y(\; \cdot \; \rel x) 
   = 
   g_* p_Z(\; \cdot\;\rel x)
  \]
  Recall that if $\Fcal$ is minimal, the underlying set of $Z$ can be
  assumed to be the product $\un X \times \un Y$. In that case
  \[
  p_Y( y \rel x ) = \frac{p_Z(x,y)}{p_X(x)}
  \]
  
  We denote the corresponding space $Y\rel x:=\big(\un
  Y,p_Y(\;\cdot\;\rel x)\big)$, as discussed at the end of
  Section~\ref{s:category-config-2fan}.

  Under some assumptions it is possible to condition a whole
  sub-configuration of $\Xcal$.  More specifically, if a configuration
  $\Xcal$ contains a sub-configuration $\Ycal$ and a probability space
  $X$ satisfying the condition that
  \begin{quote} 
    for every $Y$ in $\Ycal$ there is a fan in $\Xcal$ with terminal
    vertices $X$ and $Y$,
  \end{quote}
  then we may condition the whole of $\Ycal$ on $x \in X$ given that
  $p_X(x)>0$.
  
  For $x\in X$ with positive weight we denote by $\Ycal\rel x$ the
  configuration of spaces in $\Ycal$ conditioned on $x\in X$. The
  configuration $\Ycal\rel x$ has the same combinatorial type as
  $\Ycal$ and will be called the \term{slice} of $\Ycal$ over
  $x\in X$.  Note that the space $X$ itself may or may not belong to
  $\Ycal$. The conditioning $\Ycal\rel x$ may depend on the choice of
  a fan between $\Ycal$ and $X$, however when $\Xcal$ is complete the
  conditioning $\Ycal\rel x$ is well-defined and is independent of the
  choice of fans.
  
  Suppose now that there are two subconfiguration $\Ycal$ and
  $\Zcal$ in $\Xcal$ and in addition $\Zcal$ is a constant
  configuration, $\Zcal=Z^{\Gbf'}$ for some diagram category
  $\Gbf'$. Let $z\in\un Z$, then $\Ycal\rel z$ is well defined and is
  independent of the choice of the space in $\Zcal$, the element of
  which $z$ is to be considered.

  If $\Xcal$ is homogeneous, then $\Ycal\rel x$ is also homogeneous and
  its isomorphism class does not depend on the choice of $x\in\un X$.

\subsection{Entropy}
\label{s:config-entropy}
  For a $\Gbf$-configuration $\Xcal=\set{X_{i},f_{ij}}$ define the
  entropy function
  \[
    \ent_{*}:\prob\<\Gbf>\to\Rbb^{\size{G}},
    \quad
    \ent_{*}:\Xcal=\set{X_{i},f_{ij}}\mapsto 
    \big(\ent(X_{i})\big)\in\Rbb^{\size{G}}
  \]

  It will be convenient for us to equip the target
  $\Rbb^{\size{\Gbf}}$ with the $\ell^{1}$-norm. Thus
  \[
    |\Ent_{*}(\Xcal)|_{1}=\sum_{i=1}^{\size{\Gbf}}\ent(X_{i})
  \]

  If $\Xcal$ is a complete $\Gbf$-configuration with initial space
  $X_{0}$, then by Shannon inequality (\ref{eq:shannonineq}) there is an
  obvious estimate
  \[
    \ent(X_{0})
    \leq
    |\Ent_{*}(\Xcal)|_{1}
    \leq
    \size{\Xcal}\cdot\ent(X_{0})
  \]

	\section{Distributions and types}\label{s:disttypes}
	  In this section we recall some elementary inequalities for
  (relative) entropies and the total variation distance for
  distributions on finite sets. Furthermore, we generalize the notion
  of a probability distribution on a set to a distribution on a
  configuration of sets. Finally, we give a perspective on the theory
  of types, and also introduce types in the context of complete
  configurations.

\subsection{Distributions}
\label{s:disttypes-distributions}

\subsubsection{Single probability spaces}
\label{s:disttypes-distributions-single}
  For a finite set $S$ we denote by $\Delta S$ the collection of all
  probability distributions on $S$.  It is a unit simplex in the real
  vector space $\Rbb^S$.  We often use the fact that it is a compact,
  convex set, whose interior points correspond to fully supported
  probability measures on $S$.

  For $\pi_{1},\pi_{2}\in\Delta S$ denote by $|\pi_{1}-\pi_{2}|_1$ the total
  variation of the signed measure $(\pi_{1}-\pi_{2})$ and 
  define the entropy of the distribution $\pi_{1}$ by
  \[\tageq{entropyformula}
    h(\pi_{1})
    :=
    -\sum_{x\in X}\pi_{1}(x)\ln \pi_{1}(x)
  \]
  If, in addition, $\pi_{2}$ lies in the interior of
  $\Delta S$ define the relative entropy by
  \[\tageq{relativeentropy}
    D(\pi_{1}\sep \pi_{2})
    :=
    \sum_{x\in X}\pi_{1}(x)\ln \frac{\pi_{1}(x)}{\pi_{2}(x)}
  \]
  The entropy of a probability space is often defined through
  formula~(\ref{eq:entropyformula}).  It is a standard fact, and can
  be verified with the help of Lemma~\ref{p:typesestimates} below,
  that for $\pi\in\Delta S$ holds
  \[\tageq{entropiesequal}
    h(\pi)
    =
    \ent(S,\pi)
  \]
  which justifies the name ``entropy'' for the function $h:\Delta S\to\Rbb$.

  Define a \term{divergence ball} of radius $\epsilon>0$ centered at
  $\pi\in\Interior\Delta S$ as
  \[\tageq{divergenceball}
    B_{\epsilon}(\pi)
    :=
    \set{\pi'\in\Delta S\st D(\pi'\sep\pi)\leq \epsilon}
  \]
  For a fixed $\pi$ and $\epsilon<\!<1$ the ball $B_{\epsilon}(\pi)$
  also lies in the interior of $\Delta S$.
  \begin{lemma}{p:entropydivergence} 
  	Let $S$ be a finite set, then
    \begin{enumerate}
    \item
      \label{p:entropydivergence1}
      For any $\pi_{1},\pi_2\in\Delta S$, Pinsker's inequality holds
      \[
      |\pi_{1}-\pi_{2}|_1\leq\sqrt{2D(\pi_{1}\sep\pi_{2})}
      \]    
    \item
      \label{p:entropydivergence2} For any
      $\pi_{2}\in\Interior\Delta S$ there exists a positive constant
      $C = C_{\pi_2}$ such that for any $\pi_{1} \in \Delta S$, holds
      \[
      |\pi_{1}-\pi_{2}|_1
      \geq 
      C\sqrt{D(\pi_{1}\sep\pi_{2})}
      \]
    \item
      \label{p:entropydivergence3}  
      Suppose $\pi$ is a point in the interior of $\Delta S$ and $r>0$
      is such that $B_{r}(\pi)$ also lies in the interior of $\Delta
      S$. There exist a constant $C=C_{\pi,r}$ such that for any
      $\epsilon\leq r$ holds
      \[
      \max\set{|h(\pi_{1})-h(\pi_{2})|
        \st 
        \pi_{1},\pi_{2}\in B_{\epsilon}(\pi)}\leq C\sqrt{\epsilon}
      \]
    \end{enumerate}
  \end{lemma}
  
The first claim of the Lemma, Pinsker's inequality, is a well-known
inequality in for instance information theory, and a proof can be
found in \cite{Cover-Elements-1991}.

The second claim follows from the fact that for the fixed
$\pi_{2}\in\Interior\Delta S$ the relative entropy as the function
of the first argument is bounded, smooth on
the interior of the simplex and has a minimum at $\pi_{2}$.

To prove the last claim, note that the entropy function $h$ is
smooth in the interior of the simplex. Then this last claim follows
from the first claim.  
  
\subsubsection{Distributions on configurations}
\label{s:disttypes-distributions-config}
  A map $f:S\to S'$ between two finite sets
   induces an affine map $f_{*}:\Delta S\to\Delta S'$.

  For a configuration of sets $\Scal=\set{S_{i};f_{ij}}$ we define the
  \term[distribution on configuration of sets]{space of distributions
    on the configuration} $\Scal$ by
  \[
    \Delta\Scal
    :=
    \set{(\pi_{i})\in\prod_i\Delta S_{i}\st (f_{ij})_{*}\pi_{i}=\pi_{j}}
  \]
  Essentially, an element of $\Delta \Scal$ is a collection of
  distributions on the sets $S_i$ in $\Scal$ that is consistent with
  respect to the maps $f_{ij}$.  The consistency conditions
  $(f_{ij})_* \pi_i = \pi_j$ form a collection of linear equations
  with integer coefficients with respect to the standard convex
  coordinates in $\prod \Delta S_i$. Thus, $\Delta \Scal$ is a
  rational affine subspace in the product of simplices. In particular,
  $\Delta \Scal$ has a convex structure.
  
  If $\Scal$ is complete with initial set $S_{0}$, then specifying a
  distribution $\pi_{0}\in\Delta S_{0}$ uniquely determines
  distributions on all of the $S_{i}$'s by setting
  $\pi_{i}:=(f_{0i})_{*}\pi_{0}$. In such a situation we have
  \[\tageq{distribonconfig}
    \Delta\Scal\cong\Delta S_{0}
  \]
  If $\Scal$ is not complete and $S_{0},\ldots,S_{k}$ is a collection
  of its initial sets, then $\Delta\Scal$ is isomorphic to an affine
  subspace of the product $\Delta S_{0}\times\dots\times\Delta S_{k}$
  cut out by linear equations with integer coefficients corresponding
  to co-fans in $\Scal$ with initial sets among $S_{0},\ldots,S_{k}$.
  
  To simplify notation, for a probability space $X$ or a
  configuration $\Xcal$ we will write
  \begin{align*}
    \Delta X
    &:=
    \Delta\un X
    \\
    \Delta\Xcal
    &:=
    \Delta\un\Xcal
  \end{align*}
  
  \bigskip
  We now discuss briefly the theory of types. Types are special
  subspaces of tensor powers that consist of seqences with the same
  ``empirical distribution'' as explained in details below. For a
  more detailed discussion the reader is referred to
  \cite{Cover-Elements-1991} and \cite{Csiszar-Method-1998}. We generalize the
  theory of types to complete configurations of sets and complete
  configurations of probability spaces.

  The theory of types for configurations, that are not complete, is
  more complex and will be addressed in a
    subsequent article.

\subsection{Types for single probability spaces}
\label{s:disttypes-types-single}
  Let $S$ be a finite set. For $n\in\Nbb$ denote also
  \[
  \Deltan S
  :=
  \Delta S\cap \frac1n\Zbb^{S}
  \]
  a collection of rational points in $\Delta S$ with denominator $n$.
  (We say that a rational number $r\in\Qbb$ has denominator $n\in\Nbb$
  if $r\cdot n\in\Zbb$)

  Define the \term{empirical distribution map} $\emp:S^{n}\to\Delta S$,
  that sends $(s_{i})_{i=1}^{n}=\sbf\in S^{n}$ to the empirical
  distribution $\emp(\sbf)\in\Delta S$ given by 
  \[
  \emp(\sbf)(a)
  =
  \frac1n\cdot\big|\set{i\st s_{i}=a}\big|
  \quad
  \text{for any $a\in S$}
  \]
  Clearly the image of $\emp$ lies in $\Deltan S$.
  
  For $\pi\in\Deltan S$, the space $T^{\c n}_{\pi}S:=\emp^{-1}(\pi)$
  equipped with the uniform measure is called a \term{type} over $\pi$.
  The symmetric group $\Sbb_{n}$ acts on $S^{\otimes n}$ by permuting the
  coordinates. This action leaves the empirical distribution
  invariant and therefore could be restricted to each type, where it
  acts transitively. Thus, for $\pi\in\Deltan S$ the
  probability space $(T^{\c n}_{\pi}S,u)$ with $u$ being a uniform
  ($\Sbb_{n}$-invariant)
  distribution, is a homogeneous space.
  
  Suppose $X=(\un X,p)$ is a probability space.  Let $\tau_{n}$ be the
  pushforward of $p^{\otimes n}$ under the empirical distribution map
  $\emp:\un X^{n}\to\Delta X$ . Clearly $\supp\tau_{n}\subset\Deltan
  X$, thus $(\Delta X,\tau_{n})$ is a finite probability
  space.  Therefore we have a reduction
  \[
  \emp:X^{\otimes n}\to(\Delta X,\tau_{n})
  \]
  which we call the \term{empirical reduction}.
  If $\pi\in\Deltan X$ is such that $\tau_{n}(\pi)>0$, then 
  \[\tageq{typesiso}
  T^{\c n}_{\pi}\un X = X^{\otimes n}\rel\pi
  \]
  In particular, it follows that the right-hand side does not depend on
  the probability $p$ on $X$ as long as $\pi$ is ``compatible'' to it.
  
  The following lemma records some standard facts about types, which
  can be checked by elementary combinatorics and found in
  \cite{Cover-Elements-1991}.
  \newpage
  \begin{lemma}{p:typesestimates}
    Let $X$ be a probability space and $\xbf\in X^{\otimes n}$, then
    \begin{enumerate}
    \item 
      \[
      |\Deltan X|
      =
      \choose{n+|X|\\|X|}\leq \ebf^{|X| \cdot\ln (n+1)}
      \]
    \item
      \[
      p^{\otimes n}(\xbf)
      =
      \ebf^{-n\big[h\big(\emp(\xbf)\big)+D\big(\emp(\xbf)\sep p\big)\big]}
      \]
    \item
      \[
      \ebf^{n \cdot h(\pi) - |X| \cdot\ln(n+1)}\leq
      |T^{\c n}_{\pi} \un X|
      \leq
      \ebf^{n\cdot h(\pi)}
      \]
     \item
      \[
      \ebf^{-n\cdot D(\pi\sep p)-|X| \cdot\ln (n+1)}
      \leq
      \tau_n(\pi)=p^{\otimes n}(T^{\c n}_{\pi} \un X)
      \leq
      \ebf^{-n\cdot D(\pi\sep p)}
      \]
    \end{enumerate}
  \end{lemma}

  If $X=(\un X,p_{X})$ is a probability space with rational probability
  distribution with denominator $n$, then the type over $p_{X}$ will be
  called the true type of $X$
  \[
  T^{\c n}X
  :=
  T^{\c n}_{p_{X}}\un X
  \]
  
  As a corollary to Lemma \ref{p:typesestimates} and
  equation~(\ref{eq:entropiesequal}) we obtain the following.
  
  \begin{corollary}{p:entropy-type}
    For a finite set $S$ and $\pi\in\Deltan S$ holds
    \[
      n\cdot h(\pi)- |S| \cdot \ln (n+1)
      \leq
      \ent(T^{(n)}_{\pi}S)
      \leq
      n\cdot h(\pi)
    \]

    In particular, for a finite probability space $X=(S,p)$ with a rational
    distribution $p$ with denominator $n$ holds

    \[
      n\cdot\ent(X)- |S| \cdot \ln (n+1)
      \leq
      \ent(T^{(n)}X)
      \leq
      n\cdot \ent(X)
    \]
    
  \end{corollary}
  
  The following important theorem is known as Sanov's theorem. 
  It can be derived from Lemma
  \ref{p:typesestimates} and found in \cite{Cover-Elements-1991}.
  
  \begin{theorem}{p:sanov}{\rm (Sanov's Theorem)}
    Let $X$ be a finite probability space and let $\emp: X^{\otimes n}
    \to (\Delta X, \tau_n)$ be the empirical reduction.  Then for
    every $r>0$,
    \[
      \tau_n(\Delta X \backslash B_r(p)) 
      \leq 
      \ebf^{- n \cdot r + |X| \cdot \ln( n + 1 )}
    \]
    where $B_r(p)$ is the divergence ball (relative entropy ball)
    defined in (\ref{eq:divergenceball}).
  \end{theorem}

\subsection{Types for complete configurations}
\label{s:disttypes-types-config}
  In this subsection we generalize the theory of types for
  configurations modeled on a complete category. The theory for a
  non-complete configurations is more complex and will be addressed in
  our future work.  We will give three equivalent definitions of a
  type for a complete configuration, each of which will be useful in
  its own way. Before we describe the three approaches we need some
  preparatory material.
  \begin{lemma}{p:specialdiamond}
    Given a diamond configuration of probability spaces
    \[
    \Dcal=
    \left(
    \begin{tikzcd}[ampersand replacement=\&,row sep=small,column sep=small]
      X
      \arrow{d}{\rho_{1}}
      \arrow{r}{f}
      \&
      Y
      \arrow{d}{\rho_{2}}
      \\
      A
      \arrow{r}{g}
      \&
      B
    \end{tikzcd}
    \right)
    \]
    the following two conditions are equivalent
    \begin{enumerate}
    \item
      \label{p:specialdiamond-indep} 
      The minimal reduction of the diamond $\Dcal$ is isomorphic to the
      adhesion of its co-fan or equivalently the following independence
      condition holds
      \[
      A\indep Y\rel B
      \]
    \item
      \label{p:adhesion-cond}
      For any $a,a'\in \un A$ such that $f(a)=f(a')=b\in B$ holds
      \[
      Y\rel a=Y\rel a'=Y\rel b
      \]
    \end{enumerate}
  \end{lemma}

  Suppose $\Dcal$ is the diamond as in the
  Lemma~\ref{p:specialdiamond}.  The top row $X\to Y$ is a two-chain
  subconfiguration of $\Dcal$ and we can consider a conditioning by an
  element $a\in\un A$
  \[
  X\rel a\to Y\rel a
  \]
  If $\Dcal$ satisfies any of two conditions in
  Lemma~\ref{p:specialdiamond}, then $Y\rel a=Y\rel b$ for
  $b=g(a)$. Thus, we constructed a reduction 
  \[
  X\rel a\to Y\rel b  
  \]
  Suppose we have a reduction $f:X\to Y$ between a pair of probability
  spaces. Then for any $n\in\Nbb$ there is an induced reduction
  $f_{*}:(\Delta X,\tau_{n})\to(\Delta Y,\tau_{n})$ that can be
  included in the following diamond configuration
  \[
  \begin{tikzcd}[ampersand replacement=\&]
    X^{\otimes n} 
    \arrow{r}{f^{\otimes n}}
    \arrow{d}{\emp}
    \& 
    Y^{\otimes n} 
    \arrow{d}{\emp}
    \\
    (\Delta X,\tau_n)
    \arrow{r}{f_{*}}
    \& 
    (\Delta Y,\tau_{n})
  \end{tikzcd}
  \]
  that satisfies conditions in Lemma~\ref{p:specialdiamond}.
  It means that there is a reduction
  \[
  Tf:T^{(n)}_{\pi} X\to T^{(n)}_{\pi'}Y
  \]
  for $\pi\in\Deltan X$ and $\pi'=f_{*}\pi\in\Deltan Y$.

  Now we are ready to give the definitions of types.  Let
  $\Xcal\in\prob\<\Gbf>$ be a complete configuration,
  $\Xcal=\set{X_{i};f_{ij}}$ with initial space $X_{0}$ and let 
  $\pi\in\Deltan \Xcal$.

\subsubsection{Type of a configuration as the configuration of types.}
\label{s:type-config-analytic}
  Define the type $T^{(n)}_{\pi}\un\Xcal$ as the $\Gbf$-configuration,
  whose individual spaces are types of the individual spaces of
  $\Xcal$ over the corresponding push-forwards of $\pi$
  \[
  T^{(n)}_{\pi}\un\Xcal
  :=
  \set{T^{(n)}_{\pi_{i}}\un X_{i};T f_{ij}}
  \]
\subsubsection{Types as $\Sbb_{n}$-orbits in the tensor power}
\label{s:type-config-symmetric}
  By a \term{section} in $\Xcal$ we mean a consistent collection of points
  \[
  \xbf=(x_{0},\ldots,x_{\size{\Xcal}-1})\in \prod_{i=0}^{\size{\Xcal}-1}\un X_{i}
  \]
  such that $f_{ij}x_{i}=x_{j}$, whenever $f_{ij}$ is defined.  For
  any $j$ define the projection
  $\rho_{j}:\prod\un X_{i}\to\un X_{j}$, so that
  $\rho_{j}(\xbf)=x_{j}$.

  The symmetric group $\Sbb_{n}$ acts on the collection of sections in
  the tensor power $\Xcal^{\otimes n}$, by permuting the coordinates.
  Let $\Tcal\subset\prod\un X_{i}$ be an orbit of the action such that
  $\rho_{0}(\Tcal)=T^{(n)}_{\pi}\un X_{0}$.  Suppose
  that the pair $(i,j)$ is such, that $f_{ij}$ is defined. Since
  $f_{ij}^{\otimes n}:X_{i}^{\otimes n}\to X_{j}^{\otimes n}$ is
  $\Sbb_{n}$-equivariant, we have a map
  \[
  Tf:\rho_{i}(\Tcal)\to\rho_{j}(\Tcal)
  \]
  We can turn $\Tcal$ into a $\Gbf$-configuration, which we will call
  a type of $\Xcal$ 
  \[
  \hat T^{(n)}_{\pi}\un\Xcal:=\set{\rho_{i}(\Tcal),Tf_{ij}}
  \]
  where $\pi$ is the value of the impirical distribution on
  $\rho_{0}(\Tcal)$. 
  
  Since the initial space in $\hat T^{(n)}_{\pi}\un\Xcal$ coincides
  with the initial space in $T^{(n)}_{\pi}\un\Xcal$ and all the
  reductions coincide, we conclude that 
  \[
  \hat T^{(n)}_{\pi}\un\Xcal
  =
  T^{(n)}_{\pi}\un\Xcal
  \]

\subsubsection{Type as conditionining of the tensor power}
\label{s:type-config-synthetic}
  We can extend $\Xcal^{\otimes n}$ to a configuration $\hat{\Xcal}$ by
  adding $(\Delta \Xcal, \tau_n)$ and the empirical reduction
  $X_0^{\otimes n}\to (\Delta \Xcal, \tau_n)$.
  
  Let $\pi \in \Delta \Xcal$ with $\tau_n(\pi)>0$ and recall $\Delta
  X_0 \cong \Delta \Xcal$.  We may now define $\Xcal^{\otimes n} \rel
  \pi$ as in Section \ref{s:config-conditioning}.  
  Define a type of $\Xcal$ over $\pi\in\Deltan\Xcal$ by
  \[
  \check T^{(n)}_{\pi}\un\Xcal:=\Xcal^{\otimes n}\rel\pi
  \]
  
  By definition, it holds that
  \[
  X_0^{\otimes n} \rel \pi 
  = 
  T_\pi^{\c n} X_0
  \]
  Let $\pi_{i}=(f_{0i})_{*}\pi$.
  Using Lemma~\ref{p:specialdiamond} and discussion thereafter we
  conclude that 
  \[
  X_{i}^{\otimes n}\rel\pi=X_{i}^{\otimes n}\rel\pi_{i}
  \] 
  and therefore
  \[
  \check T^{(n)}_{\pi}\un\Xcal
  =
  T^{(n)}_{\pi}\un\Xcal
  \]
  
\subsubsection{The empirical two-fan}
\label{s:types-empirical-twofan}
  We construct a two-fan of $\Gbf$-configurations with terminal
  vertices $\Xcal^{\otimes n}$ and $(\Delta \Xcal,\tau_{n})^{\Gbf}$
  (the constant $\Gbf$-configuration in which every probability space
  is $(\Delta \Xcal, \tau_n)$, and every reduction is an identity)
  \[\tageq{empirical-2fan}
  \Rcal_{n}(\Xcal)
  =
  \left(
  \begin{tikzcd}[column sep=tiny, ampersand replacement=\&]
  \mbox{}
  \&
  \tilde\Xcal^{(n)}
  \arrow{ld}{}
  \arrow{rd}{}
  \&
  \mbox{}
  \\
  \Xcal^{\otimes n}
  \&
  \&
  (\Deltan \Xcal,\tau_{n})^{\Gbf}
  \end{tikzcd}
  \right)
  \]
  With the help of Lemma \ref{p:minimalfansconfig}, we construct
  $\Rcal_n(\Xcal)$ as the minimal reduction of the two-fan of
  $\Gbf$-configurations
  \[
  \begin{tikzcd}[column sep=tiny,ampersand replacement=\&]
  \mbox{}
  \&
  (X_{0}^{\otimes n})^{\Gbf}
  \arrow{dl}[swap]{f_{0*}^{n}}
  \arrow{dr}{\emp^{\Gbf}}
  \&
  \mbox{}
  \\
  \Xcal^{\otimes n}
  \&
  \&
  (\Delta\Xcal,\tau_{n})^{\Gbf}
  \end{tikzcd}
  \]

  Let $\pi \in (\Delta \Xcal, \tau_n)$ with $\tau_n(\pi) > 0$. 
  Then within $\Rcal_{n}(\Xcal)$ holds
  \[
  \Xcal^{\otimes n}\rel\pi=T^{(n)}_{\pi}\Xcal
  \]

\bigskip

  For every $n\in\Nbb$ and $\pi\in\Deltan\un X_{0}$ the type $T^{\c
    n}_{\pi}\Xcal$ is a homogeneous configuration. Suppose that a complete configuration
  $\Xcal$ is such that the probability distribution $p_{0}$ on the
  initial set is rational with the denominator $n$, then we call
  $T^{\c n}_{p}\un\Xcal$ the \term{true type} of $\Xcal$ and
  denote
  \[
  T^{\c n}\Xcal:=T^{n}_{p_{0}}\un\Xcal
  \]

	\section{The Kolmogorov-Sinai distance}\label{s:kolmogorov}
	  We turn the space of configurations into a pseudo-metric space by
  introducing the intrinsic Kolmogorov-Sinai distance and asymptotic
  Kolmogorov-Sinai distance. For brevity, we will usually call it the
  Kolmogorov distance and asymptotic Kolmogorov distance.  The
  intrinsic Kolmogorov-Sinai distance is obtained by taking an infimum
  of the shared information distance over all possible joint
  distributions on two probability spaces. The name is justified by
  the fact that the shared information distance (not under this name)
  appears in the proof of the theorem about generating partitions for
  ergodic systems by Kolmogorov and Sinai, see for
  example \cite{Sinai-ergodic-1976}.  Note that the Kolmogorov
  distance in statistics refers to a different notion.

\subsection{Kolmogorov distance and asymptotic Kolmogorov distance}
\label{s:kolmogorov-ikd-aikd}
\subsubsection{Kolmogorov distance in the case of single probability spaces}
\label{s:kolmogorov-single}
  For a two-fan $\Fcal=(X\ot Z\to Y)$ define a ``distance'' $\kd(\Fcal
  )$ between probability spaces $X$ and $Y$ with respect to $\Fcal$ by
  \[
  \begin{split}
 \kd(\Fcal )&:=\ent(Z\rel Y)+\ent(Z\rel X) \\
 &= 2\ent(Z)-\ent(X)-\ent(Y)
 \end{split}
  \]

  Essentially $\kd(\Fcal )$ measures the deviation of the statistical map
  defined by $\Fcal$ from being a deterministic bijection between $X$ and
  $Y$.
  
  The minimal reduction $\Fcal '$ of $\Fcal$ satisfies
  \[\tageq{kd-min-reduction}
  \kd(\Fcal ') \leq \kd (\Fcal)
  \]

  If the two-fan $\Fcal$ is minimal the ``distance''
  $\kd( \Fcal )$ can also be calculated by 
  \[ 
  \kd(\Fcal) = h( p_X ) +
  h( p_Y ) - 2 D( p_Z \sep p_X \otimes p_Y ), 
  \] 
  where $h$ and $D$ are
  respectively the entropy and relative entropy functions defined in
  (\ref{eq:entropyformula}) and (\ref{eq:relativeentropy}).

  For a pair of probability spaces $X$, $Y$ define the
  \term[Kolmogorov-Sinai distance (spaces)]{intrinsic Kolmogorov-Sinai distance} as
  \[
  \ikd(X,Y)
  :=
  \inf\set{\kd(\Fcal)\st \Fcal=(X\ot Z \to Y) 
    \text{ is a two-fan} 
  }
  \]

  The optimization takes place over all two-fans with terminal spaces
  $X$ and $Y$. In view of inequality~(\ref{eq:kd-min-reduction}) one
  could as well optimize over the space of \emph{minimal} two-fans,
  which we will also refer to as \term[coupling between probability
    spaces]{couplings} between $X$ and $Y$. The tensor product of $X$
  and $Y$ trivially provides a coupling and the set of couplings is
  compact, therefore an optimum is always achieved and it is finite.

  The bivariate function $\ikd:\Prob\times\Prob\to\Rbb_{\geq0}$
  defines a notion of pseudo-distance and it vanishes exactly on
    pairs of isomorphic probability spaces. This follows directly
  from the Shannon inequality~(\ref{eq:shannonineq}), and a more
  general statement will be proven in Proposition
  \ref{p:kolmogorovisdistance} below.

\subsubsection{Kolmogorov distance for complete configurations}
\label{s:kolmogorov-config}
  The definition of Kolmogorov distance for complete configurations
  repeats almost literally the definition for single spaces. We fix a
  complete diagram category $\Gbf$ and will be considering configurations
  from $\prob\<\Gbf>$.
  
  Consider three configurations $\Xcal=\set{X_{i},f_{ij}}$,
  $\Ycal=\set{Y_{i},g_{ij}}$ and $\Zcal=\set{Z_{i},h_{ij}}$ from
  $\prob\<\Gbf>$. Recall that a two-fan $\Fcal=(\Xcal\ot\Zcal\to\Ycal)$ is
  a $\Gbf$-configuration of two-fans
  \[
  \Fcal_{i}=(X_{i}\ot Z_{i}\to Y_{i})
  \]
  
  Define
  \begin{align*}
    \kd(\Fcal)
    &:=
    \sum_{i} \kd(\Fcal_{i})\\
    &=
    \sum_{i}\big(2\ent(Z_{i})-\ent(X_{i})-\ent(Y_{i})\big)
  \end{align*}
  
  The quantity $\kd(\Fcal)$ vanishes if and only if the fan $\Fcal$
  provides isomorphisms between all individual spaces in $\Xcal$ and
  $\Ycal$ that commute with the inner structure of the configurations,
  that is, it provides an isomorphism between $\Xcal$ and $\Ycal$ in
  $\prob\<\Gbf>$.
  
  The \term[Kolmogorov-Sinai distance (configurations)]{intrinsic
    Kolmogorov-Sinai distance between configurations} is defined in analogy
  with the case of single probability spaces
  \[
  \ikd(\Xcal,\Ycal)
  :=
  \inf\set{\kd(\Fcal) \st \Fcal=(\Xcal\ot\Zcal\to\Ycal)}
  \]
  where the infimum is over all two-fans of $\Gbf$-configurations with
  terminal vertices $\Xcal$ and $\Ycal$.
  
  The following proposition records that the intrinsic
  Kolmogorov distance is in fact a pseudo-distance on $\prob\<\Gbf>$,
  provided $\Gbf$ is a complete diagram category (that is when $\Gbf$
  has a unique initial space).
  
  \begin{proposition}{p:kolmogorovisdistance}
    Let $\Gbf$ be a complete diagram category.
    Then the bivariate function 
    \[
    \ikd:\prob\<\Gbf> \times \prob\<\Gbf>\to\Rbb
    \]
    is a pseudo-distance on $\Prob\<\Gbf>$.  \\
    Moreover, two
    configurations $\Xcal, \Ycal \in \prob\<\Gbf>$ satisfy
    $\ikd(\Xcal,\Ycal)=0$ if and only if $\Xcal$ is isomorphic to
    $\Ycal$ in $\prob\<\Gbf>$.
  \end{proposition}

  The idea of the proof is very simple.  In the case of single
  probability spaces $X, Y, Z$ a coupling between $X$ and $Z$ can be
  constructed from a coupling between $X$ and $Y$ and a coupling
  between $Y$ and $Z$ by adhesion on $Y$, see Section \ref{s:config-adhesion}.  The triangle inequality
  then follows from a Shannon inequality.  However, since we are
  dealing with configurations the combinatorial structure requires
  careful treatment.  Therefore, we provide a detailed proof on
  page \pageref{p:kolmogorovisdistance.rep}.

  It is important to note, that the proof uses the fact that $\Gbf$ is
  complete. In fact, even though the definition of $\ikd$ could be
  easily extended to some bivariate function on the space of
  configurations of any fixed combinatorial type, it fails to satisfy
  the triangle inequality in general, because the composition of
  couplings requires completeness of $\Gbf$.

\subsubsection{The asymptotic Kolmogorov-Sinai distance}
\label{kolmogorov-aikd}
Let $\Gbf$ be a complete diagram category.  We define the
\term{asymptotic Kolmogorov-Sinai distance} between two
configurations $\Xcal, \Ycal \in \prob\< \Gbf>$ by
\[\tageq{definitionaikd}
  \aikd( \Xcal, \Ycal ) 
  = 
  \lim_{n \to \infty} \frac{1}{n}
  \ikd(\Xcal^{\otimes n}, \Ycal^{\otimes n}).
\]
We will show in Corollary \ref{p:subadditivity}, that the sequence 
\[
  n \mapsto \ikd(\Xcal^{\otimes n}, \Ycal^{\otimes n})
\]
is subadditive, and therefore the limit in the definition
(\ref{eq:definitionaikd}) of $\aikd(\Xcal, \Ycal)$ always exists
and for all $n\in\Nbb$ holds
\[\tageq{aikd-ikd-bound}
  \aikd(\Xcal, \Ycal ) 
  \leq 
  \frac1n\cdot
  \ikd(\Xcal^{\otimes n}, \Ycal^{\otimes n}).
\]

As a corollary of Proposition~\ref{p:kolmogorovisdistance} and 
definition~(\ref{eq:definitionaikd}) we immediately obtain 
that also the asymptotic Kolmogorov-Sinai distance is a pseudo-distance on $\prob\<\Gbf>$.
\begin{corollary}{p:aikdisdistance}
  Let $\Gbf$ be a complete diagram category.
  Then the bivariate function 
  \[
    \aikd:\prob\<\Gbf> \times \prob\<\Gbf>\to\Rbb
  \]
  is a pseudo-distance on $\Prob\<\Gbf>$ satisfying the following
  homogeneity property.
  For any pair of configurations $\Xcal,\Ycal\in\Prob\<\Gbf>$ and any
  $n\in\Nbb_{0}$ holds
  \[
  \aikd(\Xcal^{\otimes n},\Ycal^{\otimes n})=n\cdot\aikd(\Xcal,\Ycal)
  \]
\end{corollary}

We will show in a later section, however, that there
are probability spaces $X$ and $Y$ for which $\aikd(X, Y) = 0$ that
are not isomorphic.

\bigskip 

In the rest of this section we derive some elementary properties of the intrinsic 
Kolmogorov distance and the asymptotic Kolmogorov distance.

\subsection[Lipschitz property]
{Lipschitz property for operations}
\label{s:kolmogorov-lipschitz}
  In this section we show that certain natural operations on
  configurations, namely the tensor product, entropy function and
  restriction operator, are Lipschitz continuous.  In Section
  \ref{s:extensions} we will show Lipschitz continuity of certain
  extension operations.

\subsubsection{Tensor product}\label{s:kolmogorov-lipschitz-tensor}
  We show that the tensor product on the space of configurations is
  1-Lipschitz. Later this will allow us to give a simple description
  of tropical configurations, that is of points in the asymptotic cone of
  $\prob\<\Gbf>$, as limits of certain sequences of ``classical''
  configurations.

\begin{proposition}{p:tensor1lip}
  Let $\Gbf$ be a complete diagram category.
  Then with respect to the Kolmogorov distance on $\prob\<\Gbf>$
  the tensor product
  \[
  \otimes:(\Prob\<\Gbf>,\ikd)^2\to(\Prob\<\Gbf>,\ikd)
  \]
  is 1-Lipschitz in each variable, that is, for every triple $\Xcal,
  \Ycal, \Ycal' \in \Prob\<\Gbf>$ the following
    bound holds
  \[
  \ikd(\Xcal \otimes \Ycal, \Xcal \otimes \Ycal') 
  \leq 
  \ikd(\Ycal,\Ycal')
  \]
\end{proposition}

This statement is a direct consequence of additivity of entropy with
respect to the tensor product. Details can be found on
page~\pageref{p:tensor1lip.rep}.

It follows directly from definition~(\ref{eq:definitionaikd}) and
Proposition~\ref{p:tensor1lip}, that the asymptotic Kolmogorov distance
enjoys a similar property.
\begin{corollary}{p:aikdtensor1lip}
  Let $\Gbf$ be a complete diagram category.
  Then with respect to the Kolmogorov distance on $\prob\<\Gbf>$
  the tensor product
  \[
  \otimes:(\Prob\<\Gbf>,\aikd)^2\to(\Prob\<\Gbf>,\aikd)
  \]
  is 1-Lipschitz in each variable.
\end{corollary}

As another corollary we obtain the subadditivity properties of the intrinsic
Kolmogorov distance and asymptotic Kolmogorov distance.

\begin{corollary}{p:subadditivity}
  Let $\Gbf$ be a complete diagram category and let $\Xcal, \Ycal,
  \Ucal, \Vcal \in \prob\<\Gbf>$, then
  \[
  \ikd(\Xcal \otimes \Ucal, \Ycal \otimes \Vcal)
  \leq \ikd(\Xcal, \Ycal) + \ikd( \Ucal, \Vcal ).
  \]	
  and 
  \[
  \aikd(\Xcal \otimes \Ucal, \Ycal \otimes \Vcal)
  \leq 
  \aikd(\Xcal, \Ycal) + \aikd( \Ucal, \Vcal ).
  \]	
\end{corollary}

It implies in particular that shifts are non-expanding maps in
$(\prob\<\Gbf>,\ikd)$ or $(\prob\<\Gbf>,\aikd)$. 
\begin{corollary}{p:shiftcontracting}
  Let $\Gbf$ be a complete diagram category and $\dist=\ikd,\aikd$ be
  either Kolmogorov distance or asymptotic Kolmogorov distance on
  $\prob\<\Gbf>$. Let
  $\Ucal\in\prob\<\Gbf>$. Then the shift map
  \[
  \Ucal\otimes\cdot:(\prob\<\Gbf>,\dist)\to(\prob\<\Gbf>,\dist),
  \quad
  \Xcal\mapsto \Ucal\otimes\Xcal
  \]
  is a non-expanding map with respect to either Kolmogorov distance or
  asymptotic Kolmogorov distance.
\end{corollary}

\subsubsection{Entropy}\label{s:kolmogorov-lipschitz-entropy}
Recall that we defined the entropy function 
\[
\ent_{*}:\prob\<\Gbf>\to\Rbb^{\size{\Gbf}}
\]
by evaluating the entropy of all individual spaces in a
$\Gbf$-configuration. The target space $\Rbb^{\size{\Gbf}}$ will be
endowed with the $\ell^{1}$-norm with respect to the natural coordinate
system. With such a choice, the entropy function is 1-Lipschitz with
respect to the Kolmogorov distance on $\prob\<\Gbf>$.

\begin{proposition}{p:entropy1lip}
  Suppose $\Gbf$ is a complete diagram category and $\dist=\ikd,\aikd$
  is either Kolmogorov distance or asymptotic Kolmogorov distance on
  $\prob\<\Gbf>$.  Then the entropy
  function
  \[
  \ent_{*}:(\prob\<\Gbf>,\dist)\to(\Rbb^{\size{\Gbf}},\;|\cdot|_{1}),
  \quad
  \Xcal=\set{X_{i},f_{ij}}\mapsto (\ent X_{i})_{i}\in\Rbb^{\size{\Gbf}}
  \]
  is 1-Lipschitz.
\end{proposition}
  Again, the proof of the proposition above is an application of Shannon's
  inequality, see page~\pageref{p:entropy1lip.rep} for details.

\subsubsection{Restrictions}\label{s:kolmogorov-lipschitz-restrictions}
 The restriction operators are also Lipschitz, as shown in the next
 proposition.
  \begin{proposition}{p:restriction1lip}
    Suppose $R:\Gbf'\to\Gbf$ is a functor between two complete diagram
    categories and $\dist$ stands for either Kolmogorov or asymptotic
    Kolmogorov distance. Then the restriction operator
    \[
    R^{*}:(\Prob\<\Gbf>,\dist)\to(\Prob\<\Gbf'>,\dist),
    \quad 
    \Xcal\mapsto \Xcal\circ R
    \]
    is Lipschitz.
  \end{proposition}

  As can be seen from the proof on
  page~\pageref{p:restriction1lip.rep}, the Lipschitz constant in the
  proposition above can be bounded by $\size{\Gbf'}$. In fact, a more
  careful analysis provides a better bound by the maximal number of
  objects in $\Gbf'$ that are mapped by $R$ to a single object in $\Gbf$.

\subsection{The Slicing Lemma}\label{s:kolmogorov-slicing}
  The Slicing Lemma, Proposition \ref{p:slicing} below, allows to estimate the
  Kolmogorov distance between two configurations with the integrated
  Kolmogorov distance between ``slices'', which are configurations obtained by
  conditioning on another probability space.

  The Slicing Lemma, along with the local estimate in Section~\ref{s:kolmogorov-local}, turned out to be a
  very powerful tool for estimation of the Kolmogorov distance and will be
  used below on many occasions.

  As described in Section \ref{s:config-constantconfig}, by a reduction
  of a configuration $\Xcal=\set{X_{i},f_{ij}}$ to a single space $U$ we
  mean a collection of reductions $\set{\rho_{i}:X_{i}\to U}$ from the
  individual spaces in $\Xcal$ to $U$, that commute with the reductions
  within $\Xcal$
  \[
  \rho_{j}\circ f_{ij}=\rho_{i}
  \]
  Alternatively, whenever a single probability space appears together
  with a $\Gbf$-configuration in a commutative diagram, it should be
  replaced by a constant $\Gbf$-configuration.
  \begin{proposition}{p:slicing}
    {\rm (Slicing Lemma)} Suppose $\Gbf$ is a complete diagram
    category and we are given
    $\Xcal,\hat\Xcal,\Ycal,\hat\Ycal\in\prob\<\Gbf>$ -- four
    $\Gbf$-configurations and $U,V,W\in\prob$ -- probability spaces,
    that are included into the following three-tents configuration
    
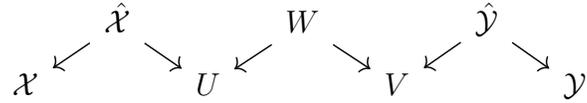
\[
    \begin{tikzcd}[column sep=small,row sep=tiny,ampersand replacement=\&]
      \mbox{}
      \&
      \hat\Xcal
      \arrow{dl}{}
      \arrow{dr}{}
      \&
      \&
      W
      \arrow{dl}{}
      \arrow{dr}{}
      \&
      \&
      \hat\Ycal
      \arrow{dl}{}
      \arrow{dr}{}
      \&
      \mbox{}
      \\
      \Xcal
      \&
      \&
      U
      \&
      \&
      V
      \&
      \&
      \Ycal
    \end{tikzcd}
    \]
    such that the two-fan $(U\ot W\to V)$ is minimal.
    Then the following estimate holds
    \begin{align*}
      \ikd(\Xcal,\Ycal) 
      &\leq 
      \int_{W}\ikd(\Xcal\rel u,\Ycal\rel v)\d p_{W}(u,v)\\
      &\quad+
      \size{\Gbf}\cdot\kd(U\ot W\to V)\\
      &\quad+
      \sum_{i}\big[\ent(U\rel X_{i})+\ent(V\rel Y_{i})\big]
    \end{align*}
  \end{proposition}
  The idea of the proof of the Slicing Lemma (page~\pageref{p:slicing.rep})
  is as follows. For every pair $(u,v)\in\un W$ we consider an optimal
  two-fan $\Gcal_{uv}$ coupling $\Xcal\rel u$ and $\Ycal\rel v$.
  These fans have the same underlying configuration of sets.
  Then we construct a coupling between $\Xcal$ and $\Ycal$ as a convex combination of distributions of
  $\Gcal_{uv}$'s weighted by $p_{W}(u,v)$. The estimates on the
  resulting two-fan then imply the proposition.

  Various implications of the Slicing Lemma are summarized in the next
  corollary.
  
  \begin{corollary}{p:slicingcorollary}
    Let $\Gbf$ be a complete diagram category, $\Xcal,\Ycal\in\prob\<\Gbf>$ and
    $U\in\prob$.
    \begin{enumerate}
    \item
      \label{i:slicing2tents}
      Given a ``two-tents'' configuration
      \[
      \Xcal\ot\hat\Xcal\to U\ot\hat\Ycal\to\Ycal
      \]
      the following inequality holds
      \begin{align*}
      \ikd(\Xcal,\Ycal) 
      &\leq 
      \int_{U}\ikd(\Xcal\rel u,\Ycal\rel u)\d p_{U}(u)
      +2\cdot\size{\Gbf}\cdot\ent(U)
      \end{align*}
    \item 
    \label{p:slicingtwofan}
      Given a fan
      \[
      \Xcal\ot\hat\Xcal\to U
      \]
      the following inequality holds
      \begin{align*}
        \ikd(\Xcal,\Ycal)
        &\leq 
        \int_{UV}\ikd(\Xcal\rel u,\Ycal)\d p_{U}(u)
        + 
        2\cdot\size{\Gbf}\cdot\ent(U)
      \end{align*}
    \item
      \label{p:slicingreduction}      
      Let $\Xcal\to U$ be a reduction, then
      \begin{align*}
      \ikd(\Xcal,\Ycal) 
      &\leq 
      \int_{U}\ikd(\Xcal\rel u,\Ycal)\d p_{U}(u)+
      \size{\Gbf}\cdot\ent(U)
      \end{align*}
    \item
      \label{p:slicingcofan}
      For a co-fan $\Xcal\to U\ot\Ycal$ holds
      \[
      \ikd(\Xcal,\Ycal)
      \leq
      \int_{U}\ikd(\Xcal\rel u,\Ycal\rel u)\d p_{U}(u)
      \]
    \end{enumerate}
  \end{corollary}
  
\subsection{Local estimate}\label{s:kolmogorov-local}
  Fix a complete diagram category $\Gbf$ and consider a
  $\Gbf$-configuration of sets $\Scal\in\Set\<\Gbf>$ with $S_{0}$
  being an initial set in $\Scal$. As discussed in Section
  \ref{s:disttypes-distributions-config}, the space of distributions
  on $\Scal$ could be identified with the space of distributions on the
  initial set
  \[
  \Delta\Scal\cong\Delta S_{0}
  \]
  Therefore, all $\Gbf$-configurations of probability spaces with the
  underlying configuration of sets equal to $\Scal$ are in one-to-one
  correspondence with the interior points of $\Delta S_{0}$.  The set
  $\Interior\Delta S_{0}$ consists of fully supported measures on the
  set $S_{0}$ and carries a total variation distance, which is just an
  $\ell^{1}$-distance with respect to the convex coordinates on the
  simplex $\Delta S_{0}$.  Our task presently is to compare the total
  variation distance with the Kolmogorov distance on the space of
  configurations with the fixed underlying configuration of sets.

  The upper bound on Kolmogorov distance, that we derive below, has two
  summands. One is linear in the total variation distance with the
  slope proportional to the $\log$-cardinality of $S_{0}$. The second one is
  super-linear in the total variation distance, but it does not depend
  on $\Scal$. So we have the following interesting observation: of
  course, the super-linear summand always dominates the linear one
  locally. However as the cardinality of $\Scal$ becomes large it is
  the linear summand that starts playing the main role.

\subsubsection{The estimate}
\label{p:kolmogorov-local-estimate}
  Suppose we are given a configuration of sets
  $\Scal=\set{S_{i},f_{ij}}\in\Set\<\Gbf>$ modeled on a complete diagram
  category $\Gbf$ with the initial set $S_{0}$.  
  We use once again the isomorphism
  \[ 
  \Delta\Scal\stackrel{\cong}{\to}\Delta S_{0}
  \]      
  that sends $p\in\Delta\Scal$ to its component in the initial space
  $p_{0}\in\Delta S_{0}$, while its inverse is given by
  $p=\set{(f_{0i})_{*}p_{0}}$.  For a pair of distributions
  $p_{0},q_{0}\in\Delta S_{0}$ denote by $|p_{0}-q_{0}|$ the
  total variation of the difference.

  For $\alpha\in[0,1]$ consider a binary probability space with the
  weight of one of the atoms equal to $\alpha$
  \[
  \Lambda_{\alpha}
  :=
  \big(
    \set{\square,\blacksquare};\,
    p_{\Lambda_{\alpha}}(\square)=1-\alpha,\,
    p_{\Lambda_{\alpha}}(\blacksquare)=\alpha
  \big)
  \]

  \newpage        
  \begin{proposition}{p:kolmogorovlocal}
    Let $\Scal=\set{S_{i},f_{ij}}\in\Set\<\Gbf>$ be a configuration of
    sets modeled on a complete diagram category $\Gbf$ with the
    initial set $S_{0}$. Let $p,q\in\Delta \Scal$ be two probability
    distributions.
    Denote $\Xcal:=(\Scal,p)$, $\Ycal:=(\Scal,q)$ and
    $\alpha=\frac12|p_{0}-q_{0}|_1$. Then
    \[
    \ikd(\Xcal,\Ycal)
    \leq
    2\cdot\size{\Gbf}\cdot\big(\alpha\cdot\ln|S_{0}|+\ent(\Lambda_{\alpha})\big)
    \]
  \end{proposition}
  To prove the local estimate we decompose both $p$ and $q$ into
  a convex combination of a common part $\hat p$ and rests $p^{+}$ and $q^+$. 
  The coupling between the common parts gives no contribution to the distance, and the worst possible estimate on the other parts is still enough to get the bound in the lemma, by using Corollary \ref{p:slicingcorollary} part (\ref{i:slicing2tents}). 
  Details of the proof can be found on page \pageref{p:kolmogorovlocal.rep}.
  
  In fact, the lower bound also holds, more specifically given a
  complete $\Gbf$-configuration of sets $\Scal$ and $p\in\Delta\Scal$
  there is a constant $C>0$ such that for any $q\in\Delta\Scal$ with
  $|p-q|_1<\!<1$ holds
  \[
  C\cdot\ent(\Lambda_{\alpha})\leq\ikd(\Xcal,\Ycal)
  \]      
  where $\Xcal=(\Scal,p)$, $\Ycal=(\Scal,q)$ and
  $\alpha=\frac12|p-q|_1$.
  We will not use this fact and therefore do not include a proof.

  Once the $\Gbf$-configuration of sets $\Scal$ is fixed, there is a map from $\Delta \Scal$ to $\Prob\<\Gbf>$.
  As can be seen from the discussion above, even though the map is continuous it is not Lipschitz.

	\section{Distance between types}\label{s:lagging}
	  As explained in Section~\ref{s:disttypes-types-config}, given a
  complete $\Gbf$-configuration $\Scal$ of sets and a rational
  distribution $\pi\in\Delta\Scal$ we construct a homogeneous
  configuration $T^{(n)}_{\pi}\Scal$, which is called the type of
  $\Scal$ over $\pi$. Our goal in this section is to estimate the
  Kolmogorov distance between two types over two different
  distributions $\pi_{1},\pi_{2}\in\Deltan\Scal$ in terms of the total
  variation distance $|\pi_{1}-\pi_{2}|_1$.

  For this purpose we use a ``lagging'' technique which is explained
  below.
 
\subsection{The lagging trick}\label{s:lagging-lagging}
  Let $\Lambda_\alpha$ be a binary probability space, 
  \[ 
  \Lambda_\alpha
  = 
  \big( \set{\square,\blacksquare};
        p_{\Lambda_\alpha}(\blacksquare)=\alpha 
  \big) 
  \] 
  and let $\Xcal=\set{(\un X_{i},p_{i});f_{ij}}$, $\Zcal=\set{(\un
    Z_{i},q_{i});g_{ij}}$ be two configurations modeled on a complete
  diagram category $\Gbf$ and included in a minimal two-fan
  \[ 
  \Lambda_\alpha \oot[\lambda] \Zcal \too[\rho] \Xcal 
  \]
  Recall that the left terminal vertex in this two-fan should be
  interpreted as a constant $\Gbf$-configuration
  $\Lambda_\alpha^{\Gbf}$.
  
  Assume further that the distribution $q$ on $\Zcal$ is rational with
  denominator $n\in\Nbb$, that is $q \in\Deltan\un\Zcal$.  It follows that
  $p$ and $p_{\Lambda_\alpha}$ are also rational with the same
  denominator $n$.

  We construct a \term{lagging two-fan} 
  \[\tageq{laggingtwofan}
    \Lcal
    :=
    \big(
      T^{\c{(1-\alpha) n}}(\Xcal\rel\square)
      \oot[l]
      T^{\c{n}}\Zcal
      \too[T\rho]
      T^{\c{n}}\Xcal 
    \big)
  \]
  as follows. 
  The right leg $T\rho$ of $\Lcal$ is induced by the right leg
  $\rho$ of the original two-fan.  The left leg
  \[ 
  l:T^{\c{n}}\Zcal
  \to 
  T^{\c{(1-\alpha) n}}(\Xcal\rel\square)
  \]
  is obtained by erasing symbols that reduce to $\blacksquare$ and
  applying $\rho$ to the remaining symbols.  The target space for the
  reduction $l$ is the true type of $\Xcal \rel \square$ which is
  ``lagging'' behind $T^{(n)} \Zcal$ by a factor of $(1-\alpha)$.
  More specifically, the reduction $l$ is constructed as follows.  Let
  $\lambda_{j}:Z_{j}\to \Lambda_\alpha$ be the components of the
  reduction $\lambda:\Zcal\to\Lambda_\alpha$.  

  Given $\bar z=(z_{i})_{i=1}^{n}\in T^{\c{n}}Z_{j} $ define the
  subset of indexes
  \[
    I_{\bar z}
    :=
    \set{i\st \lambda_{j}(z_{i})=\square }
  \]
  and define the $j^{\text{th}}$ component of $l$ by 
  \[
    l_{j}\big((z_{i})_{i=1}^{n}\big)
    :=
    (\rho(z_{i}))_{i\in I_{\bar z}}
  \]

  By equivariance each $l_{j}$ is a reduction of homogeneous spaces,
  since the inverse image of any point has the same
  cardinality. Moreover the reductions $l_{j}$ commute with the reductions in
  $T^{\c{n}}\Zcal$ as explained in Section \ref{s:disttypes-types-config}
  and therefore $l$ is a reduction of configurations.

  The next lemma uses the lagging two-fan to estimate the
  Kolmogorov distance between its terminal configurations.
 
  \begin{lemma}{p:lagging-kd}
    Let $\Xcal,\Zcal\in\prob\<\Gbf>$ be
    two configurations modeled on a complete diagram category $\Gbf$
    and included in a minimal two-fan
    \[
      \Lambda_\alpha
      \oot[\lambda]
      \Zcal
      \too[\rho]
      \Xcal
    \]
    where distribution on $\Zcal$ is rational with denominator
    $n\in\Nbb$.  Then
    \begin{align*}
      \ikd\Big(T^{\c{(1-\alpha) n}}&(\Xcal\rel\square)\,,\,
               T^{\c{n}}\Xcal\Big)
      \\
	&
	\leq
	n\cdot\size{G}\cdot
	\left[
	2\ent(\Lambda_\alpha)+
	\alpha\cdot\ln|X_{0}| 
	\right]
	+
	2 \cdot \size{\Gbf}\cdot|X_{0}|\cdot
	\ln(n+1)
      \\
      &= 
      n\cdot\size{G}\cdot
      \big[
        2\ent(\Lambda_\alpha)+
        \alpha \cdot\ln|X_{0}|
      \big]
      + \O\left(|X_{0}|\cdot\ln n \right)
    \end{align*}
  \end{lemma}

  The Lemma (and Proposition \ref{p:dist-types-complete} below) are
  closely related to the local estimate,
  Proposition \ref{p:kolmogorovlocal}. It is an immediate consequence
  of the Slicing Lemma, in particular
  Corollary \ref{p:slicingcorollary} part (\ref{p:slicingtwofan}) that
  \[          
  \ikd\Big(\Xcal\rel\square\,,\,
  \Xcal\Big)
  \leq
  \size{G}\cdot
  \left[
    2\ent(\Lambda_\alpha)+
    \alpha\cdot\ln|X_{0}| 
    \right]
  \]
  This is a tacit ingredient in the proof of the local estimate. 
  By the subadditivity of the Kolmogorov distance, 
  \[
  \ikd\Big((\Xcal\rel\square)^{\otimes n}\,,\,
  \Xcal^{\otimes n}\Big)
  \leq
  n\cdot
  \size{G}\cdot
  \left[
    2\ent(\Lambda_\alpha)+
    \alpha\cdot\ln|X_{0}| 
    \right]
  \]

  This bound is almost the estimate in Lemma \ref{p:lagging-kd},
  except Lemma \ref{p:lagging-kd} estimates the distance between types
  rather than tensor powers. We will soon see that tensor powers and
  types are very close in the Kolmogorov distance. However, for the
  purpose of the proof of Lemma \ref{p:lagging-kd}, it suffices to
  know that their entropies are close, an estimate that is provided by
  Corollary \ref{p:entropy-type}.

\medskip
\noindent{\it Proof} (of Lemma~\ref{p:lagging-kd}): 
  We will use the lagging two-fan constructed in
  Equation~(\ref{eq:laggingtwofan}), namely
  \[
    \Lcal
    :=
    \big(
    T^{\c{(1-\alpha) n}}(\Xcal\rel\square)
    \oot[l]
    T^{\c{n}}\Zcal
    \too[T\rho]
    T^{\c{n}}\Xcal 
    \big)
  \]
  as a coupling to estimate the Kolmogorov distance
  \[
    \ikd\big(T^{\c{(1-\alpha) n}}(\Xcal\rel\square)\,,\,
    T^{\c{n}}\Xcal\big)
    \leq
    \kd(\Lcal)
  \]
	
  Recall that by Corollary \ref{p:entropy-type} for a probability
  space $X$ with a rational distribution we have
  \[  
    n \cdot\ent(X) - |X| \cdot \ln (n+1) 
    \leq 
    \ent(T^{\c n} X) \leq n \cdot\ent(X)
  \]  
  
  Thus we can estimate $\kd(\Lcal)$ as follows
  \begin{align*}
    \kd(\Lcal)
    &=
    \sum_{i}
    \Big[
      \big(
      \ent(T^{\c n}Z_{i}) 
      - 
      \ent(T^{\c{n}} X_{i})
      \big)
      \\
      &\qquad \qquad +  
      \big(
      \ent(T^{\c n}Z_{i}) 
      - 
      \ent(T^{((1-\alpha)n)} (X_{i}\rel\square))
      \big)
    \Big]
    \\
    &\leq
    n\cdot
    \sum_{i}
    \Big[
      \big(
      \ent(Z_{i}) 
      - 
      \ent(X_{i})
      \big)
      +
      \big(
      \ent(Z_{i}) 
      - 
      (1-\alpha)\ent(X_{i}\rel\square)
      \big)
    \Big] 
    \\
    &\quad + 
    2 \cdot \size{\Gbf}\cdot |X_0|\cdot\ln(n+1) 
  \end{align*}
     
  By minimality of the original two-fan and Shannon inequality
  (\ref{eq:shannonineq}) we have a bound
  \[
    \ent(Z_{i}) - \ent(X_{i}) 
    \leq
    \ent(\Lambda_\alpha)
  \]
	
  The second part in the sum can be estimated using relation
  (\ref{eq:relentinteg}) as follows
  \begin{align*}
    \ent(Z_{i}) - (1-\alpha)\ent(X_{i}\rel\square)
    &=
    \ent(\Lambda_\alpha) + 
    \ent(X_{i}\rel\Lambda_\alpha) - 
    (1-\alpha)\ent(X_{i}\rel\square)
    \\
    &=
    \ent(\Lambda_\alpha) + 
    (1-\alpha)\ent(X_{i}\rel\square) + 
    \alpha\Ent(X_{i}\rel\blacksquare) -
    \\
    &\quad-
    (1-\alpha)\ent(X_{i}\rel\square)
    \\
    &\leq
    \ent(\Lambda_\alpha)
    +
    \alpha\cdot\ln|X_{i}|
  \end{align*}
  Combining all of the above we obtain the estimate in the conclusion
  of the lemma.  
\eproof

\subsection{Distance between types}
  In this section we use the lagging trick as described above to
  estimate the distance between types over two different distributions
  in $\Delta \Scal$ where $\Scal$ is a complete configuration of sets.

  \begin{proposition}{p:dist-types-complete}
    Suppose $\Scal$ is a complete $\Gbf$-configuration of sets with initial
    set $S_{0}$. 
    Suppose $p, q \in\Deltan \Scal$ and let
    $\alpha=\frac12|p_0-q_0|_1$. Then
    \begin{align*}  
      \ikd(T^{\c n}_{p}\Scal,T^{\c n}_{q}\Scal)
      &\leq 
      2n\cdot\size{\Gbf}\cdot
      \left[
        \alpha\cdot\ln |S_{0}|
        + 
        2\ent(\Lambda_{\alpha})
      \right]
      + 
      4 \size{\Gbf}\cdot|X_{0}|\cdot\ln(n+1)
      \\
      &=
      2n\cdot\size{\Gbf}\cdot
      \left[
        \alpha\cdot\ln |S_{0}|
        + 
        2\ent(\Lambda_{\alpha})
        \right]
      + 
      \O(|X_{0}|\cdot\ln n)
    \end{align*} 
  \end{proposition}

  As in the local estimate, the idea of the proof is to write $p$ and
  $q$ as a convex combination of a common distribution $\hat p$ and
  ``small amounts'' of $p^{+}$ and $q^{+}$, respectively. Then we use
  the lagging trick to estimate distances between types over $p$ and
  $\hat p$, as well as between types over $q$ and $\hat p$.
  We now present details of the proof.
  
\medskip
\noindent\textit{Proof} (of Proposition~\ref{p:dist-types-complete}):
  	Recall that for a complete configuration $\Scal$ with initial set
  	$S_{0}$ we have 
  	\[\tageq{distr-complete-again}
  	\Delta\Scal\cong\Delta S_{0}
  	\]

  	Our goal now is to write $p$ and $q$ as the convex combination of
  	three other distributions $\hat p$, $p^{+}$ and $q^{+}$ as in
  	\begin{align*}
  	p
  	&=
  	(1-\alpha)\cdot\hat p + \alpha\cdot p^{+}
  	\\
  	q
  	&=
  	(1-\alpha)\cdot\hat p + \alpha\cdot q^{+}
  	\end{align*}
  	
  	We could do it the following way.  Let
  	$\alpha:=\frac12|p_{0}-q_{0}|_1$. If $\alpha=1$ then the proposition
  	follows trivially by constructing a tensor-product fan, so
  	from now on we assume that $\alpha<1$.  Define three probability
  	distributions $\hat p_{0}$, $p_{0}^{+}$ and $q_{0}^{+}$ on $S_{0}$
  	by setting for every $x\in S_{0}$
  	\begin{align*}
  	\hat p_{0}(x) 
  	&:= \frac1{1-\alpha}
  	\min\set{p_{0}(x),q_{0}(x)}
  	\\ 
  	p_{0}^{+} 
  	&:=
  	\frac{1}{\alpha}\big(p_{0}-(1-\alpha)\hat p_{0}\big)
  	\\ 
  	q_{0}^{+} 
  	&:= 
  	\frac{1}{\alpha}\big(q_{0}-(1-\alpha)\hat p_{0}\big)
  	\end{align*}
  	
  	Denote by $\hat p,p^{+},q^{+}\in\Delta\Scal$ the distributions
        corresponding to $\hat p_{0},p_{0}^{+},q_{0}^{+}\in\Delta
        S_{0}$ under the affine
        isomorphism~(\ref{eq:distr-complete-again}). Thus we have
  	\begin{align*}
  	p&=(1-\alpha)\hat p+\alpha\cdot p^{+}\\
  	q&=(1-\alpha)\hat p+\alpha\cdot q^{+}
  	\end{align*}
  	
  	Now we construct a pair of two-fans of
  	$\Gbf$-configurations 
  	\begin{align*}
  	\tageq{twofanfortypes}
  	\Lambda_{\alpha}\ot\tilde\Xcal\to\Xcal\\
  	\Lambda_{\alpha}\ot\tilde\Ycal\to\Ycal
  	\end{align*}
  	by setting 
  	\begin{align*}
  	\Xcal
  	&:=
  	(\Scal,p)
  	\\
  	\Ycal
  	&:=
  	(\Scal,q)
  	\\
  	\tilde X_{i}
  	&:=
  	\Big(S_{i}\times\un\Lambda_{\alpha};\;
  	\tilde p_{i}(s,\square)=(1-\alpha)\hat p_{i}(s),\,
  	\tilde p_{i}(s,\blacksquare)=\alpha\cdot p^{+}_{i}(s)
  	\Big)
  	\\
  	\tilde Y_{i}
  	&:=
  	\Big(S_{i}\times\un\Lambda_{\alpha};\;
  	\tilde q_{i}(s,\square)=(1-\alpha)\hat p_{i}(s),\,
  	\tilde q_{i}(s,\blacksquare)=\alpha\cdot q^{+}_{i}(s)
  	\Big)
  	\\
  	\end{align*}
  	and
  	\begin{align*}
  	\tilde\Xcal
  	&:=
  	\set{\tilde X_{i};\,f_{ij}\times\Id}\\
  	\tilde\Ycal
  	&:=
  	\set{\tilde Y_{i};\,f_{ij}\times\Id}\\
  	\end{align*}
  	
  	The reductions in~(\ref{eq:twofanfortypes}) are given by coordinate
  	projections. We have the following isomorphisms
  	\begin{align*}
  	\Xcal\rel\square
  	&\cong
  	\Ycal\rel\square
  	\cong
  	(\Scal,\hat p)
  	\end{align*}
  	
  	To estimate the distance between types we now apply
  	Lemma~\ref{p:lagging-kd} to the fans
        in~(\ref{eq:twofanfortypes}) 
  	\begin{align*}
  	\ikd(T^{(n)}_{p}\Scal,T^{(n)}_{q}\Scal)
  	&=
  	\ikd(T^{(n)}\Xcal,T^{(n)}\Ycal)
  	\\
  	&\leq
  	\ikd\big(T^{(n)}\Xcal,T^{((1-\alpha)n)}(\Xcal\rel\square)\big)  
  	+
  	\ikd\big(T^{((1-\alpha)n)}(\Ycal\rel\square),T^{(n)}\Ycal\big)  
  	\\
        &\leq
        2n\cdot\size{\Gbf}\cdot
        \left[
          \alpha\cdot\ln |S_{0}|
          + 
          2\ent(\Lambda_{\alpha})
        \right]
        + 
        4 \size{\Gbf}\cdot|X_{0}|\cdot\ln(n+1)
  	\end{align*}	
\eproof

  The reason for the similarity between the local estimate and the
  distance estimate between types will become clear in the next
  section, when we establish the asymptotic equivalence between the
  Bernoulli sequence of probability spaces and sequence of types over
  rational distributions approximating the true distribution.

	\section[AEP for configurations]{Asymptotic equipartition 
		property for configurations}
	\label{s:ac-aep-conf}
	  Below we prove that any Bernoulli sequence can be approximated by a
  sequence of homogeneous configurations.  This is essentially the
  \term{Asymptotic Equipartition Theorem for configurations}.
 
  \begin{theorem}{p:aep-complete}
    Suppose $\Xcal\in\prob\<\Gbf>$ is a complete configuration of
    probab ility spaces.  Then there exists a sequence
    $\bar\Hcal=(\Hcal_{n})_{n=0}^{\infty}$ of homogeneous
    configurations of the same c ombinatorial type as $\Xcal$ such
    that
    \[
      \frac{1}{n}\ikd (\Xcal^{\otimes n},\Hcal_{n}) 
      = 
      \O\left( \sqrt{\frac{\ln^3 n}{n}} \right)
    \]

  More precisely, the sequence $\bar\Hcal$ may be chosen such that for
  all $n \geq |X_0|$
  \[\tageq{quantaep} 
  \frac{1}{n} 
  \ikd (\Xcal^{\otimes n},\Hcal_{n}) 
  \leq 
  C(|X_0|,\size{\Gbf}) \cdot \sqrt{\frac{\ln^3 n}{n}} 
  \]
  where $C(|X_0|, \size{\Gbf})$ is a constant only depending on
  $|X_0|$ and $\size{\Gbf}$.
\end{theorem}

\begin{Proof}
  Denote by $\Scal = \un \Xcal$ the underlying configuration of sets
  and by $p_\Xcal$ the true distribution on $\Scal$, such that
  \[
  \Xcal = (\Scal, p_X)
  \]

  We will construct the approximating homogeneous sequence by taking
  types over rational approximations of $p_\Xcal$ in $\Delta\Scal$, that
  converge sufficiently fast to the true distribution $p_{\Xcal}$.

  More specifically, we select rational distributions $p_n \in \Deltan
  \Scal$ such that
  \[
  |p_n - p_\Xcal|_{1} \leq \frac{|S_{0}|}{n}
  \]

  As homogeneous spaces $\Hcal_n$ we set $\Hcal_n = T_{p_{n}}^{\c
    n}\Scal$.  We will show that the Kolmogorov distance between
  $\Hcal_n$ and $\Xcal^{\otimes n}$ satisfies the required estimate
  (\ref{eq:quantaep}).

  First we apply slicing along the empirical two-fan 
  \[
  \Rcal_{n}(\Xcal)
  =
  \left(
    \Xcal^{\otimes n}
    \ot
    \tilde\Xcal^{(n)}
    \to
    (\Delta\Scal,\tau_{n})^{\Gbf} 
  \right)
  \]
  defined in Section~\ref{s:disttypes-types-config},
  Equation~(\ref{eq:empirical-2fan}) on
  page~\pageref{eq:empirical-2fan}.

  For the estimate below we use the fact that
  \[
    \ent(\Delta\Scal,\tau_{n})
    \leq
    \ln|\Deltan\Scal|
    \leq |S_0|\cdot\ln (n+1)
  \]

  By slicing (see
  Corollary~\ref{p:slicingcorollary}(\ref{p:slicingtwofan})) along the
  empirical two-fan we have
  \begin{align*}
    \ikd(T^{(n)}_{p_{n}}\Scal,\Xcal^{\otimes n})
    &\leq
    2 \cdot \size{\Gbf}\cdot\ent(\Delta \Scal,\tau_{n})
    +
    \int_{\Delta \Scal}
    \ikd(T^{(n)}_{p_{n}}\Scal,T^{(n)}_{\pi}\Scal)
    \d\tau_{n}(\pi) 
    \\
    &\leq
    2 \cdot \size{\Gbf} \cdot|S_0|\cdot \ln (n + 1)
    +
    \int_{\Delta \Scal}
    \ikd(T^{(n)}_{p_{n}}\Scal,T^{(n)}_{\pi}\Scal)
    \d\tau_{n}(\pi) 
  \end{align*}
  To estimate the integral we split the domain into a small divergence
  ball $B_{\epsilon_n} = B_{\epsilon_n}(p_\Xcal)$ around the ``true''
  distribution and its complement
  \begin{align*}
    \int_{\Delta \Scal}
    \ikd(T^{(n)}_{p_{n}}\Scal,T^{(n)}_{\pi}\Scal)
    \d\tau_{n}(\pi) 
    &=
    \int_{\Delta \Scal\setminus B_{\epsilon_{n}}}
    \ikd(T^{(n)}_{p_{n}}\Scal,T^{(n)}_{\pi}\Scal)
    \d\tau_{n}(\pi)\\
    &\quad +\tageq{aepinball}
    \int_{B_{\epsilon_{n}}}
    \ikd(T^{(n)}_{p_{n}}\Scal,T^{(n)}_{\pi}\Scal)
    \d\tau_{n}(\pi)
  \end{align*}
  and we set the radius $\epsilon_n$ equal to
  \[
  \epsilon_n 
  := 
  (|S_0| + 1) \frac{\ln (n+1)}{n}
  \]
  
  To estimate the first integral on the right-hand side of equality (\ref{eq:aepinball}) note that the distance between two
  types over the same configuration of sets can always be crudely estimated
  by 
  \[
  2 \cdot \ln |S_0| \cdot \size{\Gbf} \cdot n
  \]
  Moreover, by Sanov's theorem, Theorem \ref{p:entropydivergence}, we can estimate the empirical measure of the complement of the divergence ball
  \[
  \tau_{n}(\Delta\Scal\setminus
  B_{\epsilon_{n}}) \leq \ebf^{-n\cdot\epsilon_{n}+|S_0| \cdot \ln (n+1)} \leq \frac{1}{n}
  \] 
  where we used the definition of $\epsilon_n$ to conclude the last inequality.
  Therefore we obtain
  \begin{align*}
    \int_{\Delta \Scal\setminus B_{\epsilon_{n}}}
    \ikd(T^{(n)}_{p_{n}}\Scal,T^{(n)}_{\pi}\Scal) \d\tau_{n}(\pi) 
    &\leq
    2 \cdot\ln |S_0| \cdot \size{\Gbf} \cdot n \cdot
    \tau_{n}(\Delta\Scal\setminus
    B_{\epsilon_{n}})
    \\
    &\leq 
    2 \cdot\ln |S_0| \cdot \size{\Gbf}
  \end{align*}
  
  Define
  \[
  \alpha_n 
  = 
  \frac{|S_0|}{n} + \sqrt{2 \epsilon_n }
  \]
  if the right-hand side is smaller than $1$ and set $\alpha_n=1$
  otherwise.  Then every $\pi \in B_{\epsilon_n}(p_\Xcal)$ satisfies
  $|p_n - \pi| \leq 2 \alpha_n$ by Pinsker's inequality (Lemma
  \ref{p:entropydivergence}, (\ref{p:entropydivergence1})), and the
  triangle inequality.  Consequently, by the estimate on the distance
  between types in Proposition \ref{p:dist-types-complete}
  
  \begin{align*}
  \int_{B_{\epsilon_{n}}}
  &\ikd(T^{(n)}_{p_{n}}\Scal,T^{(n)}_{\pi}\Scal)
  \d\tau_{n}(\pi)
  \\
  &\leq
  2n \cdot
  \size{\Gbf}
  \cdot
  \left( \alpha_n \ln |S_0| + 2 \ent(\Lambda_{\alpha_n}) \right)
    + 4\cdot \size{\Gbf} \cdot |S_{0}| \cdot \ln(n+1)
  \end{align*}
  
  Using the definition of $\alpha_n$ and $\epsilon_n$ we find that 
  \[
    \int_{B_{\epsilon_{n}}}
    \ikd(T^{(n)}_{p_{n}}\Scal,T^{(n)}_{\pi}\Scal)
    \d\tau_{n}(\pi) 
    = 
    \O\left(\sqrt{n \cdot \ln^{3} n} \right)
  \]
  and hence combining the above estimates
  \[
    \frac{1}{n}\ikd(T^{(n)}_{p_{n}}\Scal,\Xcal^{\otimes n}) 
    = 
    \O\left( \sqrt{\frac{\ln^3 n }{n}} \right).
  \]
  A more precise check shows that for $n \geq |S_0|$, the constants
  appearing in $\O$ only depend on $|S_0|$ and $\size{\Gbf}$.
\end{Proof}

\bigskip

  It is worth noting that each type considered as a subspace of the
  tensor power takes up only small probability. In fact its
  probability converges to zero with growing $n$. But as the
  calculation above shows, most (in terms of probability) of the
  configuration $\Xcal^{\otimes n}$ consists of \emph{polynomially
  many} types that are $\ikd$-similar to each other. Relative to the
  exponential growth of sizes of all the parts, ``polynomially many''
  is as good as one.  This is the difference with the setup used in
  Gromov's \cite{Gromov-Search-2012}

	\section{Extensions}
	\label{s:extensions}
	  In the introduction we have already emphasized the close
  relationship between relative entropic sets and
  Information-Optimization problems.  There, our definitions were
  restricted to extensions of two-fans to full configurations
  corresponding to three random variables.  We will now generalize
  these definitions, and make the relationship between relative
  en\-tro\-pic sets and Information-Optimization problems explicit.

  Further we will prove the Extension Lemma and use it to show that
  the relative entropic set associated to a full configuration depends
  continuously on the configuration.

\subsection{Information-Optimization and the relative entropic set}

  In Section \ref{s:config-restrictions-fullfull} we introduced (for
  $k \leq l$) the restriction operator 
  \[
  R_{k,l}^*: \prob\<\Lambdabf_l> 
  \to 
  \prob\<\Lambdabf_k>
  \]
  as follows.  For a minimal full configuration
  $\Ycal=\langle Y_{i}\rangle_{i=1}^{l}$ we denote by
  \[
    R_{k,l}^*\Ycal
    =
    \langle Y_{i}\rangle_{i\in \{1, \ldots, k\}}
  \]
  the restriction of $\Ycal$ to a minimal full configuration generated
  by $Y_{i}$, $i\in \{1, \dots, k\}$.

  We call a minimal configuration $\Ycal \in \Prob\<\Lambda_k>$ an
  $l$-extension of a configuration $\Xcal$ if
  \[
    R_{k,l}^* \Ycal = \Xcal
  \]
  and we denote the class of all $l$-extensions by $\Ext_l(\Xcal)$.

  Recall that for a full configuration $\Ycal \in \Prob\<\Lambda_l>$,
  we record the entropies of all its probability spaces in a vector in
  $\Rbb^{2^{\set{1,\dots,l}}\setminus\set{\emptyset}}$ that we denote by
  \[
  \Ent_*(\Ycal) := \big(\Ent(Y_I)\big)_{I \in 2^{\set{1, \dots,
        l}}\setminus\set{\emptyset}}
  \]
  The entries in this vector are all nonnegative. To simplify
  notations we set 
  \[
  \Ebb_{l}:=
  (\Rbb^{2^{\set{1,\ldots,l}}\setminus\set{\emptyset}},|\,\cdot\,|_{1})
  \]
   and denote by $\Ebb_{l}^{*}$ its dual vector-space.

  As in the introduction, we introduce the \term{unstabilized relative
    entropic set}
  \[
  \res_l(\Xcal) 
  := 
  \set{\Ent_*(\Ycal) \st \Ycal \in \Ext_l(\Xcal)}
  \]
  By the additivity property of the entropy with respect to tensor
  powers, there is the inclusion
  \[\tageq{superadditivity-res}
  \res_l(\Xcal^{\otimes m}) 
  + 
  \res_l(\Xcal^{\otimes n}) 
  \subset 
  \res_l(\Xcal^{\otimes (m + n)})
  \]
  where the sum on the left hand side is the Minkowski sum. 
  This
  allows us to define the limit
  \[
  \lim_{n\to \infty} \frac{1}{n} \res_l(\Xcal^{\otimes n}) 
  := \bigcup_{n \in \Nbb} \frac{1}{n} \res_l (\Xcal^{\otimes n})
  \]
  and we define the \term{stabilized relative entropic set} by
  \[
  \sres_l(\Xcal) 
  := 
  \closure 
    \left( 
      \lim_{n \to \infty} 
        \frac{1}{n} \res_l(\Xcal^{\otimes n}) 
    \right)
  \]
    which is a closed convex subset of $\Ebb_{l}$ by
  property~(\ref{eq:superadditivity-res}).

  For a vector $c \in \Ebb_{l}^*$, we define the
  \term{Information-Optimization problem}
  \[
  \IO_{c}(\Xcal) 
  := 
  \inf_{\Ycal \in \Ext_l(\Xcal)} 
    \big\langle\, c\, ,\, \Ent_*(\Ycal)\,\big\rangle
  = 
  \inf_{\Ycal \in \Ext_l(\Xcal)} 
    \sum_{I \subset \{1, \dots, l\}} c_I \cdot\Ent (Y_I)
  \]
  where $c_{I}$'s are the coordinates of the vector $c$ with respect
  to the basis in $\Ebb^{*}_{l}$ dual to the standard basis in
  $\Ebb_{l}$.  Note that, equation~(\ref{eq:superadditivity-res})
  implies that the sequence
  \[
  n \mapsto \IO_{c}{(\Xcal^{\otimes n})}
  \]
  is subadditive.  Hence, the limit
  \[
  \lim_{n \to \infty} \frac{1}{n} \IO_{c}(\Xcal^{\otimes n})
  \]
  always exists (but may be equal to $-\infty$). If for all $n \in \Nbb$, 
  \[
  \frac{1}{n} \IO_{c}(\Xcal^{\otimes n}) = \IO_{c}(\Xcal)
  \]
  we call the optimization problem associated
  to $c$ \term[stable optimization problem]{stable}.  In general, we
  define the stabilized optimization problem
  \[
  \IO^s_{c}(\Xcal) := 
  \lim_{n \to \infty} \frac{1}{n} \IO_{c}(\Xcal_n^{\otimes n}).
  \]
  
  As the stabilized relative entropic set is convex, it is the
  intersection of half-spaces that are defined by linear inequalities
  on entropies
  \[
  \sres_k(\Xcal) 
  = 
  \bigcap_{c \in \Ebb_{l}^{*}} 
  \set{ x \in \Ebb_{l} \, 
    \st 
    \, \langle c,x\rangle \geq \IO^s_{c}(\Xcal) }
  \]

  In other words, the stabilized information optimization problems,
  that occur so often in practice, identify supporting hyper-planes of
  the convex set.  The solution of all such linear problems determine
  the shape of the relative entropic set and vice versa.

\subsection{The entropic set and the entropic cone}
  The definitions of relative entropic sets are motivated by the more
  classical notion of the entropic
  cone, which we will briefly discuss now.  For $l \in \Nbb$, the
  entropic set is defined as
  \[
  \res_l 
  := 
  \set{ \ent_*(\Ycal) 
    \st 
    \Ycal \in \prob\<\Lambdabf_l>,\,\Ycal \text{ is minimal} 
  }
  \]
  Its closure is usually referred to as the entropic cone
  \[
  \sres_l 
  := 
  \closure(\res_l)
  \]
  Indeed, the entropic cone $\sres_l$ is a closed, convex cone in
  $\Rbb^{2^l-1}$ \cite{Yeung-First-2012}. For $l \leq 3$, the entropic
  cone $\sres_l$ is polyhedral and completely described by Shannon
  inequalities. However, for $l \geq 4$, the situation is much more
  complicated. It is known that $\sres_l$ is not polyhedral for $l
  \geq 4$ \cite{Matus-Infinitely-2007}. The shape of the entropic cone
  is not known as of the time of writing this article. It is an
  important open problem in information theory to find tight bounds on
  the entropic cone for $l \geq 4$. We hope that the techniques
  developed in this article will eventually lend itself to finding a
  useful characterization.

  In fact, the entropic cone can be considered as the relative
  entropic set of an empty configuration     
  $\emptyconfig\in\prob\<\emptycat>$, that corresponds to the
  empty diagram category $\emptycat=\Lambdabf_{0}$
  \[
  \sres_{l}=\sres_{l}(\emptyconfig)
  \]

  For a diagram category $\Gbf$ let us denote by
  $\set{\bullet}=\set{\bullet}^{\Gbf}$ the constant
  $\Gbf$-configuration of one-point spaces.  Given an $l$-extension
  $\Ycal\in\prob\<\Lambdabf_{l}>$ of $\set{\bullet}^{\Lambdabf_{k}}$ the restriction to
  the last $l-k$ terminal spaces induces a linear isomorphism
  \[\tageq{entcone-and-entset}
  \sres_{l}(\set{\bullet}^{\Lambdabf_{k}})
  \cong
  \sres_{l-k}
  \]

\subsection{Extension lemma}

  The Lipschitz continuity of relative entropic sets will follow from
  the following important proposition, which we will refer to as the
  Extension Lemma.

\begin{proposition}{p:extensionlemma}{\rm(Extension Lemma)}
  Let $k,l\in\Nbb$, $k\leq l$ and let $\Xcal, \Xcal' \in
  \prob\<\Lambdabf_k>$ be minimal full configurations.  For
  every $\Ycal \in \Ext_l \Xcal$ there exists a $\Ycal' \in \Ext_l
  \Xcal'$ such that
  \[
  \ikd(\Ycal', \Ycal ) 
  \leq 2^{l-k} 
  \ikd (\Xcal', \Xcal)
  \]
\end{proposition}

The key behind the proof of the Extension Lemma, is that there is a full configuration $\Zcal$ that extends both $\Ycal$ and the optimal coupling between $\Xcal$ and $\Xcal'$. The configuration $\Ycal'$ can be chosen to be the restriction of $\Zcal$ to the full configuration generated by $\Xcal'$ and the terminal spaces in $\Ycal$ which are not in $\Xcal$. The estimate directly follows from Shannon inequalities.  
We present details at page \pageref{p:extensionlemma.rep}.

It follows immediately from the Extension Lemma and the Lipschitz
property of the entropy function $\Ent_*$ that asymptotically
equivalent configurations have the same solutions to all
  Information-Optimization problems and, consequently, they have the
same stabilized relative entropic set.

  In fact, we have a much stronger statement.  Both the unstabilized
  and stabilized relative entropic sets have a Lipschitz dependence on
  the configuration, if the distance between sets is measured by the
  Hausdorff distance.

  Let us endow the collection of subsets of $\Ebb_{l}$ with the
  Hausdorff metric with respect to the $\ell_1$-distance. For two
  subsets $S_1, S_2$ of $\Ebb_{l}$, define the Hausdorff distance
  between them by
  \[
  \d_H \left( S_1, S_2 \right) 
  = 
  \inf\set{ \epsilon > 0 
            \st 
            S_1 \subset S_2 + B_\epsilon \text{ and } S_2 \subset
            S_1 + B_\epsilon 
          }
  \]
  where $B_\epsilon$ is the $\ell_1$-ball of size $\epsilon$ around
  the origin in $\Ebb_l$.

  In fact, at this point the Hausdorff distance is only an extended
  pseudo-metric, in the sense, that it may take infinite values and it
  may vanish on pairs of non-identical points.

  Suppose now that we are given two minimal full configurations
  $\Xcal, \Xcal' \in \prob\<\Lambdabf_k>$, and suppose a point $y \in
  \Ebb_l$ lies in the unstabilized relative
  entropic set of $\Xcal$, that is
  \[
  y \in \res_l(\Xcal )
  \]
  This means that there is an extension $\Ycal \in \Ext_l(\Xcal)$ such that 
  \[
  \ent_*(\Ycal) = y
  \]
  By the Extension Lemma, there exists a configuration $\Ycal' \in
  \Ext_l(\Xcal')$ such that
  \[
  \ikd(\Ycal, \Ycal') 
  \leq 
  2^{l-k} \ikd (\Xcal, \Xcal')
  \]
  and by the $1$-Lipschitz property of the entropy function the point
  $y' := \ent_*(\Ycal')$ is close to the point $y$, that is
  \[
  | y - y' |_{1} 
  = 
  | \ent_*(\Ycal) - \ent_*(\Ycal')|_{1} \leq 2^{l-k} \ikd(\Xcal, \Xcal')
  \]
  We have thus obtained the following corollary to the Extension Lemma.

  \begin{corollary}{p:unstabilized-rel-ent-set-Lipschitz}
    Let $k \in \Nbb$ and $\Xcal, \Xcal' \in \prob\<\Lambdabf_k>$.
    Then the Hausdorff distance between their unstabilized relative
    entropic sets satisfies the following Lipschitz estimate
    \[
    \d_H\big( \res_{l}(\Xcal), \res_l(\Xcal') \big) 
    \leq 2^{l-k} \ikd(\Xcal, \Xcal')
    \]
  \end{corollary}

  Note that in particular, the distance between unstabilized relative
  entropic sets is always finite and
  \[
  \d_H\big( 
        \res_l(\Xcal), 
        \res_l(\set{\bullet}^{\Lambdabf_{k}} )
      \big) 
  \leq 
  2^{l-k} \ikd(\Xcal, \set{\bullet}^{\Lambdabf_{k}}) 
  = 
  2^{l-k}|\ent_*(\Xcal)|_1.
  \]
  Let us denote by $\convconend$ the metric space of closed convex
  sets $K$ in $\Ebb_{l}$ such that
  \[
  d_H( K, \sres_l(\set{\bullet}^{\Lambdabf_{k}})  ) < \infty
  \]
  endowed with the Hausdorff distance.

  \begin{theorem}{p:stabilized-rel-ent-set-Lipschitz}
    Let $k \in \Nbb$ and $\Xcal, \Xcal' \in \prob\<\Lambdabf_k>$.
    Then for all $l \in \Nbb$, the Hausdorff distance between their
    stabilized relative entropic sets satisfies the Lipschitz estimate
    \[
    \d_H\big( \sres_l(\Xcal), \sres_l(\Xcal') \big) 
    \leq 
    2^{l-k} \aikd(\Xcal, \Xcal')
    \]
    In other words, the map $\sres_l$ from minimal full
    configurations in $\prob\<\Lambdabf_k>$ to $\convconend$ is
    $2^{l-k}$-Lipschitz.
\end{theorem}

  Finally, as a primer to Section \ref{s:tropical}, note that for any
  set $K \in \convconend$ the sequence
  \[
  n \mapsto \frac{1}{n} K
  \]
  converges in the Hausdorff distance to
  $\sres_l(\set{\bullet}^{\Lambdabf_{k}})$. The set $K\subset
  \Ebb_{l}$ can be viewed as a metric space itself, by just restricting
  the $\ell^1$-metric to it. The above convergence can then be
  expressed by saying that the asymptotic cone of $K$ equals
  $\sres_l(\set{\bullet}^{\Lambdabf_{k}})$ and is isomorphic to $\sres_{l-k}$.

	\section{Mixtures}
	\label{s:mixtures}
	  Mixtures provide some technical tools, which we will use in
  Section~\ref{s:tropical}. The input data for the mixture operation
  is a family of $\Gbf$-con\-fi\-gu\-ra\-tions, parametrized by a
  probability space. As result one obtains another
  $\Gbf$-con\-fi\-gu\-ra\-tion with the pre-specified
  conditionals. One particular instance of a mixture is when one mixes
  two configurations $\Xcal$ and $\set{\bullet}^{\Gbf}$, the latter
    being a constant $\Gbf$-configuration of one-point probability spaces. This
  operation will be used as a substitute for taking radicals
  ``$\Xcal^{\otimes(1/n)}$'' in Section~\ref{s:tropical} below.

\subsection{Definition and elementary properties}
  Let $\Gbf$ be a complete diagram category and $\Theta$ be a
  probability space. Let $\set{\Xcal_{\theta}}_{\theta\in\un\Theta}$
  be a family of $\Gbf$-configurations parametrized by $\Theta$. The
  \term{mixture} of the family $\set{\Xcal_{\theta}}$ is the reduction
  \[      
  \mix\set{\Xcal_{\theta}}=
  \left(
    \Ycal
    \too
    \Theta^{\Gbf}
  \right)
  \]      
  such that
  \[\tageq{mixture}
  \Ycal\rel\theta\cong\Xcal_{\theta}
  \]

  The mixture exists and is uniquely defined by
  property~(\ref{eq:mixture}) up to an isomorphism which is identity
  on $\Theta^{\Gbf}$.

  We denote the top configuration of the mixture
  \[
  \Ycal=\bigoplus_{\theta\in\Theta}\Xcal_{\theta}
  \]
  and also call it the mixture of the family $\set{\Xcal_{\theta}}$.

  When 
  \[
  \Theta=\Lambda_{\alpha}
  =
  \big(\set{\square,\blacksquare};
  p(\blacksquare)=\alpha
  \big) 
  \]
  is a binary space we write simply
  \[
  \Xcal_{\blacksquare}\oplus_{\Lambda_{\alpha}}\Xcal_{\square}
  \]
  for the mixture. The configuration subindexed by the $\blacksquare$
  will always be the first summand.

  The entropy of the mixture can be evaluated by the following formula
  \[
  \ent_{*}\left(\bigoplus_{\theta\in\Theta}\Xcal_{\theta}\right)
  =
  \int_{\Theta}\ent_{*}(\Xcal_{\theta})\d
  p(\theta) + \ent_{*}(\Theta^{\Gbf})
  \]
  Mixtures satisfy the distributive law with respect to the tensor
  product
  \begin{align*}
    \mix(\set{\Xcal_{\theta}}_{\theta\in\Theta})
    \otimes
    \mix(\set{\Ycal_{\theta'}}_{\theta'\in\Theta'})
    &\cong
    \mix(\set{\Xcal_{\theta}\otimes\Ycal_{\theta'}}
    _{(\theta,\theta')\in\Theta\otimes\Theta'})
    \\
    \left(\bigoplus_{\theta\in\Theta}\Xcal_{\theta}\right)
    \otimes
    \left(\bigoplus_{\theta'\in\Theta'}\Ycal_{\theta'}\right)
    &\cong
    \bigoplus_{(\theta,\theta')\in\Theta\otimes\Theta'}(\Xcal_{\theta}\otimes\Ycal_{\theta'})
  \end{align*}
  
\subsection{The distance estimates}
  Recall that for a diagram category $\Gbf$ we denote by
  $\set{\bullet}=\set{\bullet}^{\Gbf}$ the constant
  $\Gbf$-configuration of one-point spaces. 
  
  The mixture of a $\Gbf$-configuration with $\set{\bullet}^{\Gbf}$
  may serve as an ersatz of taking radicals of the configuration.  The
  following lemma provides a justification of this by some distance
  estimates related to mixtures and will be used in
  Section~\ref{s:tropical}.

  \begin{lemma}{p:mixtures-dist1}
    Let $\Gbf$ be a complete diagram category and
    $\Xcal,\Ycal\in\prob\<\Gbf>$. Then
    \begin{enumerate}
    \item
      $\displaystyle
      \aikd(\Xcal,\Xcal^{\otimes n}\oplus_{\Lambda_{1/n}}\set{\bullet})
      \leq
      \ent(\Lambda_{1/n})
      $
    \item
      $\displaystyle
      \aikd\big(\Xcal,(\Xcal\oplus_{\Lambda_{1/n}}\set{\bullet})^{\otimes n}\big)
      \leq
      n\cdot\ent(\Lambda_{1/n})
      $
    \item
      $\displaystyle
      \aikd\big(
      (\Xcal\otimes\Ycal)\oplus_{\Lambda_{1/n}}\set{\bullet},
      (\Xcal\oplus_{\Lambda_{1/n}}\set{\bullet})
      \otimes
      (\Ycal\oplus_{\Lambda_{1/n}}\set{\bullet})
      \big)
      \leq
      3\ent(\Lambda_{1/n})
      $
    \item
      $\displaystyle
      \aikd\big((\Xcal\oplus_{\Lambda_{1/n}}\set{\bullet}),
      (\Ycal\oplus_{\Lambda_{1/n}}\set{\bullet})\big)
      \leq
      \frac1n\aikd(\Xcal,\Ycal)
      $
    \end{enumerate}
  \end{lemma}
  
  The proof can be found on page \pageref{p:mixtures-dist1.rep}.  Note,
  that the distance estimates in the lemma above are with respect to
  the asymptotic Kolmogorov distance. This is essential, since from
  the perspective of the intrinsic Kolmogorov distance mixtures are
  very badly behaved.

	\section[Tropical probability]{Tropical probability spaces 
		and their configurations}
	\label{s:tropical}
	  In this section we introduce the notion of \term[tropical
  probability space]{tropical probability spaces} and
  their \term[configuration of tropical probability
  spaces]{configurations}. Configurations of tropical
  probability spaces are points in the asymptotic cone of the space
  $\prob\<\Gbf>$, that is they are ``limits'' of certain divergent
  sequences of ``normal'' configurations. We will first give
  the construction of an asymptotic cone in an abstract context. Next,
  we will apply the construction to the particular case of
  configurations of probability spaces. For some background on asymptotic cones, see for instance \cite{Burago-Course-2001}.
     
\subsection{Asymptotic cones of metric spaces.}
\label{s:tropical-ac}
  The \term[asymptotic cone]{asymptotic cone} captures large-scale
  geometry of a metric space.  Abstractly, the asymptotic cone of a
  pointed metric space is the pointed Gromov-Hausdorff limit of the
  sequence of spaces obtained from the given one by scaling down the
  metric. Of course, convergence is in general by no means
  assured. Sometimes a weaker type of convergence (using ultrafilters)
  is considered.  Since, in our case, the asymptotic cone can be
  evaluated relatively explicitly we do not give the definition of
  Gromov-Hausdorff convergence or convergence with respect to an
  ultrafilter here, but instead give a construction.
  
  We would like to understand asymptotic cones of the space of
  configurations of probability spaces, considered as a metric space
  with the pseudo-metric $\ikd$ or $\aikd$. For a fixed complete
  diagram category $\Gbf$ the space $\prob\<\Gbf>$ is a monoid with
  operation $\otimes$.  It has the additional property that shifts are
  non-expanding maps.  This simplifies the construction and analysis
  of its asymptotic cone.  In fact, as we will see later the metric
  $\aikd$ is already \emph{asymptotic} relative to $\ikd$. The
  application of the asymptotic cone construction to the metric
  $\aikd$ allows us to obtain a complete metric space with a simple
  description of points in it.
  
  Note that even though the monoid
  $(\prob\<\Gbf>,\otimes)$ is not Abelian it has the property that for
  any $\Xcal_{1},\Xcal_{2}\in\prob\<\Gbf>$ one
  has \[ \ikd(\Xcal_{1}\otimes\Xcal_{2},\Xcal_{2}\otimes\Xcal_{1})=0 \]
  Thus, from a metric perspective it is as good as being Abelian.

\subsubsection{Metrics versus pseudo-metrics}
\label{s:tropical-ac-metrics} 
  A pseudo-metric $\dist$ on a set $X$ is a bivariate function
  satisfying all the axioms of a distance function except it is
  non-negative definite rather than positive definite. That is, the
  pseudo-distance function is allowed to vanish on pairs of
  non-identical points. A set equipped with a pseudo-metric will be
  called a pseudo-metric space. An isometry of such spaces is a
  distance-preserving map, such that for any point in the target space
  there is a point in the image, which is distance zero from it.
  Given such an pseudo-metric space $(X,\dist)$ one could always
  construct an isometric metric space $(X/_{\dist=0}\,,\dist)$ by
  identifying all pairs of points that are distance zero apart.

  Any property formulated in terms of the pseudo-metric holds
  simultaneously for a pseudo-metric space and its metric quotient.
  It will be convenient for us to construct pseudo-metrics on spaces
  instead of passing to the quotient spaces.
 
\subsubsection{Asymptotic cone of a metric Abelian monoid}
 Let $(\Gamma,\otimes,\dist)$ be a monoid with a pseudo-metric
 $\dist$, which satisfies the following properties
 \begin{enumerate}
 \item 
   The shifts
   \[
   \cdot \otimes \gamma' : 
   \Gamma \to \Gamma, \quad 
   \gamma \mapsto \gamma \otimes \gamma'
   \]
   are non-expanding for any $\gamma' \in \Gamma$
 \item
   For any $\gamma,\gamma'\in\Gamma$ holds
  \[
  \dist(\gamma \otimes \gamma' , \gamma' \otimes \gamma) = 0
  \]   
 \end{enumerate}
  We will call a monoid with pseudo-metric that satisfies these
  conditions a \term{metric Abelian monoid}.  It follows from the
  shift-invariance property that for any
  $\gamma_{1},\gamma_{2},\gamma_{3}\in\Gamma$ holds
  \[\tageq{shifts}
    \dist(\gamma_{1}\otimes\gamma_{3},\gamma_{2}\otimes\gamma_{3})
    \leq
    \dist(\gamma_{1},\gamma_{2})
  \]
  and for any quadruple
  $\gamma_{1},\gamma_{2},\gamma_{3}, \gamma_{4} \in\Gamma$ holds
  \[\tageq{subadditive}
    \dist(\gamma_{1}\otimes\gamma_{2},\gamma_{3}\otimes\gamma_{4})
    \leq
    \dist(\gamma_{1},\gamma_{3}) + \dist(\gamma_{2},\gamma_{4})
  \]
  and, in particular, the monoid operation is 1-Lipschitz with
  respect to each argument.

  As a direct consequence, for every $n \in \Nbb$, and $\gamma_1,
  \gamma_2 \in \Gamma$ also holds
  \[\tageq{contractiondist}
    \dist(\gamma_1^{\otimes n}, \gamma_2^{\otimes n} ) 
    \leq 
    n\dist(\gamma_1, \gamma_2)
  \]

  For a sequence $\bar\gamma=\set{\gamma(i)}\in\Gamma^{\Nbb_{0}}$
  define its \term{defect} with respect to the distance function
  $\dist$ by
  \[
    \defect_{\dist}(\bar\gamma)
    =
    \sup_{i,j\in\Nbb_{0}}\dist\big(\gamma(i+j),\gamma(i)\otimes\gamma(j)\big)
  \] 
 
  The sequence $\bar\gamma$ will be called \term[linear
  sequence]{$\dist$-linear} if $\defect_{\dist}(\bar\gamma)=0$,
  and \term[quasi-linear sequence]{$\dist$-quasi-linear} if
  $\defect_{\dist}(\bar\gamma)<\infty$.  Denote by
  $\lin_{\dist}(\Gamma)$ and $\qlin_{\dist}(\Gamma)$ the sets of all
  linear and, respectively, quasi-linear sequences in $\Gamma$ with
  respect to the distance $\dist$.

  For two elements $\bar\gamma_1,\bar\gamma_{2}\in\qlin_{d}(\Gamma)$, define an
  asymptotic distance between them by
  \[
    \dista(\bar\gamma_{1},\bar\gamma_{2})
    :=
    \lim_{n\to\infty}\frac1n\dist\big(\gamma_{1}(n),\gamma_{2}(n)\big)
  \]

  \begin{lemma}{p:adistonql}
    For a pair $\bar\gamma_{1},\bar\gamma_{2}\in\qlin_{\dist}(\Gamma)$
    the limit
    \[
      \lim_{n\to\infty}\frac1n\dist\big(\gamma_{1}(n),\gamma_{2}(n)\big)
    \]
    exists and is finite.
  \end{lemma}
  We provide the proof in Section~\ref{s:technical} on
  page~\pageref{p:adistonql.rep}. 

  The bivariate function $\dista$ is a pseudo-distance on the set
  $\qlin_{\dist}(\Gamma)$. We call two sequences
  $\bar\gamma_{1},\bar\gamma_{2}\in\qlin_{\dist}(\Gamma)$
  \term[asymptotic equivalence]{asymptotically equivalent} if
  $\dista(\bar\gamma_{1},\bar\gamma_{2})=0$ and write 
  \[
  \bar \gamma_{1} \aseq{\dista}  \bar \gamma_2
  \]

  We will call a sequence $\bar\gamma$ \term[weakly quasi-linear
  sequence]{weakly quasi-linear}, if it is asymptotically equivalent
  to a quasi-linear sequence. Note that the space of all weakly
  quasi-linear sequences can also be endowed with the asymptotic
  distance and it is isometric to the space of quasi-linear
  sequences. As we will see later all the natural operations we consider are
  $\dista$-Lipschitz and therefore coincide for the asymptotically
  equivalent sequences. Thus given a weakly quasi-linear sequence we
  could always replace it by an equivalent quasi-linear sequence
  without any visible effect. Thus, we take the liberty to omit the
  adverb ``weakly''. Whenever we say quasi-linear sequence, we mean
  a weakly quasi-linear sequence, that is silently replaced
  by an asymptotically equivalent genuine quasi-linear sequence, if
  necessary.

  The validity of the following constructions is very easy to verify,
  so we omit the proofs.

  The set $\qlin_{\dist}(\Gamma)$ admits an action of the
  multiplicative semigroup $(\Rbb_{\geq0},\,\cdot\,)$
  defined in the following way. Let $\lambda\in\Rbb_{\geq 0}$ and
  $\bar\gamma=\set{\gamma(n)}\in\qlin_{\dist}(\Gamma)$. Then define
  the action of $\lambda$ on $\bar\gamma$
  by 
  \[\tageq{R-action} 
  {\bar\gamma}^\lambda
  := 
  \set{\gamma(\lfloor\lambda\cdot n\rfloor)}_{n\in\Nbb_{0}} 
  \]
  This is only an action up to asymptotic equivalence. 
  Similarly, in the constructions that follow we are tacitly assuming they are valid up to asymptotic equivalence.

  The action 
  \[
  \cdot:\Rbb_{\geq 0}\times(\qlin_{\dist}(\Gamma),\dista)
  \to
  (\qlin_{\dist}(\Gamma),\dista)
  \]
  is continuous with respect to $\dista$ and, moreover it is
  a homothety (dilation), that is
  \[
  \dista(\bar\gamma_{1}^\lambda,\bar\gamma_{2}^\lambda)
  =
  \lambda\cdot\dista(\bar\gamma_{1},\bar\gamma_{2})
  \]

  The group operation $\otimes$ on $\Gamma$ induces a
  $\dista$-continuous (in fact, 1-Lipschitz) group
  operation on $\qlin_{\dist}(\Gamma)$ by multiplying sequences
  element-wise.  The semigroup structure on $\qlin_{\dist}(\Gamma)$
  is distributive with respect to the $\Rbb_{\geq 0}$-action
  \begin{align*}
    (\bar\gamma_{1}\otimes\bar\gamma_{2})^\lambda
    &=
    \bar\gamma_{1}^{\lambda}\otimes\bar\gamma_{2}^{\lambda}\\
   \bar\gamma^{ \lambda_{1}+\lambda_{2}}
    &\aseq{\dista}
    \bar\gamma^{\lambda_1}\otimes \bar\gamma^{\lambda_2}
  \end{align*}

  In particular for $n \in \mathbb{N}$
  \[
  \bar{\gamma}^n \aseq{\dista} \bar\gamma^{\otimes n}
  \]

  The path
  \[
  [0,1] \ni \lambda \mapsto \bar\gamma_1^{1-\lambda} \otimes \bar\gamma_2^{\lambda}
  \]
  will be called a convex interpolation and is a constant-speed
  $\dista$-geodesic between $\bar\gamma_1$ and $\bar\gamma_2$, that is
  for $\lambda \in [0,1]$,
  \begin{align*}
  \dista(\bar\gamma_1^{(1-\lambda)} \otimes \bar\gamma_2^{\lambda} , 
         \bar\gamma_1) 
  &= 
  \lambda \dista(\bar\gamma_1,\bar\gamma_2)
  \\
    \dista(\bar\gamma_1^{(1-\lambda)} \otimes \bar\gamma_2^{\lambda}, 
           \bar\gamma_2) 
    &= 
    (1-\lambda) \dista(\bar\gamma_1,\bar\gamma_2)
  \end{align*}

\subsubsection{Conditions for completeness}
  We would like to call
  \[
  \Gamma^{(\infty)}_{\dist}
  :=
  (\qlin_{\dist}(\Gamma),\otimes,\cdot,\dista)
  \]    
  the asymptotic cone of $(\Gamma,\otimes,\dist)$. However it is not
  clear in general, whether $\Gamma^{(\infty)}_{\dist}$ is a complete
  space.
 
  We can simply consider the metric completion, and call it the asymptotic
  cone of $(\Gamma,\otimes,\dist)$. We feel, however, that it adds
  just another level of obscurity as to what the points of
  $\Gamma^{(\infty)}_{\dist}$ are.

  Under some circumstances, however, the completeness of the space of
  quasi-linear sequences comes for free. This is the subject of the
  proposition below.

  Suppose the metric Abelian monoid $(\Gamma,\otimes,\dist)$ has an
  additional property: There exists a constant $C>0$, such that for
  any quasi-linear sequence $\bar\gamma\in\qlin_{\dist}(\Gamma)$,
  there exists an asymptotically equivalent quasi-linear sequence
  $\bar\gamma'$ with defect bounded by $C$.  If this is the case, we
  say that the metric monoid $(\Gamma,\otimes,\dist)$ has the
  ($C$-)\term{uniformly bounded defect property}.

  \begin{proposition}{p:boundeddefect}
    Suppose $(\Gamma,\otimes,\dist)$ is a metric Abelian monoid such that 
    \begin{enumerate}
    \item
      \label{pi:homogeneity} 
      the distance function $\dist$ is homogeneous, that is
      for any $\gamma_{1},\gamma_{2}\in\Gamma$ and $n\in\Nbb_{0}$ 
      \[
      \dist(\gamma_{1}^{\otimes n},\gamma_{2}^{\otimes n})
      =
      n\cdot\dist(\gamma_{1},\gamma_{2})
      \]
    \item
      \label{pi:boundeddefect} 
      $(\Gamma,\otimes,\dist)$ has the uniformly bounded defect property.
    \end{enumerate}
    Then the space $(\qlin_{\dist}(\Gamma),\dista)$ is complete.
  \end{proposition}
  The proof of the proposition can be found on
  page~\pageref{p:boundeddefect.rep}.

\subsubsection{On the density of linear sequences}

  In Section~\ref{s:ac-aep-conf} we have shown that Bernoulli sequences of
    configurations can be approximated by sequences of homogeneous
    configurations. The proposition below will allow us to extend this
    statement to a wider class of sequences.  It gives a sufficient condition under which
  the linear sequences are dense in the quasi-linear sequences.

  \begin{proposition}{p:eps-linear-dense}
    Suppose $(\Gamma, \otimes, \dist)$ has the $\epsilon$-uniformly
    bounded defect property for every $\epsilon > 0$. Then
    $\lin_{\dist}(\Gamma)$ is dense in $\qlin_{\dist}(\Gamma)$
  \end{proposition}
  See page~\pageref{p:eps-linear-dense.rep} for the proof.

\subsubsection{Asymptotic metric on original semigroup}
\label{suse:metric-original-group}
  Starting with an element $\gamma\in\Gamma$ one can construct a
  linear sequence $\linseq\gamma=\set{\gamma^{\otimes
      i}}_{i\in\Nbb_{0}}$. In view of inequality
  (\ref{eq:contractiondist}), this map is a contraction
  \[\tageq{inclusions}
  \big(\Gamma, \dist \big) \to \big(\lin_{\dist}(\Gamma), \dista\big)
  \]
 
  By the inclusions in (\ref{eq:inclusions}) we have an induced metric
  $\dista$ on $\Gamma$, satisfying for any
  $\gamma_{1},\gamma_{2}\in\Gamma$
  \[\tageq{deltalessd}
  \dista(\gamma_{1},\gamma_{2})\leq\dist(\gamma_{1},\gamma_{2})
  \]
  and the following scale-invariance condition is gained
  \[\tageq{deltalin}
  \dista(\gamma_{1}^{\otimes n},\gamma_{2}^{\otimes n})
  =
  n\cdot\dista(\gamma_{1},\gamma_{2})
  \]
   for all $n \in \Nbb_0$.

  Note moreover that if $\dist$ was scale-invariant to begin with,
  then $\dista$ coincides with $\dist$ on $\Gamma$.

\subsubsection{Iteration of construction}
\label{suse:iterate-construction}
  We may now iterate the constructions above, that is, we may apply
  them to $(\Gamma, \dista)$ instead of $(\Gamma, \dist)$.  One may
  wonder what is the purpose.  However, we have already observed that
  $\dista$ satisfies the scale-invariance condition
  (\ref{eq:deltalin}), which is one of the conditions going into a
  proof of completeness in Proposition \ref{p:boundeddefect}.
  Moreover, when we will later apply the theory in this section to the
  particular case of $\Gamma = \prob\<\Gbf>$, we will see that
  \[
  (\Gamma, \dista) = ( \prob\<\Gbf> , \aikd )
  \]
  and we will show that the latter space has the $\epsilon$-uniformly
  bounded defect property for every $\epsilon > 0$.
 
  By virtue of the bound $\dista \leq \dist$, sequences that are
  quasi-linear with respect to $\dista$, are also quasi-linear with
  respect to $\dist$.  Since $\dista$ is scale-invariant, the
  associated asymptotic distance $\distaa$ coincides with $\dista$ on
  $\Gamma$. We will show (in Lemma \ref{p:dist-dista-isometry} below)
  that $\distaa$ also corresponds to $\dista$ on $\dist$-quasi-linear
  sequences.
 
  In order to organize all these statements, and to be more precise,
  let us include the spaces in the following commutative diagram.

  \[\tageq{tropical-diagram-0}
  \begin{tikzcd}[ampersand replacement=\&,row sep=tiny]
    \mbox{}
    \&
    \big(\lin_{\dist}(\Gamma), \dista)
    \arrow[hookrightarrow]{dd}{\i_{1}}
    \arrow[hookrightarrow]{r}{\j_{1}}
    \&
    \big(\qlin_{\dist}(\Gamma), \dista )
    \arrow[hookrightarrow]{dd}{\i_{2}}
    \\
    (\Gamma,\dista)
    \arrow{ru}{f}
    \arrow{rd}{\phi}
    \\
    \mbox{}
    \&
    \big(\lin_{\dista}(\Gamma),\distaa)
    \arrow[hookrightarrow]{r}{\j_{2}}
    \&
    \big(\qlin_{\dista}(\Gamma),\distaa)
  \end{tikzcd}
  \]

  The maps $f, \phi$ and $\i_1$ are isometries.  The maps $\j_1$ and
  $\j_2$ are isometric embeddings.  The next lemmas show that $\i_2$
  is also an isometric embedding, and it has dense image.

  \begin{lemma}{p:dist-dista-isometry}
    The natural inclusion 
    \[
    \i_{2}:(\qlin_{\dist}(\Gamma),\dista)
    \into
    (\qlin_{\dista}(\Gamma),\distaa)
    \]
    is an isometric embedding.
  \end{lemma}
  
  \begin{lemma}{p:dist-dista-dense}
    The image of the isometric embedding
    \[
    \i_{2}:(\qlin_{\dist}(\Gamma),\dista)
    \into
    (\qlin_{\dista}(\Gamma),\distaa)
    \]
    is dense in $(\qlin_{\dista}(\Gamma),\distaa)$
  \end{lemma}
  The proofs of the two lemmas above are to be found on
  page~\pageref{p:dist-dista-isometry.rep}.

\subsection{Tropical probability spaces and configurations}
  Now we apply the above construction to the space of complete
  configurations with fixed combinatorial type $\Gbf$.

  Fix a complete diagram category $\Gbf$ and consider the space $\prob\<\Gbf>$
  of configurations modeled on $\Gbf$. 
  It carries the following structures:
  \begin{enumerate}
  \item A pseudo-metric $\ikd$ or $\aikd$.
  \item A 1-Lipschitz tensor product $\otimes$.
  \item A 1-Lipschitz entropy function
    $\ent_{*}:\prob\<\Gbf>\to\Rbb^{\size{\Gbf}}$.
  \end{enumerate}

  The tensor product of configurations is commutative from a metric
  perspective.  Recall that in Corollary~\ref{p:subadditivity} the
  subadditivity of both $\ikd$ and $\aikd$ was established, namely for
  any $\Xcal, \Ycal, \Ucal, \Vcal \in \prob\<\Gbf>$ holds
  \[
  \ikd(\Xcal \otimes \Ucal, \Ycal \otimes \Vcal)
  \leq 
  \ikd(\Xcal, \Ycal) + \ikd( \Ucal, \Vcal ).
  \]	
  and 
  \[
  \aikd(\Xcal \otimes \Ucal, \Ycal \otimes \Vcal)
  \leq 
  \aikd(\Xcal, \Ycal) + \aikd( \Ucal, \Vcal ).
  \]

  The space $(\prob\<\Gbf>, \otimes, \ikd)$ is a metric Abelian
  monoid.  Note also that $\hat{\ikd} = \aikd$ on $\prob\<\Gbf>$,
  along the lines of Section \ref{suse:metric-original-group}.
  
  However, the metric $\ikd$ is not scale-invariant. Moreover, it is
  unclear whether the metric semigroup $(\prob\<\Gbf>, \otimes ,\ikd)$
  has the uniformly bounded defect property. This is why we iterate
  the construction, as announced in Section
  \ref{suse:iterate-construction}, and consider the space of
  $\aikd$-quasi-linear sequences instead.
  
  \begin{lemma}{p:unifsmalldefectaikd}
    For a complete diagram category, and for every $\epsilon > 0$, the
    space $(\prob\<\Gbf>,\otimes,\aikd)$ has the $\epsilon$-uniformly
    bounded defect property, that is for any $\aikd$-quasi-linear
    sequence $\bar\Xcal\in\qlin_{\aikd}(\prob\<\Gbf>)$ there exists an
    asymptotically equivalent sequence $\bar\Ycal$ with defect not
    exceeding $\epsilon$.
  \end{lemma}

  By applying the general setup in the previous section to the metric
  semigroups $(\prob\<\Gbf>, \otimes, \ikd)$ and $(\prob\<\Gbf>,
  \otimes, \aikd)$ and as a corollary to Lemma
  \ref{p:unifsmalldefectaikd} we obtain the following theorem.

  \begin{theorem}{p:corollary-acone-ikd-aikd}
    Consider the commutative diagram
    \[\tageq{tropical-diagram}
    \begin{tikzcd}[ampersand replacement=\&,row sep=tiny]
      \mbox{}
      \&
      \big(\lin_{\ikd}(\prob\<\Gbf>) , \hat\ikd \big)
      \arrow[hookrightarrow]{dd}{\i_{1}}
      \arrow[hookrightarrow]{r}{\j_{1}}
      \&
      \big(\qlin_{\ikd}(\prob\<\Gbf>) , \hat\ikd\big)
      \arrow[hookrightarrow]{dd}{\i_{2}}
      \\
      (\prob\<\Gbf>,\aikd)
      \arrow{ru}{f}
      \arrow{rd}{\phi}
      \\
      \mbox{}
      \&
      \big(\lin_{\aikd}(\prob\<\Gbf>), \hat\aikd \big)
      \arrow[hookrightarrow]{r}{\j_{2}}
      \&
      \big(\qlin_{\aikd}(\prob\<\Gbf>), \hat\aikd \big)
    \end{tikzcd}
    \]
    Then the following statements hold:
    \begin{enumerate}
    \item 
      The maps $f, \phi, \i_1$ are isometries.
    \item 
      The maps $\i_2, \j_1,\j_2$ are isometric embeddings and each map
      has a dense image in the corresponding target space.
    \item 
      The space in the lower-right corner,
      $\big(\qlin_{\aikd}(\prob\<\Gbf>), \hat\aikd \big)$, is
      complete.
    \end{enumerate}
  \end{theorem} 
  
  We may finally define the space of \term[tropical
    configuration]{tropical $\Gbf$-configurations}, as the space in
  the lower-right corner of the diagram
  \[
  \prob\<\Gbf>^{(\infty)}
  :=
  \big(\qlin_{\aikd}(\prob\<\Gbf>),\otimes,\cdot,\hat\aikd\big)
  \]
  By the Theorem~\ref{p:corollary-acone-ikd-aikd} above, this space is
  complete.
  
  The entropy function $\ent_{*}:\prob\<\Gbf>\to\Rbb^{\size{\Gbf}}$
  extends to a linear functional
  \[
  \ent_{*}:\prob\<\Gbf>^{(\infty)}\to(\Rbb^{\size{\Gbf}},|\,\cdot\,|_{1})
  \]
  of norm one, defined by
  \[
  \ent_{*}(\bar\Xcal) = \lim_{n\to \infty } \frac1n \ent_{*} \big(\Xcal(n)\big)
  \]
  
  \subsubsection{Sequences of homogeneous configurations are dense.}
  Let $\tilde\lin_{\ikd}(\prob\<\Gbf>_{\hbf})$ stand for the weakly linear
  sequences of homogeneous configurations, that is those sequences,
  that are asymptotically equivalent to a linear sequence (not
  necessarily of homogeneous spaces).
  
  For a sequence of homogeneous spaces
  $\bar\Hcal \in \tilde\lin_{\ikd}(\prob\<\Gbf>_{\hbf})$ define
  $\aep( \bar\Hcal )$ to be a $\ikd$-linear sequence asymptotically
  equivalent to $\bar\Hcal$.
  
  Now we can extend the commutative
  diagram~(\ref{eq:tropical-diagram}) as follows
  \[\tageq{tropical-diagram-hom}
  \begin{tikzcd}[ampersand replacement=\&,row sep=tiny]
  \tilde\lin_{\ikd}(\prob\<\Gbf>_{\hbf})
  \arrow{r}{\aep}
  \&
  \lin_{\ikd}(\prob\<\Gbf>)
  \arrow[hookrightarrow]{dd}{\i_{1}}
  \arrow[hookrightarrow]{r}{\j_{1}}
  \&
  \qlin_{\ikd}(\prob\<\Gbf>)
  \arrow[hookrightarrow]{dd}{\i_{2}}
  \\
  (\prob\<\Gbf>,\aikd)
  \arrow{ru}{f}
  \arrow{rd}{\phi}
  \\
  \mbox{}
  \&
  \lin_{\aikd}(\prob\<\Gbf>)
  \arrow[hookrightarrow]{r}{\j_{2}}
  \&
  \prob\<\Gbf>^{(\infty)}    
  \end{tikzcd}
  \]

  By the Asymptotic Equipartition Property for configurations,
  Theorem~\ref{p:aep-complete}, the map $\aep$ is an isometry, hence
  we have the following theorem.
    \begin{theorem}{p:homo-dense}
    The map
    \[
    \j_2\circ\i_1\circ\aep:
    \tilde\lin_{\ikd}(\prob\<\Gbf>_{\hbf})\to\prob\<\Gbf>^{(\infty)}
    \]
    is an isometric embedding with dense image.
  \end{theorem}

  Let $\prob\<\Gbf>_{\hbf}^{(\infty)} \subset \prob\<\Gbf>^{(\infty)}$
  denote the space of weakly quasi-linear sequences of configurations
  $\bar{\Hcal} \in \prob\<\Gbf>^{(\infty)}$, such that for every
  $n \in \Nbb_0$, the configuration $\Hcal(n)$ is homogeneous.  We
  will refer to $\prob\<\Gbf>_{\hbf}^{(\infty)}$ as the space of
  homogeneous tropical configurations.

  Denote by $\aep$ the embedding 
  \[
    \aep: 
    \prob\<\Gbf>_{\hbf}^{(\infty)} 
    \hookrightarrow 
    \prob\<\Gbf>^{(\infty)}
  \]  
  
  \begin{theorem}{p:aep-tropical-config}{\rm(Asymptotic Equipartition Theorem for tropical configurations)}
  Let $\Gbf$ be a complete diagram category. Then the map 
  \[
  \aep: \prob\<\Gbf>_{\hbf}^{(\infty)} \hookrightarrow \prob\<\Gbf>^{(\infty)}
  \]
  is an isometry.
  \end{theorem}
  
  \begin{Proof}
  	 We need to show that for every tropical configuration $\bar\Xcal \in \prob\<\Gbf>^{(\infty)}$, there 
  	exists a homogeneous tropical configuration $\bar \Hcal \in \prob\<\Gbf>_\hbf^{(\infty)}$ such that
  	\[
  	\hat{\aikd}(\bar\Hcal, \bar\Xcal) = 0
  	\]
  By Lemma \ref{p:unifsmalldefectaikd} and Proposition \ref{p:eps-linear-dense}, for every $j \in \Nbb$ there exists a sequence $\bar\Ycal_j \in \lin_{\aikd}(\prob\<\Gbf>) \cong \lin_{\ikd}(\prob\<\Gbf>)$	such that 
  \[
  \hat{\aikd}(\bar{\Ycal}_j, \bar\Xcal) \leq \frac{1}{j}
  \]
  By the Asymptotic Equipartition Property for configurations, Theorem \ref{p:aep-complete},
  there are sequences of homogeneous configurations $\bar\Hcal_j$ such that 
  \[
  \hat{\aikd}(\bar\Ycal_j, \bar\Hcal_j) = 0
  \]
  Define $\mathbf{i}(j)$ such that for all $k \geq \mathbf{i}(j)$
  \[
  \frac{1}{k} \aikd(\Ycal_j(k), \Hcal_j(k)) \leq \frac{1}{j}
  \]
  and  moreover
  \[
  \frac{1}{k} \aikd(\Ycal_j(k), \Xcal(k)) \leq \frac{2}{j}
  \]
  The function $\mathbf{i}$ can be chosen monotonically increasing. 
  For every $i \in \Nbb_0$ there is a unique $\jbf(i) \in \Nbb_0$ such that 
  \[
  \mathbf{i}\big(\jbf(i)\big) \leq i < \mathbf{i}\big(\jbf(i)+1\big)
  \]
  Define then
  \[
  \Hcal(k) = \Hcal_{\jbf(k)}(k)
  \]
  It follows that for $k > \mathbf{i}(1)$,
  \[
  \frac{1}{k}\aikd(\Hcal(k), \Xcal(k)) \leq \frac{3}{\jbf(k)}
  \]
  Since $\jbf$ is a non-decreasing, divergent sequence, the theorem follows.
  \end{Proof}  

  \bigskip
  
  Thus we have shown that all arrows in
  diagram~(\ref{eq:tropical-diagram-hom}) are isometric embeddings
  with dense images.  We would like to conjecture, that, in fact, they
  are all isometries.  In any case, the difference between metrics
  $\hat\aikd$ and $\aikd$ is so small ($\aikd$ is defined on the dense
  subset of the domain of definition of $\hat\aikd$ and they coincide
  whenever both are defined), that we will not write the hat anymore
  and just use notation $\aikd$ for the metric on
  $\prob\<\Gbf>^{(\infty)}$.

\subsection{Tropical probability spaces and tropical chains}
\label{se:tropprobtropchains}
  In this section we evaluate the spaces $\prob^{(\infty)}$ and
  $\prob\<\Cbf_{n}>^{(\infty)}$, where $\Cbf_n$ is a chain, which is
  the diagram category introduced in \ref{s:config-examples-chain} on
  page \pageref{s:config-examples-chain}.

  Recall that a finite probability space $U$ is homogeneous if
  $\Aut(U)$ acts transitively on the support of the measure. The property
  of being homogeneous is invariant under isomorphism and every
  homogeneous space is isomorphic to a probability space with the
  uniform distribution. 

  Homogeneous chains also have a very simple description. A chain of
  reductions is homogeneous, if and only if all the individual spaces
  are homogeneous.

  This simple description allows us to evaluate explicitly the
  Kolmogorov distance on the spaces of weakly linear sequences of
  homogeneous chains and consequently the space of tropical chains.

  \begin{theorem}{p:tropical-simple}\rule{0mm}{1mm}\\
    \begin{enumerate}
    \item 
      $\displaystyle
      \prob^{(\infty)}
      \cong
      (\Rbb_{\geq0},|\cdot - \cdot|,+,\cdot)
      $
      \vspace{3mm}
    \item 
      $\displaystyle
      \prob\<\Cbf_{n}>^{(\infty)}
      \cong
      \set{
        \left(
        \begin{matrix}
          x_{1}\\\vdots\\x_{n}
        \end{matrix}
        \right) \in\Rbb^{n} \st 0\leq x_{n}\leq \cdots\leq x_{1}} 
      $
      \vspace{3mm}
      \\ 
      where the right-hand side is a cone in $(\Rbb^{n},|\,\cdot\,|_{1})$.
    \end{enumerate}
  \end{theorem}

  To prove the Theorem~\ref{p:tropical-simple} we evaluate first
  the isometry class of the space of weakly linear sequences of
  homogeneous spaces (or chains). We will only present an argument for
  single spaces, since the argument for chains is very similar.

  \begin{lemma}{p:tropical-homo-simple}
    \[
    \tilde\lin_{\ikd}(\prob_{\hbf})
    \cong
    (\Rbb_{\geq0},|\cdot - \cdot|,+,\cdot)
    \]
  \end{lemma}

  Note that the right-hand side is a complete metric space, thus the
  Asymptotic Equipartition Theorem for tropical configurations,
  Theorem ~\ref{p:homo-dense}, together with
  Lemma~\ref{p:tropical-homo-simple}, imply
  Theorem~\ref{p:tropical-simple}.
 
  To prove Lemma~\ref{p:tropical-homo-simple} we need to evaluate the
  Kolmogorov distance between two homogeneous spaces, or chains of
  homogeneous spaces. This is the
  subject of the next lemma, from which
  Lemma~\ref{p:tropical-homo-simple} follows immediately.

\begin{lemma}{p:ikduniform}
  Denote by $U_{n}$ a finite uniform probability space of cardinality
  $n$, then
  \begin{enumerate}
  \item 
    \[
    \ikd(U_{n},U_{m})\leq 2 \ln 2 + \left|\ln\frac nm\right|
    \]
  \item
    \[
    \aikd(U_{n},U_{m})=|\ent(U_{n})-\ent(U_{m})|
    \]
  \end{enumerate}
\end{lemma}

\subsection{Stochastic processes}
  Often, stochastic processes naturally give rise to
  $\aikd$-quasi-linear sequences.  We include this last subsection as
  an indication that our statements, together with the construction of
  the tropical cone, have a much larger reach than sequences of
  independent random variables.  We will be brief, and come back to
  the topic in a subsequent article.

  For a minimal diamond configuration
  
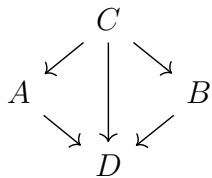
\[
  \begin{tikzcd}[row sep=small, column sep=small]
    \mbox{}
    & 
    C
    \arrow{dl}{}
    \arrow{dr}{}
    \arrow{dd}{}
    &
    \\
    A
    \arrow{dr}{}
    & 
    &
    B
    \arrow{dl}{}
    \\
    \mbox{}
    &
    D
    &
  \end{tikzcd}
  \]
  we define the conditional mutual information between $A$ and $B$ given $D$ by
  \[
  \MI(A; B \, \rel \, D) := \ent(A) + \ent(B) - \ent(C) - \ent(D)
  \]
  Shannon's inequality (\ref{eq:shannonineq}) says that the
  conditional mutual information is always non-negative.  Any minimal
  two-fan $A \ot C \to B$ can be completed to a diamond with the
  one-point probability space $\set{\bullet}$ as the terminal vertex,
  and the mutual information between $A$ and $B$ is defined as
  \[
  \MI(A; B) = \MI(A; B \, \rel \, \set{ \bullet }) = \ent(A) + \ent(B) - \ent(C)
  \]
  
  Let 
  \[
  \ldots, X_{-1}, X_0, X_1, \ldots
  \]
  be a stationary stochastic process with finite state space $\un
  X$. Thus, for any $I=\set{k,k+1,\ldots,l}\subset\Zbb$ we have
  jointly distributed random variables $X_{k},\ldots,X_{l}$, that
  generate a full configuration
  \[
  \Xcal_{I}
  =
  \langle X_{k},\ldots,X_{l}\rangle
  =
  \set{X_{J}}_{J\subset I}
  \]
  as explained in Section~\ref{s:config-examples-full}.
  This collection of full configurations is consistent in the sense that
  for $k\leq k'\leq l'\leq l$ there are canonical isomorphisms  
  \[
  \Xcal_{I'}
  \cong 
  R^{*}_{I',I}\Xcal_{I}
  \]
  where $I:=\set{k,\ldots,l}$, $I':=\set{k',\ldots,l'}$ and
  $R^*_{I',I}$ is the restriction operator introduced in
  Section~\ref{s:config-restrictions-fullfull}.
  
  The property of being
  stationary means that there are canonical isomorphisms for any finite
  subset $I\subset\Zbb$ and $l\in\Zbb$
  \[
  \Xcal_{I}\too[\cong]\Xcal_{I+l}
  \]
  
  For $I=\set{k,k+1,\ldots,l}$ we call the initial space $X_{I}$ of the
  configuration $\Xcal_{I}$ the
  space of trajectories of the process over $I$ and denote it $X_{k}^{l}$.
  
  Note that by stationarity, for every $m \in \Zbb$, $k\in
  \Nbb$ and $l \in \Nbb_0$,
  \[
  \MI(X_m^{m+k-1}; X_{m+k}^{m+k + l-1}) = \MI(X_{-k+1}^{0}; X_{1}^l )
  \]
  Moreover the right-hand side is an increasing function of both $k$
  and $l$.  We make the following important observation.  The defect
  of the sequence $n \mapsto X_1^n$ is equal to
  \begin{align*}
    \defect_{\aikd}\left( \{ X_1^n\} \right) 
    &= 
    \sup_{m, n \in \Nbb_0} 
    \aikd\left( X_1^{m+n}, X_1^m \otimes X_1^n\right) 
    \\
    &= 
    \sup_{m, n \in \Nbb_0} 
    \left| \ent(X_1^m \otimes X_1^n) - \ent(X_1^{m+n}) \right| 
    \\
    &= 
    \sup_{m, n \in \Nbb_0} 
    \left| \ent(X_{-m+1}^0 \otimes X_{1}^{n}) - \ent(X_1^{m+n})
    \right| 
    \\  
    &= 
    \sup_{m, n \in \Nbb_0} \MI(X_{-m+1}^0, X_1^n)
  \end{align*}
  Therefore, the sequence $n \mapsto X_1^n$ is $\aikd$-quasi-linear if
  and only if
  \[\tageq{stochastic-ql}
    \lim_{k,l\to \infty}\MI(X_{-k+1}^{0}, X_1^l) < \infty
  \]
  
  Once condition~(\ref{eq:stochastic-ql}) is satisfied for a stochastic
  process, it defines a tropical probability space $\bar X\in\prob$. 
  
  Note that condition~(\ref{eq:stochastic-ql}) is satisfied for
  any stationary, finite-state Markov chains.

	\let\standardthesection=\thesection
	\renewcommand{\thesection}{T}
	\section{Technical proofs}\label{s:technical}
	\let\thesubsectionstandard=\thesubsection
This section contains some proofs that did not make it into the
  main text.

\def\thesubsection{\thesection.\ref{s:config}}
\subsection{Statements from the section ``Configurations''}
\repeatclaim{p:minimalfansconfig}

\begin{Proof}
  We will need the following lemma
  \begin{tlemma}{p:minimising-reduction}
    Suppose we are given two two-fans of probability spaces
    \begin{align*}
      \Fcal
      &=
      (X\oot[\alpha]Z\too[\beta]Y)\\
      \Fcal''
      &=
      (X''\oot[\alpha''] Z''\too[\beta''] Y'')
    \end{align*}
    such that $\Fcal''$ is minimal. Let
    \[
    \Fcal\,\too[\mu]\,\Fcal'\!\!=\!(X\!\oot[\alpha']\! Z'\!\too[\beta']\! Y)
    \]
    be a minimal reduction of
    $\Fcal$. Then for any reduction $\rho:\Fcal\to\Fcal''$, there exists
    a reduction $\rho':\Fcal'\to\Fcal''$ such that $\rho=\rho'\circ \mu$
  \end{tlemma}
  \begin{Proof}
  	We define $\rho'$ on the terminal spaces of $\Fcal'$ to coincide with $\rho$.
   
    To prove the lemma we just need to provide a dashed arrow that makes the
    following diagram commutative
    
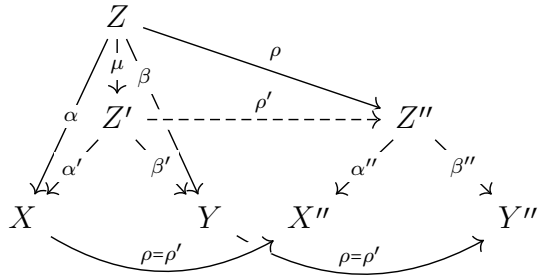
\[
    \begin{tikzcd}[row sep=normal, column sep=small, ampersand replacement=\&]
      \mbox{}
      \&
      Z
      \arrow{d}[description]{\mu}
      \arrow{ddl}[description]{\alpha}
      \arrow{ddr}[description, near start]{\beta}
      \arrow{drrr}{\rho}
      \\
      \&
      Z'
      \arrow{dl}[description]{\alpha'}
      \arrow{dr}[description]{\beta'}
      \arrow[dashrightarrow, crossing over]{rrr}{\rho'}
      \&
      \&
      \&
      Z''
      \arrow{dl}[description]{\alpha''}
      \arrow{dr}[description]{\beta''}
      \\
      X
      \&
      \&
      Y
      \arrow[bend right, crossing over]{rrr}{\rho=\rho'}
      \&
      X''
      \arrow[bend left, leftarrow, crossing over]{lll}[swap]{\rho=\rho'}
      \&
      \&
      Y''
    \end{tikzcd}
    \]
    The reduction $\rho'$ is constructed by simple diagram chasing and by
    using the minimality of $\Fcal''$.  
    Suppose $z'\in \un {Z'}$ and $z_{1},z_{2}\in\un Z$ are such that
    $z' = \mu(z_{1})=\mu(z_{2})$. By commutativity of the solid arrows in the
    diagram above, we have 
    \[ 
    \alpha''\circ\rho(z_{1})
    =
    \rho\circ\alpha'\circ\mu(z_{1})
    =
    \rho\circ\alpha'\circ\mu(z_{2})
    =
    \alpha''\circ\rho(z_{2})
    \]
    Similarly
    \[
    \beta''\circ\rho(z_{1})
    =
    \beta''\circ\rho(z_{2})
    \]
    Thus by minimality of $\Fcal''$ it follows that
    $\rho(z_{1})=\rho(z_{2})$. Hence, $\rho'$ can be constructed by setting $\rho'(z') = \rho(z_1)$.
  \end{Proof}
  Now we proceed to prove claim~(\ref{p:minimalfansconfig1}) of
  Lemma~\ref{p:minimalfansconfig}. Let $\Gbf=\set{O_{i};m_{ij}}$ be a
  diagram category, $\Xcal,\Ycal,\Zcal\in\Prob\<\Gbf>$ be three
  $\Gbf$-configurations and  
  $\Fcal=(\Xcal\ot\Zcal\to\Ycal)$ be a two-fan.
  Recall that it can also be considered as a $\Gbf$-configuration of
  two-fans 
  \[
  \Fcal=\set{\Fcal_{i};f_{ij}}
  \]
  Any minimizing reduction
  \[
  \Fcal\!=\!(\Xcal\!\ot\!\Zcal\!\to\!\Ycal)
  \too
  \Fcal'\!=\!(\Xcal\!\ot\!\Zcal'\!\to\!\Ycal)
  \]
  induces reductions 
  \[
  \Fcal_{i}\!=\!(X_{i}\!\ot\! Z_{i}\!\to\! Y_{i})
  \too
  \Fcal_{i}\!=\!(X_{i}\!\ot\! Z_{i}'\!\to\! Y_{i})
  \]
  for all $i$ in the index set $I$.
  It follows that if all $\Fcal_{i}$'s are minimal, then so is
  $\Fcal$. 

  Now we prove the implication in the other direction. Suppose $\Fcal$ is
  minimal. We have to show that all $\Fcal_{i}$ are minimal as well.
  Suppose there exist a non-minimal fan among $\Fcal_{i}$'s.  For an
  index $i \in I$ let
  \begin{align*}
    \check J(i)
    &:=
    \set{j\in I \st\Hom_{\Gbf}(O_{j},O_{i})\neq\emptyset}
    \\
    \hat J(i)
    &:=
    \set{j\in I \st\Hom_{\Gbf}(O_{i},O_{j})\neq\emptyset}
  \end{align*}

  Choose an index $i_{0}$ such that 
  \begin{enumerate}
  	\item  $\Fcal_{i_{0}}$ is not minimal
  	\item for any $j\in \hat J(i_{0}) \backslash \{ i_0 \}$ the two-fan $\Fcal_{j}$ is
  minimal. 
  \end{enumerate}
  Consider now the minimal reduction
  $\mu:\Fcal_{i_{0}}\to\Fcal'_{i_{0}}$ and construct a two-fan
  $\Gcal=\set{\Gcal_{i};g_{ij}}$ of $\Gbf$-configurations by setting
  \[
  \Gcal_{i}
  :=
  \begin{cases}
    \Fcal'_{i}
    &
    \text{if $i=i_{0}$}
    \\
    \Fcal_{i}
    &
    \text{otherwise}
  \end{cases}
  \]
  and 
  \[
  g_{ij}
  :=
  \begin{cases}
    \mu\circ f_{ij}
    &
    \text{if $j=i_{0}$ and $i\in\check J(i_{0})$}
    \\
    f'_{ij}
    &
    \text{if $i=i_{0}$ and $j\in\hat J(i_{0})$}
    \\
    f_{ij}
    &
    \text{otherwise}
  \end{cases}
  \]
  where $f'_{i_{0}j}$ is the reduction provided by the
  Lemma~\ref{p:minimising-reduction} applied to the diagram
  \[
  \begin{tikzcd}[ampersand replacement=\&]
    \Fcal_{i_{0}}
    \arrow{d}{\mu}
    \arrow{dr}{f_{i_{0}j}}
    \\
    \Fcal'_{i_{0}}
    \arrow[dashrightarrow]{r}{f'_{i_{0}j}}
    \&
    \Fcal_{j}
  \end{tikzcd}
  \]
  We thus constructed a non-trivial reduction $\Fcal\to\Gcal$ which is
  identity on the terminal  $\Gbf$-configurations $\Xcal$ and
  $\Ycal$. This contradicts the minimality of $\Fcal$.
  
  To address the second assertion of the
  Lemma~\ref{p:minimalfansconfig} observe that the argument above
  gives an algorithm for the construction of a minimal reduction of any
  two-fan of $\Gbf$-configurations.
\end{Proof}



\def\thesubsection{\thesection.\ref{s:kolmogorov}}
\subsection{Statements from the section ``Kolmogorov distance''}
\repeatclaim{p:kolmogorovisdistance}

\begin{Proof}
  The symmetry of $\ikd$ is immediate. The non-negativity of $\ikd$
  follows from the fact that entropy of the target space of a reduction
  is not greater then the entropy of the domain, which is a
  particular instance of the Shannon inequality~(\ref{eq:shannonineq}).

  We proceed to prove the triangle inequality. We will make use of the
  following lemma
  \begin{tlemma}{p:kdtriangle}
    For a minimal full configuration of probability spaces
    \[
    \<X,Y,Z>=
    \left(
    \begin{tikzcd}[row sep=small,column sep=tiny,ampersand replacement=\&]
      \mbox{}
      \&
      XYZ
      \arrow{dl}{}
      \arrow{d}{}
      \arrow{dr}{}
      \&
      \\
      XY
      \arrow{d}{}
      \&
      XZ
      \arrow{dr}{}
      \arrow{dl}{}
      \&
      YZ
      \arrow{d}{}
      \\
      X
      \&
      Y
      \arrow[leftarrow,crossing over]{ur}{}
      \arrow[leftarrow,crossing over]{ul}{}
      \&
      Z
    \end{tikzcd}
    \right)
    \] 
    holds
    \[
    \kd(X\ot XZ\to Z)
    \leq
    \kd(X\ot XY\to Y)
    +
    \kd(Y\ot YZ\to Z)
    \]
  \end{tlemma}
  
  \begin{Proof}
    By Shannon inequality~(\ref{eq:shannonineq}) on page
    \pageref{eq:shannonineq} we have
    \[ 
    \ent( X \rel Z ) 
    \leq 
    \ent(XY \rel Z) 
    \leq 
    \ent( X \rel Y ) + \ent (Y \rel Z) 
    \] 
    Similarly, 
    \[ 
    \ent(Z \rel X ) 
    \leq 
    \ent(Z \rel Y) + \ent (Y \rel X) 
    \] 
    and therefore 
    \[ 
    \kd(X\ot XZ\to Z) 
    \leq 
    \kd(X\ot XY\to Y) 
    + 
    \kd(Y\ot YZ\to Z)  
    \]
  \end{Proof}

  Now we continue with the proof of
  Proposition~\ref{p:kolmogorovisdistance}.

  Let $\Gbf$ be an arbitrary complete reduction category.
  Suppose $\Xcal = \set{X_i ; f_{ij}}$, $\Ycal = \set{Y_i; g_{ij}}$
  and $\Zcal = \set{Z_i; h_{ij}}$ are $\Gbf$-configurations, with
  initial spaces being $X_{0}$, $Y_{0}$ and $Z_{0}$, respectively. Let
  \begin{align*}
    \hat\Fcal
    &=
    (\Xcal\ot\Fcal\to\Ycal)
    \\
    \hat\Gcal
    &=
    (\Ycal\ot\Gcal\to\Zcal)
  \end{align*}
  be two optimal minimal two-fans
  satisfying
  \begin{align*} 
    \ikd(\Xcal, \Ycal)
    &= 
    \kd(\hat\Fcal) 
    \\ 
    \ikd(\Ycal,\Zcal)
    &= 
    \kd(\hat\Gcal) 
  \end{align*}
   
  Recall that each two-fan of $\Gbf$-configurations is a
  $\Gbf$-configuration of two-fans between the individual spaces, that
  is 
  \[ 
  \begin{split} 
    \Fcal 
    &= 
    \set{\rule{0em}{2ex} \Fcal_i = (X_i \ot F_i \to Y_i)} 
    \\ 
    \Gcal 
    &= 
    \set{ \rule{0em}{2ex} \Gcal_i = (Y_i \ot G_i \to Z_i)} 
  \end{split} 
  \] 

  We construct a coupling $\Hcal$ between $\Xcal$ and $\Zcal$ in the following
  manner. Starting with the two-tents configuration between the
  initial spaces, we use adhesion to extend it to a full
  configuration, thus constructing a coupling between $X_{0}$ and
  $Z_{0}$. This full configuration could then be ``pushed down'' and
  provides full extensions of two-tents on all lower levels. Thus we
  could ``compose'' couplings $\Fcal$ and $\Gcal$ and use a Shannon
  inequality to establish the triangle inequality for the Kolmogorov
  distance. Details are as follows. 
  
  Consider a two-tents configuration
  \[
  X_{0}\ot F_{0}\to Y_{0}\ot G_{0}\to Z_{0}
  \]
  and extend it by adhesion, as described in Section
  \ref{s:config-adhesion} to a $\Lambdabf_{3}$-configuration
  \[
  \begin{tikzcd}[row sep=small, column sep=small, ampersand replacement=\&]
    \mbox{}
    \&
    A_{0}
    \arrow{dl}{}
    \arrow{d}{}
    \arrow{dr}{}
    \\
    F_{0}
    \arrow{d}
    \&
    H_{0}
    \arrow{dl}
    \arrow{dr}  
    \&
    G_{0}
    \arrow{d}
    \\
    X_{0}
    \&
    Y_{0}
    \arrow[leftarrow, crossing over]{ul}
    \arrow[leftarrow, crossing over]{ur}
    \&
    Z_{0}
  \end{tikzcd}
  \]
  Together with the reductions 
  \begin{align*}
  (X_{0})^{\Gbf}
    &\to 
    \Xcal
    \\
  (Y_{0})^{\Gbf}
    &\to 
    \Ycal
    \\
  (Z_{0})^{\Gbf}
    &\to 
    \Zcal
  \end{align*}
  it gives rise to a $\Lambdabf_{3}$-configuration of
  $\Gbf$-configurations
  \[\tageq{bigfull}
  \begin{tikzcd}[row sep=small, column sep=small, ampersand replacement=\&]
    \mbox{}
    \&
    (A_{0})^{\Gbf}
    \arrow{dl}{}
    \arrow{d}{}
    \arrow{dr}{}
    \\
    (F_{0})^{\Gbf}
    \arrow{d}
    \&
    (H_{0})^{\Gbf}
    \arrow{dl}
    \arrow{dr}  
    \&
    (G_{0})^{\Gbf}
    \arrow{d}
    \\
    \Xcal
    \&
    \Ycal
    \arrow[leftarrow, crossing over]{ul}
    \arrow[leftarrow, crossing over]{ur}
    \&
    \Zcal
  \end{tikzcd}
  \]

  Note that the minimal reductions of the two-fan subconfigurations of
  (\ref{eq:bigfull})
  \begin{align*}
    \Xcal\ot(F_{0})^{\Gbf}\to\Ycal\\
    \Ycal\ot(G_{0})^{\Gbf}\to\Zcal
  \end{align*}
  are the two-fans $\hat\Fcal$ and $\hat\Gcal$, respectively, by
  Lemma~\ref{p:minimalfansconfig}.

  Now consider the ``minimization'' of the above configuration
  \[\tageq{minfull}
  \begin{tikzcd}[row sep=small,column sep=small, ampersand replacement=\&]
    \mbox{}
    \&
    \Acal
    \arrow{dl}{}
    \arrow{d}{}
    \arrow{dr}{}
    \\
    \Fcal
    \arrow{d}
    \&
    \Hcal
    \arrow{dl}
    \arrow{dr}  
    \&
    \Gcal
    \arrow{d}
    \\
    \Xcal
    \&
    \Ycal
    \arrow[leftarrow, crossing over]{ul}
    \arrow[leftarrow, crossing over]{ur}
    \&
    \Zcal
  \end{tikzcd}
  \]
  It could also be viewed as a $\Gbf$-configuration of $\Lambdabf_{3}$
  configurations, 
  \[
  \begin{tikzcd}[row sep=small, column sep=small, ampersand replacement=\&]
    \mbox{}
    \&
    A_{i}
    \arrow{dl}{}
    \arrow{d}{}
    \arrow{dr}{}
    \\
    F_{i}
    \arrow{d}
    \&
    H_{i}
    \arrow{dl}
    \arrow{dr}  
    \&
    G_{i}
    \arrow{d}
    \\
    X_{i}
    \&
    Y_{i}
    \arrow[leftarrow, crossing over]{ul}
    \arrow[leftarrow, crossing over]{ur}
    \&
    Z_{i}
  \end{tikzcd}  
  \]
  each of which is minimal by Corollary \ref{p:minimalfullconfig}.
  
  Now we can apply Lemma~\ref{p:kdtriangle}
  to each level to conclude that
  \begin{align*}
    \ikd(\Xcal,\Zcal)
    &\leq
    \kd(\Xcal\ot\Hcal\to\Zcal)
    \\
    &\leq
    \kd(\Xcal\ot\Fcal\to\Ycal)
    +
    \kd(\Xcal\ot\Gcal\to\Ycal)
    \\
    &=
    \ikd(\Xcal,\Ycal)
    +
    \ikd(\Ycal,\Zcal)
  \end{align*}

  Finally, if $k(\Xcal, \Ycal) = 0$, then there is a two-fan $\Fcal$
  of $\Gbf$-configurations between $\Xcal$ and $\Ycal$ with
  $\kd(\Fcal)= 0$, from which it follows that $\Xcal$ and $\Ycal$ are
  isomorphic.
\end{Proof}

\repeatclaim{p:tensor1lip}
\begin{Proof}
  The claim follows easily from the additivity of
  entropy in equation (\ref{eq:entropyadditive}). Suppose that
  $\Xcal=\set{X_{i};f_{ij}}$, $\Ycal=\set{Y_{i};g_{ij}}$ and
  $\Ycal'=\set{Y'_{i};g'_{ij}}$ are three $\Gbf$-configurations and
  \[
  \Fcal
  =
  (\Ycal\ot \Zcal\to \Ycal')
  \]
  is an optimal fan, so that
  \[
  \ikd(\Ycal,\Ycal')
  =
  \sum_{i} \big[2\ent(Z_{i})-\ent(Y_{i})-\ent(Y'_{i})\big]
  \]
  Consider the fan
  \[
  \Gcal
  =
  (\Xcal\otimes \Ycal\ot \Xcal\otimes \Zcal\to \Xcal\otimes \Ycal')
  \]
  Then, by additivity of entropy, in equation (\ref{eq:entropyadditive}), we have
  \begin{align*}
    \kd(\Fcal)
    &=
    \sum_{i}\big[2\ent(X_{i}\otimes Z_{i})-\ent(X_{i}\otimes Y_{i}) -
      \ent(X_{i}\otimes Y'_{i})\big]
    \\
    &=
    \sum_{i}\big[2\ent(Z_{i})-\ent(Y_{i}) -
                 \ent(Y'_{i})\big]
    \\
    &=
    \kd(\Gcal)
  \end{align*}
  and, therefore,
  \begin{align*}
    \ikd(\Xcal\otimes \Ycal,\Xcal\otimes \Ycal')
    &\leq
    \kd(\Gcal)=\kd(\Fcal)=\ikd(\Ycal,\Ycal')
  \end{align*}
  Thus, the tensor product of probability spaces is 1-Lipschitz with
  respect to each argument.
\end{Proof}
  
\repeatclaim{p:entropy1lip}
\begin{Proof}

  Let $\Xcal,\Ycal\in\prob\<\Gbf>$ and let
  \[
  \Gcal=(\Xcal\ot\Zcal\to\Ycal)
  \]
  be an optimal fan with components
  \[
  \Gcal_{i}=(X_{i}\ot Z_{i}\to Y_{i})
  \]
  
  For a fixed index $i$ we can estimate
  the difference of entropies
  \[
    \ent(X_{i})-\ent(Y_{i})
    =
    2\big(\ent(X_{i})-\ent(Z_{i})\big) + \kd(\Gcal_{i})
    \leq
    \kd(\Gcal_{i})
  \]
  By symmetry we then have
  \[
  |\ent(X_{i})-\ent(Y_{i})|
  \leq
  \kd(\Gcal_{i})    
  \]

  Adding above inequalities for all $i$ we have
  \[
  |\ent_{*}(\Xcal)-\ent_{*}(\Ycal)|_{1}
  \leq
  \kd(\Gcal)=\ikd(\Xcal,\Ycal)
  \]
  By the additivity of entropy we also obtain the $1$-Lipschitz
  property of the entropy function with respect to the asymptotic
  Kolmogorov distance $\aikd$.
\end{Proof}

\repeatclaim{p:restriction1lip}
\begin{Proof}
  The claim  follows from the functoriality of the restriction
  operator. We argue as follows.

  Suppose that $R:\Gbf'\to\Gbf$ is a functor and $R^{*}$ is the
  corresponding restriction operator. For
  $\Xcal_{1},\Xcal_{2}\in\prob\<\Gbf>$ let
  \[
  \Fcal=(\Xcal_{1}\ot\Ycal\to\Xcal_{2})
  \]
  be an optimal fan. Then
  \[
  \Fcal':=(R^{*}\Xcal_{1}\ot R^{*}\Ycal\to R^{*}\Xcal_{2})
  \]
  is a fan with the terminal vertices being the restrictions of
  $\Xcal_{1}$ and $\Xcal_{2}$. It can be considered as a
  $\Gbf'$-configuration of two-fans over individual spaces in
  $R^{*}\Xcal_{1}$ and $R^{*}\Xcal_{2}$ each of which also appears as
  a fan in $\Fcal$. 
  
  Thus, we obtain the rough estimate
  \[
  \ikd(R^{*}\Xcal_{1},R^{*}\Xcal_{2})
  \leq
  \size{\Gbf'}\cdot\ikd(\Xcal_{1},\Xcal_{2})
  \]
  Since the restriction operator commutes with tensor powers, the same
  estimate also holds for the asymptotic Kolmogorov distance $\aikd$.
\end{Proof}


\repeatclaim{p:slicing}
\begin{Proof}
  Since the two-fan $(U\ot W\to V)$ is minimal the probability space
  $W$ could be considered having underlying set to be a subset of the
  Cartesian product of the underlying sets of $U$ and $V$. For any
  pair $(u,v)\in \un W$ with a positive weight consider an optimal
  two-fan
  \[\tageq{opt4slices}
  \Gcal_{uv}
  =
  (\Xcal\rel u
  \stackrel{\pi_{\Xcal}}{\oot}
  \Zcal_{uv}
  \stackrel{\pi_{Y}}{\too}
  \Ycal\rel v)
  \]
  where $\Zcal_{uv}=\set{Z_{uv,i};\rho_{ij}}$.  Let $p_{uv,i}$ be the
  probability distributions on $Z_{uv,i}$ -- the individual spaces in
  the configuration $\Zcal_{uv}$.  The next step is to take a convex
  combination of distributions $p_{uv,i}$ weighted by $p_{W}$ to
  construct a coupling $\Xcal\ot\Zcal\to\Ycal$.
  
  First we extend the 7-vertex configuration to a full
  $\Lambdabf_{4}$-configuration of $\Gbf$-configurations, such
  that the top vertex has the distribution
  \[\tageq{definitiondbinslicing}
  p_{i}(x,y,u,v):=p_{uv,i}(x,y)p_{W}(u,v)
  \]
  as described in the Section~\ref{s:config-examples-full}. 
  
  If we integrate over $y$, we obtain
  \[
  \sum_{y} p_{i}(x,u,v,y) = ((\pi_{\Xcal,i})_{*}p_{uv,i})(x) p_{W}(u,v)
  \]
  Then we use that by (\ref{eq:opt4slices}) it holds that $
  (\pi_{\Xcal,i})_{*}p_{uv,i}=p_{X_i}(\,\cdot\,\rel u)$ and therefore
  \[
  \sum_y p_i(x,y,u,v) = p_{X_i}(x \rel u) p_{W}(u,v).
  \]
  In the same way,
  \[
  \sum_x p_i(x,y,u,v) = p_{Y_i}(y \rel v) p_{W}(u,v).
  \]
  Note that this exactly corresponds to adhesion, as described in
  Section \ref{s:config-adhesion}.  It follows that
  \[\tageq{reluvisrelu}
  \Xcal\rel uv = \Xcal \rel u 
  \quad \text{and} \quad 
  \Ycal\rel uv = \Ycal \rel v
  \]
  and
  \[\tageq{slicingentad}
  \ent(X_{i}\rel UV)=\ent(X_{i}\rel U)
  \quad \text{and} \quad 
  \ent(Y_{i}\rel UV)=\ent(Y_{i}\rel V)
  \]
  
  The extended configuration contains a two-fan of configurations 
  $\Fcal=(\Xcal\ot\Zcal\to\Ycal)$ with
  terminal vertices $\Xcal$ and $\Ycal$. We call its initial vertex
  $\Zcal=\set{XY_{i},f_{ij}}$.
  
  The following estimates conclude the proof the the Slicing Lemma.
  First we use the definitions of intrinsic Kolmogorov distance $\ikd$
  and of $\kd(\Fcal)$ to estimate
  \begin{align*}
    \ikd(\Xcal,\Ycal)
    &\leq
    \kd(\Fcal)\\
    &=
    \sum_{i}\kd(\Fcal_{i})\\
    &=
    \sum_{i}\big[
      2\ent(XY_{i}) - \ent(X_{i}) - \ent(Y_{i})
      \big]
  \end{align*}
  Next, we apply the definition of the conditional entropy to rewrite
  the right-hand side
  \begin{align*}
    \ikd(\Xcal, \Ycal)
    &\leq  
    \sum_{i}\big[
      2\ent(XY_{i}\rel UV)+2\ent(UV)-2\ent(UV\rel XY_{i})\\
      &\quad\quad\;
      -\ent(X_{i}\rel U)-\ent(U)+\ent(U\rel X_{i})\\
      &\quad\quad\;
      -\ent(Y_{i}\rel V)-\ent(V)+\ent(V\rel Y_{i})\big]
  \end{align*}
  We now use (\ref{eq:slicingentad}) and rearrange terms to obtain
  \begin{align*}
    \ikd(\Xcal, \Ycal)
    &\leq
    \sum_{i}\big[
      2\ent(XY_{i}\rel UV)-\ent(X_{i}\rel UV)-\ent(Y_{i}\rel UV)\\ 
      &\quad\quad\;
      + 2\ent(UV) - \ent(U) - \ent(V)\\
      &\quad\quad\;
      - 2\ent(UV\rel XY_{i})
      +\ent(U\rel X_{i})
      +\ent(V\rel Y_{i})\big]
  \end{align*}
  By the integral formula for conditional
  entropy (\ref{eq:relentinteg}) applied to the first
  three terms we get
  \begin{align*}
    \sum_i 
    \big[
      2\ent(XY_{i}\rel UV)&-\ent(X_{i}\rel UV)-\ent(Y_{i}\rel UV)
    \big]
    \\
    &=
    \int_{UV}\ikd(\Xcal\rel{uv},\Ycal\rel uv)\d p_{W}(u,v)
  \end{align*}
  However, because of (\ref{eq:reluvisrelu}) this simplifies to
\[
\int_{UV}\ikd(\Xcal\rel{uv},\Ycal\rel uv)\d p_{W}(u,v)
= \int_{UV}\ikd(\Xcal\rel{u},\Ycal\rel v)\d p_{W}(u,v)
\]
Therefore,
  \begin{align*}
  \ikd(\Xcal, \Ycal)
  &\leq \int_{UV}\ikd(\Xcal\rel{u},\Ycal\rel v)\d p_{W}(u,v)+
  \size{\Gbf}\cdot\kd(U \ot W \to V) \\
  &\quad+ 
  \sum_{i}\big[\ent(U\rel X_{i})+\ent(V\rel Y_{i})\big]
  \end{align*}
  
\end{Proof}


\repeatclaim{p:kolmogorovlocal}
\begin{Proof}
  We will need the following obvious rough estimate of the Kolmogorov
  distance that holds for any $p,q\in\Delta\Scal$:
  \[\tageq{roughlocalestimate}
  2 \cdot \ikd(\Xcal,\Ycal)\leq 2\size{\Gbf}\cdot\ln|S_{0}|
  \]
  It can be obtained by taking a tensor product for the coupling
  between $\Xcal$ and $\Ycal$.

  Our goal now is to write $p$ and $q$ as the convex combination of
  three other distributions $\hat p$, $p^{+}$ and $q^{+}$ as in
  \begin{align*}
    p
    &=
    (1-\alpha)\cdot\hat p + \alpha\cdot p^{+}
    \\
    q
    &=
    (1-\alpha)\cdot\hat p + \alpha\cdot q^{+}
  \end{align*}
  with the smallest possible $\alpha\in[0,1]$.

  We could do it the following way.  Let
  $\alpha:=\frac12|p_{0}-q_{0}|$. If $\alpha=1$ then the proposition
  follows from the rough estimate~(\ref{eq:roughlocalestimate}), so
  from now on we assume that $\alpha<1$.  Define three probability
  distributions $\hat p_{0}$, $p_{0}^{+}$ and $q_{0}^{+}$ on $S_{0}$
  by setting for every $x\in S_{0}$
  \begin{align*}
    \hat p_{0}(x) 
    &:= \frac1{1-\alpha}
    \min\set{p_{0}(x),q_{0}(x)}
    \\ 
    p_{0}^{+} 
    &:=
    \frac{1}{\alpha}\big(p_{0}-(1-\alpha)\hat p_{0}\big)
    \\ 
    q_{0}^{+} 
    &:= 
    \frac{1}{\alpha}\big(q_{0}-(1-\alpha)\hat p_{0}\big)
  \end{align*}
  
  Denote by $\hat p,p^{+},q^{+}\in\Delta\Scal$ the distributions
  corresponding to $\hat p_{0},p_{0}^{+},q_{0}^{+}\in\Delta S_{0}$
  under isomorphism~(\ref{eq:distribonconfig}). Thus we have
  \begin{align*}
    p&=(1-\alpha)\hat p+\alpha\cdot p^{+}\\
    q&=(1-\alpha)\hat p+\alpha\cdot q^{+}
  \end{align*}

  Now we construct a ``two-tents'' configuration of
  $\Gbf$-configurations 
  \[\tageq{2tentsforlocal}
  \Xcal\ot\tilde\Xcal\to\Lambda_{\alpha}\ot\tilde\Ycal\to\Ycal
  \]
  by setting 
  \begin{align*}
    \tilde X_{i}
    &:=
    \Big(S_{i}\times\un\Lambda_{\alpha};\;
         \tilde p_{i}(s,\square)=(1-\alpha)\hat p_{i}(s),\,
         \tilde p_{i}(s,\blacksquare)=\alpha\cdot p^{+}_{i}(s)
    \Big)
    \\
    \tilde Y_{i}
    &:=
    \Big(S_{i}\times\un\Lambda_{\alpha};\;
         \tilde q_{i}(s,\square)=(1-\alpha)\hat p_{i}(s),\,
         \tilde q_{i}(s,\blacksquare)=\alpha\cdot q^{+}_{i}(s)
    \Big)
    \\
  \end{align*}
  and
  \begin{align*}
    \tilde\Xcal
    &:=
    \set{\tilde X_{i};\,f_{ij}\times\Id}\\
    \tilde\Ycal
    &:=
    \set{\tilde Y_{i};\,f_{ij}\times\Id}\\
  \end{align*}
  The reductions in the ``two-tents'' sub-configurations
  of~(\ref{eq:2tentsforlocal}) are given by coordinate projections.  Note
  that the following isomorphisms hold
  \begin{align*}
    \Xcal\rel\square
    &\cong
    (\Scal,\hat p)
    \\
    \Xcal\rel\blacksquare
    &\cong
    (\Scal, p^{+})
    \\
    \Ycal\rel\square
    &\cong
    (\Scal, \hat p)
    \cong
    \Xcal\rel\square
    \\
    \Ycal\rel\blacksquare
    &\cong
    (\Scal,q^{+})
  \end{align*}

  Now we apply part (\ref{i:slicing2tents}) of Corollary
  \ref{p:slicingcorollary} to obtain the desired inequality
  \begin{align*}
    \ikd(\Xcal,\Ycal)
    &\leq
    (1-\alpha)\ikd(\Xcal\rel\square,\Ycal\rel\square)
    +
    \alpha\cdot\ikd(\Xcal\rel\blacksquare,\Ycal\rel\blacksquare)
    \\
    &\quad
    +\sum_{i}\big[\ent(\Lambda_{\alpha}\rel X_{i})
      +\ent(\Lambda_{\alpha}\rel Y_{i})\big]
    \\
    &\leq
    2\cdot\size{\Gbf}\cdot
    \big(\alpha\cdot\ln|S_{0}|+\ent(\Lambda_{\alpha})\big)
  \end{align*}
\end{Proof}

\def\thesubsection{\thesection.\ref{s:extensions}}
\subsection{Statements from the section ``Extensions''}
\repeatclaim{p:extensionlemma}
\begin{Proof}
  Denote by $X_{0}$, $X_{0}'$ the initial spaces in the
  configurations $\Xcal$, $\Xcal'$, respectively.
  Let $Y_{0}$ be the initial space of the full sub-configuration of
  $\Ycal$ generated by $Y_{k+1},\ldots,Y_{l}$. Let $K_{0}$ be the
  initial space in the optimal coupling between $\Xcal$ and $\Xcal'$
  \[
  \Fcal=(\Xcal'\ot\Kcal\to\Xcal)
  \]
  Recall that $X_{0}$ could be considered as the Cartesian product of
  the underlying sets of spaces generating $\Xcal$ with some
  distribution on it. A similar view holds for $X_{0}'$, $Y_{0}$ and
  $K_{0}$.
  Thus we have in particular $\un K_{0}=\un X_{0}'\times\un X_{0}$
  
  Define a full minimal configuration
  $\Zcal\in\prob\<\Lambdabf_{2k+l}>$ by providing a distribution on
  \[
  \un X_{0}'\times\un X_{0}\times\un Y_{0}=\un K_{0}\times\un Y_{0}
  \]
  as explained in the Section \ref{s:config-examples-full}.
  The distribution will be defined by
  \[
  p(\xbf',\xbf,\ybf)
  :=
  p_{\Fcal}(\xbf',\xbf)\cdot p_{\Ycal}(\xbf,\ybf)/p_{\Xcal'}(\xbf)
  \]
  
  It is clear that $\Zcal$ contains both the coupling $\Fcal$ and
  configuration $\Ycal$ as restrictions. It also contains the minimal full
  configuration 
  \[
  \Ycal'=\<X_{1}',\ldots,X_{k}',Y_{k+1},\ldots,Y_{l}>
  \]
  and a coupling $\Gcal$ between $\Ycal$ and $\Ycal'$. 
  
  For a pair of spaces $A$ and $B$ in $\Zcal$ we denote by $AB$
  the initial space of a minimal fan in $\Zcal$ with the terminal
  spaces $A$ and $B$. 
  The two-fan of $\Lambdabf_{k+l}$-con\-fi\-gu\-rations $\Gcal$ can be considered as a $\Lambdabf_{k+l}$-configuration
  of two-fans
  \[
  \Gcal_{IJ}:=(X_{I}Y_{J}\oot G_{IJ}\too X'_{I}Y_{J})
  \]
  Using this notation we estimate for $I\subset\set{1,\ldots,k}$ and
  	$J\subset\set{k+1,\ldots,l}$
  \begin{align*}
    \ikd(\Ycal,\Ycal')
    &\leq
    \kd{\Gcal}
    =
    \sum_{I,J}\kd(\Gcal_{IJ})
    \\
    &=
    \sum_{I,J}
      \big[
        2\ent(G_{IJ})-\ent(X_{I}Y_{J})-\ent(X_{I}'Y_{J})
      \big]
    \\
    &\leq
    \sum_{I,J}
    \Big[
      2\ent(X_{I}X'_{I})-\ent(X_{I})-\ent(X_{I}')+\\
      &\quad\quad+
      \big(
        2\ent(Y_{J}\rel X_{I}X_{I}')-\ent(Y_{J}\rel
        X_{I})-\ent(Y_{J}\rel X_{I}')
      \big)
    \Big]
    \\
    &\leq
    2^{k-l}
    \sum_{I}\kd(\Fcal_{I})
    \\
    &\leq
    2^{k-l}\kd(\Fcal)
    =
    2^{k-l}\ikd(\Xcal,\Xcal')
  \end{align*}
\end{Proof}

\repeatclaim{p:stabilized-rel-ent-set-Lipschitz}
\begin{Proof}
  Note that by Corollary \ref{p:unstabilized-rel-ent-set-Lipschitz},
  or more directly by Proposition \ref{p:extensionlemma}, for $n \in
  \Nbb$
  \[
  \frac1n \d_H 
  \left(
    \res_l\big((\Xcal)^{\otimes n}\big), 
    \res_l \big((\Xcal')^{\otimes n}\big)
  \right)
  \leq 
  2^{l-k}\frac1n \ikd \big( (\Xcal)^{\otimes n}, (\Xcal')^{\otimes n} \big)
  \]	
  Hence, by the scaling properties of the Hausdorff distance
  \[
  \d_H 
  \left(
    \frac{1}{n}\res_l\big((\Xcal)^{\otimes n}\big), 
    \frac{1}{n}\res_l \big((\Xcal')^{\otimes n}\big)
  \right)
  \leq 
  2^{l-k} \frac1n \ikd \big( (\Xcal)^{\otimes n}, (\Xcal')^{\otimes n} \big)
  \]
  For convenience, we introduce the notation
  \begin{align*}
    K_n 
    &= 
    \closure\left( \frac{1}{n}\res_l\big((\Xcal)^{\otimes n}\big) \right) 
    & 
    K
    &= 
    \res_l(\Xcal) 
    \\
    K_n' 
    &= 
    \closure
    \left( 
       \frac{1}{n}\res_l\big((\Xcal')^{\otimes n}\big) 
    \right) 
    &
    K'
    &= 
    \sres_l(\Xcal')
  \end{align*}
  Recall that by definition,
  \[
  K 
  = 
  \closure \left(\bigcup_{n \in \Nbb} K_n\right) 
  \qquad 
  K'  
  = 
  \closure \left( \bigcup_{n \in \Nbb} K_n'\right)
  \]

  Note that by the superadditivity property of the unstabilized
  relative entropic sets (see inclusion
  (\ref{eq:superadditivity-res})) the sequences $n \mapsto K_{n!}$ and
  $n \mapsto K_{n!}'$ are monotonically increasing sequences of sets,
  and
  \[
  \bigcup_{i =1}^n K_i \subset K_{n!} 
  \qquad 
  \bigcup_{i =1}^n K'_i \subset K'_{n!}
  \]

  Now select a large radius $R > 0$. Let $B_R(0)$ denote the ball of
  radius $R$ around the origin in $\Rbb^{2^{\{1,\ldots,l\}}}$.  By
  compactness and the definition of the stabilized relative entropic
  set
  \begin{align*}
    \d_H( K_{n!} \cap B_R(0) , K \cap B_R(0) ) 
    &\to 
    0
    \\
    \d_H( K_{n!}' \cap B_R(0), K' \cap B_R(0) ) 
    & \to 
    0
  \end{align*}
  as $n \to \infty$. Therefore also
  \[
  \d_H
  \big( 
    K \cap B_R(0) , 
    K' \cap B_R(0) 
  \big) 
  \leq 
  2^{l-k} \aikd\big( \Xcal, \Xcal' \big)
  \]
  Because this inequality holds for every $R>0$, the estimate in the
  lemma follows.
\end{Proof}

\def\thesubsection{\thesection.\ref{s:mixtures}}
\subsection{Statements from the section ``Mixtures''}
\repeatclaim{p:mixtures-dist1}
\begin{Proof}
  Recall that for the empirical reduction
  \[
  \emp:\Lambda_{1/n}^{\otimes N}\to\Delta\Lambda_{1/n}
  \]
  the quantity $N\cdot\emp(\lambda)(\blacksquare)$ counts the number of
  black squares in the sequence $\lambda$. It is a binomially
  distributed random variable with the mean $N/n$ and variance
  $\frac{N}{n}(1-\frac1n)$.
  
  The first claim is then proven by the following calculation
  \begin{align*}
    \aikd(\Xcal&,\Xcal^{\otimes n}\oplus_{\Lambda_{1/n}}\set{\bullet})
    \\
    &=
    \lim_{N\to\infty}
      \frac1N
        \ikd\left(
          \Xcal^{\otimes N},
          (\Xcal^{\otimes n}\oplus_{\Lambda_{1/n}}\set{\bullet})^{\otimes N}
        \right)
    \\
    &=
    \lim_{N\to\infty}
      \frac1N
        \ikd\left(
          \Xcal^{\otimes N},
          \bigoplus_{\lambda\in\Lambda_{1/n}^{\otimes N}}
              \Xcal^{\otimes n\cdot N\cdot \emp(\lambda)(\blacksquare)}
        \right)
    \\
    &\leq
    \ent(\Lambda_{1/n}) +
    \lim_{N\to\infty}
      \frac1N
      \int_{\lambda\in\Lambda^{\otimes n}_{1/n}}
        \ikd(\Xcal^{\otimes N},
             \Xcal^{\otimes(N\cdot n\cdot \emp(\lambda)(\blacksquare))})
      \d p(\lambda)
    \\
    &\leq
    \ent(\Lambda_{1/n}) +
    |\ent_{*}(\Xcal)|_{1}\cdot
    \lim_{N\to\infty}
      \frac{n}{N}
      \cdot
      \int_{\lambda\in\Lambda_{1/n}^{\otimes N}}
         \big|N/n- N\cdot \emp(\lambda)(\blacksquare)\big|
         \d p(\lambda)
    \\
    &\leq
    \ent(\Lambda_{1/n}) +
    |\ent_{*}(\Xcal)|_{1}\cdot
    \lim_{N\to\infty}
      \frac{n}{N}\cdot\sqrt{N\cdot\frac1n(1-\frac1n)}
    \\
    &=
    \ent(\Lambda_{1/n})
  \end{align*}

  The second claim is proven similarly and the third follows from the
  second and the $1$-Lipschitz property of the tensor
  product. Finally, the fourth follows from
  Corollary~\ref{p:slicingcorollary}(\ref{p:slicingcofan}), by slicing
  both arguments along $\Lambda_{1/n}$.
\end{Proof}

\def\thesubsection{\thesection.\ref{s:tropical}}
\subsection{Statements from the section ``Tropical Probability''}
\begin{tlemma}{p:subadditive}
  Suppose the sequence $\set{a(i)}_{i\in\Nbb_{0}}$ of  real numbers is
  bounded from below and is
  quasi-subadditive, that is there is a constant $C\in\Rbb$ such that
  for any $i,j\in\Nbb_{0}$ holds
  \[
  a(i+j)\leq a(i) + a(j) + C
  \]
  Then the limit
  \[
  \lim_{i\to\infty}\frac1i a(i)
  \]
  exists and is finite.
\end{tlemma}

\begin{Proof}
  The lemma is standard and is sometimes refered to as Fekete's subadditive lemma. We include a proof for the convenience of the reader.
  Assume first that $C=0$.  Then the sequence satisfies $a(k\cdot
  i)\leq k\cdot a(i)$ and in particular $a(i)\leq i\cdot a(1)$.  Let
  $l:=\liminf\frac1i a(i)\in[0,\infty)$.  Choose $\epsilon>0$. Then we can
  find $k\in\Nbb$ such that $\frac1k a(k)\leq l+\epsilon$.  For
  $n\in\Nbb$ let $q,r$ be the quotient and the reminder of the integer
  division of $n$ by $k$, that is
  \[
  n=q\cdot k+r,\quad0\leq r<k
  \]
  Then
  \begin{align*}
    \frac1n a(n)
    &\leq
    \frac1n (q\cdot a(k)+a(r))
    \leq
    \frac1{q\cdot k+r}(q\cdot a(k))+
    \frac1n a(r)
    \leq
    l+\epsilon+\epsilon=l+2\epsilon
  \end{align*}
  The last inequality holds once $n$ is sufficiently large, specifically when
  \[
  n\geq\frac1\epsilon\max_{0\leq i\leq k}\;a(i)
  \]
  Therefore 
  \[
  \lim_{i \to \infty}\frac1i a(i)=l
  \]
  
  Now if $C>0$ then the sequence $b(i):=a(i)+C$ is subadditive and $\frac1i b(i)$
  converges by the previous argument. Thus we have
  \[
  \lim_{i \to \infty} \frac1i b(i)=\lim_{i \to \infty}\frac1i (a(i)+C)=\lim_{i \to \infty}\frac1i a(i)
  \]
\end{Proof}

\repeatclaim{p:adistonql} 
\begin{Proof} 
  Suppose $\bar\gamma_{1}$ and $\bar\gamma_{2}$ are two quasi-linear
  sequences of elements of $\Gamma$, then for any $i,j\in\Nbb_{0}$

  \begin{align*}
    \dist\big(\gamma_{1}(i+j),&\gamma_{2}(i+j)\big)
    \\
    &\leq
    \dist\big(\gamma_{1}(i+j),\gamma_{1}(i)\otimes\gamma_{1}(j)\big) +
    \dist\big(\gamma_{2}(i+j),\gamma_{2}(i)\otimes\gamma_{2}(j)\big) \\
    &\;\;\;+
    \dist\big(\gamma_{1}(i)\otimes\gamma_{1}(j),
            \gamma_{2}(i)\otimes\gamma_{2}(j)\big)\\
            &\leq
            \defect_{\dist}(\bar\gamma_{1}) +
            \defect_{\dist}(\bar\gamma_{2}) +
            \dist\big(\gamma_{1}(i),\gamma_{2}(i)\big) +
            \dist\big(\gamma_{1}(j),\gamma_{2}(j)\big)
  \end{align*}
  
  Thus the sequence $\dist(\gamma_{1}(i),\gamma_{2}(i))$ is
  quasi-subadditive and by Lemma~\ref{p:subadditive} the limit 
  \[
  \lim_{i\to\infty} \frac1i \dist\big(\gamma_{1}(i),\gamma_{2}(i)\big)
  \]
  exists and is finite.
\end{Proof}

\repeatclaim{p:boundeddefect} 
\begin{Proof}
  Given a Cauchy sequence $\set{\bar\gamma_{i}}$ of elements in
  $(\qlin_{\dist}(\Gamma),\dista)$ we need to find a limiting element
  $\bar\phi\in\qlin_{\dist}(\Gamma)$.  We will do that by a version of
  the diagonal process, that is we define $\phi(n)$ to
  have value $\gamma_{i}(n)$ for $i$ sufficiently large depending on
  $n$. The quasi-linearity of $\bar\phi$ would follow from the fact
  that for a fixed $n$ and all sufficiently large $i$ the set
  $\set{\gamma_{i}(n)}$ is uniformly bounded.

  Now we give the detailed argument.  First we replace each element
  of the sequence $\set{\bar\gamma_{i}}$ by an asymptotically
  equivalent element with defect bounded by the constant $C$
  according to assumption~(\ref{pi:boundeddefect}) of the lemma. We
  will still call the new sequence $\set{\bar\gamma_{i}}$. The Cauchy
  sequence $\set{\bar\gamma_{i}}$ satisfies
  \[
  \sup_{i,j\geq\ibf}\dista(\bar\gamma_{i},\bar\gamma_{j})\to0
  \quad\text{as}\quad
  \ibf\to\infty
  \]
    
  By assumption~(\ref{pi:boundeddefect}) of the lemma for any
  $n,k\in\Nbb_{0}$ holds
  \begin{align*}
    k \cdot \dist\big(\gamma_{i}(n),\gamma_{j}(n)\big)
    &=
    \dist\big(\gamma_{i}(n)^{\otimes k},\gamma_{j}(n)^{\otimes k}\big)\\
    &\leq
    \dist\big(\gamma_{i}(kn),\gamma_{j}(kn)\big) + 2 k \cdot C
  \end{align*}
  Dividing by $k$ we obtain
  \[
  \dist(\gamma_{i}(n),\gamma_{j}(n)) \leq \frac1k\dist(\gamma_{i}(kn),\gamma_{j}(kn)) + 2C
  \]
  Now we pass to the limit sending $k$ to infinity, while keeping $n$
  fixed:
  \[
  \dist(\gamma_{i}(n),\gamma_{j}(n))
  \leq
  n\cdot\dista(\bar\gamma_{i},\bar\gamma_{j}) + 2C
  \]
  Given $n$ let $\ibf(n)$ be a number such that for any $i,j\geq\ibf(n)$
  holds
  \[
  \dista(\bar\gamma_{i},\bar\gamma_{j})\leq \frac1n
  \]
  
  We may assume that $\ibf(n)$ is nondecreasing as a function of $n$.
  Then for any $i,j,n\in\Nbb$ with $i,j\geq \ibf(n)$ we have
  the following bound
  \[\tageq{boundedmembers}
    \dist\big(\gamma_{i}(n),\gamma_{j}(n)\big)
  \leq
  2C+1
  \]
  
  Now we are ready to define the limiting sequence $\bar\phi$ by
  setting
  \[
  \phi(n):=\gamma_{\ibf(n)}(n)
  \]
  First we verify that $\bar\phi$ is quasi-linear
  \begin{align*}
    \dist
    &\big(
      \phi(n+m),
      \phi(n)\otimes\phi(m)
    \big)
    =
    \dist
    \big(
      \gamma_{\ibf(n+m)}(n+m),
      \gamma_{\ibf(n)}(n)\otimes\gamma_{\ibf(m)}(m)
    \big)
    \\
    &\leq
    \dist
    \big(
      \gamma_{\ibf(n+m)}(n+m),
      \gamma_{\ibf(n+m)}(n)\otimes\gamma_{\ibf(n+m)}(m)
    \big)
    \\
    &\quad+
    \dist
    \big(
      \gamma_{\ibf(n+m)}(n)\otimes\gamma_{\ibf(n+m)}(m),
      \gamma_{\ibf(n)}(n)\otimes\gamma_{\ibf(m)}(m)
    \big)
    \\
    &\leq  
    C+
    \dist
    \big(
      \gamma_{\ibf(n+m)}(n),
      \gamma_{\ibf(n)}(n)
    \big)
    +
    \dist
    \big(
      \gamma_{\ibf(n+m)}(m),
      \gamma_{\ibf(m)}(m)
    \big)
    \\
    &\leq
    C+2(2C+1)=5C+2=:C'
  \end{align*}
  
  The convergence of $\bar\gamma_{i}$ to $\bar\phi$ is shown as
  follows. For $n,k\in\Nbb$ let $q,r\in\Nbb_{0}$ be the quotient and
  the remainder of the division of $n$ by $k$, that is $n=q\cdot k+r$
  and $0\leq r<k$.  Fix $k\in\Nbb$ and let $i\geq\ibf(k)$, then
  \begin{align*}
    \dista(\bar\gamma_{i},\bar\phi)
    &=
    \lim_{n\to\infty}\frac1n\dist\big(\gamma_{i}(n),\phi(n)\big)\\
    &=
    \lim_{n\to\infty}
    \frac1n
    \dist\big(\gamma_{i}(q\cdot k+r),\gamma_{\ibf(n)}(q\cdot k+r)\big)\\
    &\leq
    \lim_{n\to\infty}
    \frac1n
    \left(\rule{0mm}{5mm}
    q\cdot\dist\big(\gamma_{i}(k),\gamma_{\ibf(n)}(k)\big)+
    \dist\big(\gamma_{i}(r),\gamma_{\ibf(n)}(r)\big) +
    2qC' + 2C'
    \right)\\
    &\leq
    \lim_{n\to\infty}
    \frac1n
    \big((3q+3)\cdot C'\big)\\
    &=
    \frac{3C'}{k}      
  \end{align*}
  Since $k\in\Nbb$ is arbitrary we have 
  \[
  \lim_{i\to\infty}\dista(\bar\gamma_{i},\bar\phi)=0
  \]
\end{Proof}

\repeatclaim{p:eps-linear-dense}
\begin{Proof}
	Let $\bar\gamma=\set{\gamma(n)}$ be a $\dist$-quasi-linear
	sequence. We need to approximate it with linear sequences. For
	$i\in\Nbb$, let
	$\bar\gamma_{i}$ be a sequence asymptotically equivalent to $\gamma$
	and satisfying
	\[
	\defect_{\dist}\bar\gamma_{i}\leq 1/i
	\]
	as provided by the $1/i$-uniformly bounded defect property.
	
	Define a $\dist$-linear sequence $\bar\eta_{i}$ by
	\[
	\eta_{i}(n):=\gamma_{i}(1)^{\otimes n}
	\]
	Then
	\begin{align*}
	\dista(\bar\gamma,\bar\eta_{i})
	&=
	\dista(\bar\gamma_{i},\bar\eta_{i})
	\\
	&=
	\lim_{n\to\infty}
	\frac1n
	\dist(\gamma_{i}(n),\eta_{i}(n))
	\\
	&=
	\lim_{n\to\infty}
	\frac1n
	\dist(\gamma_{i}(n),\gamma_{i}(1)^{\otimes n})
	\\
	&\leq
	\lim_{n\to\infty}
	\frac1n
	\cdot n\cdot\defect_{\dist}(\bar\gamma_{i})
	\\
	&\leq
	\frac1i    
	\end{align*}
	Thus $\lim\bar\eta_{i}=\gamma$.
\end{Proof}


\repeatclaim{p:dist-dista-isometry}

\begin{Proof}
  Let $\bar\gamma_{1},\bar\gamma_{2}\in\qlin_{\dist}(\Gamma)$ be two
  sequences of $\dist$-quasi-linear sequences.  We have to show that the
  two numbers
  \[
  \dista(\bar\gamma_{1},\bar\gamma_{2})
  =
  \lim_{n\to\infty}\frac1n \dist\big(\gamma_{1}(n),\gamma_{2}(n)\big)
  \]
  and
  \[
  \distaa(\bar\gamma_{1},\bar\gamma_{2})
  =
  \lim_{n\to\infty}\frac1n \dista\big(\gamma_{1}(n),\gamma_{2}(n)\big)
  \]
  are equal.  Since shifts are non-expanding maps, we have
  $\dista\leq\dist$ and it follows immediately that
  \[
  \distaa(\bar\gamma_{1},\bar\gamma_{2})
  \leq
  \dista(\bar\gamma_{1},\bar\gamma_{2})
  \]
  and we are left to show the opposite inequality.
  We will do it as follows. Fix $n>0$, then
  \begin{align*}
    \dista(\bar\gamma_{1},\bar\gamma_{2})
    &=
    \lim_{k\to\infty}\frac{1}{kn}\dist\big(\gamma_1(kn),\gamma_{2}(kn)\big)\\
    &\leq
    \lim_{k\to\infty}\frac{1}{kn}
    \bigg(
    \dist\big(\gamma_1(n)^{\otimes k},\gamma_{2}(n)^{\otimes k}\big)
    +
    k\cdot\big(
    \defect_{\dist}(\bar\gamma_{1})
    +
    \defect_{\dist}(\bar\gamma_{2})
    \big)
    \bigg)\\
    &\leq
    \frac1n\dista\big(\gamma_{1}(n),\gamma_{2}(n)\big)
    +
    \frac1n\big(\defect_{\dist}(\bar\gamma_{1})
    +
    \defect_{\dist}(\bar\gamma_2)\big)
  \end{align*}
  Passing to the limit with respect to $n$ gives required inequality
  \[
  \dista(\bar\gamma_{1},\bar\gamma_{2})
  \leq
  \distaa(\bar\gamma_{1},\bar\gamma_{2})
  \]
\end{Proof}


\repeatclaim{p:dist-dista-dense}
\begin{Proof}
	Given an element $\bar\gamma=\set{\gamma(n)}$ in
	$\qlin_{\dista}(\Gamma)$ we have to find a $\distaa$-approxi\-ma\-ting
	sequence $\bar\gamma_{i}=\set{\gamma_{i}(n)}$ in
	$\qlin_{\dist}(\Gamma)$.
	Define 
	\[
	\gamma_{i}(n):=\gamma(i)^{\otimes\lfloor\frac{n}{i}\rfloor}
	\]
	We have to show that each $\bar\gamma_{i}$ is $\dist$-quasi-linear and
	that $\distaa(\bar\gamma_{i},\bar\gamma)\too[i\to\infty]0$.
	These follow from 
	\begin{align*}
	\dista\big(\gamma_i(m+n), \gamma_i(m) \otimes \gamma_i(n) \big)
	&=\dista\left( \gamma(i)^{\otimes \lfloor \frac{m+n}{i}\rfloor} , 
	\gamma(i)^{\otimes \lfloor \frac{m}{i}\rfloor } \otimes \gamma(i)^{\otimes \lfloor \frac{n}{i} \rfloor}\right)\\
	&\leq \dista \left( \gamma(i), \bm 1 \right)
	\end{align*} 
	and
	\begin{align*}
	\distaa(\bar\gamma_i,\bar\gamma)
	&= \lim_{n \to \infty} \frac1n\dista\left( \gamma_i(n), \gamma(n) \right) \\
	&= \lim_{n \to \infty} \frac1n \dista \left( \gamma(i)^{\otimes \lfloor\frac{n}{i}\rfloor}, \gamma(n) \right)\\
	&\leq \lim_{n \to \infty} \left[\frac1n \dista \left( \gamma\left(i\lfloor\tfrac{n}{i}\rfloor\right), \gamma(n) \right) + \frac1n \lfloor \tfrac{n}{i} \rfloor \defect_{\dista}(\bar\gamma) \right]\\
	&\leq \lim_{n \to \infty}\left[ \frac1n \max_{k=0,\ldots,i-1} \dista\left(\bm 1, \gamma(k) \right)
	 + \frac{i}{n} \defect_{\dista} (\bar\gamma) \right]
	 + \frac{1}{i} \defect_{\dista}(\bar\gamma)\\
	 &\leq \frac{1}{i} \defect_{\dista}(\bar \gamma)
	\end{align*}
  It is worth noting, that the defect of $\bar\gamma_{i}$ need not to
	be uniformly bounded with respect to $i$.
\end{Proof}

\repeatclaim{p:unifsmalldefectaikd}
\begin{Proof}
  Let $\bar\Xcal=\set{\Xcal(i)}$ be a quasi-linear sequence and let $\epsilon > 0$. We will
  find an asymptotically equivalent sequence with defect less than $\epsilon$.
  
  Define a
  new sequence $\bar\Ycal=\set{\Ycal(i)}$ by
  \[
  \Ycal(i)
  :=
  \big[\Xcal(k\cdot i)\big]
  \oplus_{\Lambda_{1/k}}
  \set{\bullet}
  \]
  where the number $k\in\Nbb$ will be chosen later.  First we verify
  that the sequences $\bar\Xcal$ and $\bar\Ycal$ are asymptotically
  equivalent, that is 
  \begin{align*}
    \hat\aikd(\bar\Xcal,\bar\Ycal)
    &:=
    \lim_{i\to\infty}
    \frac{1}{i}
    \aikd\left(\Xcal(i), 
               \Ycal(i)
         \right)
     =
     0
  \end{align*}
  We estimate the asymptotic distance between individual members of
  sequences $\bar\Xcal$ and $\bar\Ycal$ using
  Lemma~\ref{p:mixtures-dist1} as
  follows
  \begin{align*}
    \aikd(&\Xcal(i), 
               \Ycal(i)
         )
    =
    \aikd\big(\Xcal(i), 
               \Xcal(k\cdot i)
               \oplus_{\Lambda_{1/k}}\set{\bullet}
         \big)    
    \\
    &\leq
    \aikd\left(\Xcal(i), 
               \Xcal(i)^{\otimes k}
               \oplus_{\Lambda_{1/k}}\set{\bullet}
         \right)
    +    
    \aikd\left(\Xcal(i)^{\otimes k}
               \oplus_{\Lambda_{1/k}}\set{\bullet}, 
               \Xcal(k\cdot i)
               \oplus_{\Lambda_{1/k}}\set{\bullet}
         \right)    
    \\
    &\leq
    \ent(\Lambda_{1/k})
    +    
    \defect_{\aikd}(\bar\Xcal)   
  \end{align*}
  Thus $\hat\aikd(\bar\Xcal,\bar\Ycal)=0$ and the two sequences are
  asymptotically equivalent.
  
  Next we show that the sequence $\bar\Ycal$ is $\aikd$-quasi-linear and evaluate
  its defect using Lemma~\ref{p:mixtures-dist1}. Let $i,j\in\Nbb$, then
  \begin{align*}
    \aikd&\big(\Ycal(i+j),\Ycal(i)\otimes\Ycal(j)\big)
    \\
    &=
    \aikd\Big(\Xcal(k\cdot i+k\cdot j)
                \oplus_{\Lambda_{1/k}}\set{\bullet},
              \big[
                \Xcal(k\cdot i)
                \oplus_{\Lambda_{1/k}}\set{\bullet}
              \big]
              \otimes
              \big[
                \Xcal(k\cdot j)
                \oplus_{\Lambda_{1/k}}\set{\bullet}
              \big]
         \Big)    
    \\
    &\leq
    \aikd\Big(\big[\Xcal(k\cdot i)\!\otimes\!\Xcal(k\cdot j)\big]
                   \oplus_{\Lambda_{1/k}}\!\set{\bullet},     
              \big[
                \Xcal(k\cdot i)
                \oplus_{\Lambda_{1/k}}\!\set{\bullet}
              \big]
              \otimes
              \big[
                \Xcal(k\cdot j)
                \oplus_{\Lambda_{1/k}}\!\set{\bullet}
              \big]
         \Big)    
     \\
     &\quad+
     \frac1k \defect(\bar\Xcal)
     \\
     &\leq
     3\ent(\Lambda_{1/k})+\frac1k \defect_{\aikd}(\bar\Xcal)
  \end{align*}
  Thus, by choosing $k$ to be a solution to the inequality
  \[
  3\ent(\Lambda_{1/k})+\frac1k \defect_{\aikd}(\bar\Xcal)\leq\epsilon
  \]
  we can make sure that 
  \[
  \defect_{\aikd}(\bar\Ycal)\leq\epsilon
  \]
\end{Proof}

\repeatclaim{p:ikduniform}
\begin{Proof} 
  Consider a two-fan $U_{n}\stackrel{f}{\ot} U_{nm}\stackrel{g}{\to}
  U_{m}$. To construct specific reductions $f$ and $g$ we identify
  $U_{nm}$, $U_{n}$ and $U_{m}$ with the cyclic groups of the
  corresponding order
\begin{align*}
  U_{nm}&\leftrightarrow \Zbb_{nm}\\
  U_{n}&\leftrightarrow\Zbb_{n}\\
  U_{m}&\leftrightarrow\Zbb_{m}
\end{align*}

Consider the short exact sequences
\begin{align*}
  \set{0} \too \Zbb_{n}
  &\stackrel{\times m}{\too}
  \Zbb_{nm}
  \stackrel{\mod m}{\too}
  \Zbb_{m} \too \set{0}\\
  \set{0} \too \Zbb_{m}
  & \stackrel{\times n}{\too}
  \Zbb_{nm}
  \stackrel{\mod n}{\too}
  \Zbb_{n}\too \set{0}\\  
\end{align*}
Choose for $f$ the left splitting in the first exact sequence, and for $g$ the left splitting in the second exact sequence.

Now that we constructed a two-fan $U_{n}\stackrel{f}{\ot} U_{nm}\stackrel{g}{\to}
U_{m}$, let $U_{n}\ot Z\to U_{m}$ be its minimal reduction. Now we estimate
$|Z|\leq n+m$, which implies that
\begin{align*}
  \ikd(U_{n},U_{m})
  &\leq
  2\Ent(Z)-\ent(U_{n})-\ent(U_{m})\\
  &\leq
  2\ln(n+m)-\ln n-\ln m\\
  &\leq
  2\ln2+2\ln\max\set{n,m}-\ln n-\ln m\\
  &=2\ln2+\left|\ln\frac{n}{m}\right|
\end{align*}

To prove the second assertion note that entropy is a
$\ikd$-1-Lipschitz function. Therefore we have
\[
|\ent(U_{n})-\ent(U_{m})|\leq\ikd(U_{n},U_{m})
\leq
|\ent(U_{n})-\ent(U_{m})|+2\ln2
\]
Substituting in the definition of asymptotic Kolmogorov distance we
obtain the required equality.
\end{Proof}


\let\thesubsection=\thesubsectionstandard

	\let\thesection=\standardthesection

	\let\section=\Section
	\bibliographystyle{alpha}
	\bibliography{ReferencesEntropy}

\newcommand{\etalchar}[1]{$^{#1}$}
\begin{thebibliography}{ABK{\etalchar{+}}15}

\bibitem[ABD{\etalchar{+}}08]{Ay-Predictive-2008}
Nihat Ay, Nils Bertschinger, Ralf Der, Frank G{\"u}ttler, and Eckehard Olbrich.
\newblock Predictive information and explorative behavior of autonomous robots.
\newblock {\em The European Physical Journal B}, 63(3):329--339, 2008.

\bibitem[ABK{\etalchar{+}}15]{Abramsky-Contextuality-2015}
Samson Abramsky, Rui~Soares Barbosa, Kohei Kishida, Raymond Lal, and Shane
  Mansfield.
\newblock Contextuality, cohomology and paradox.
\newblock {\em arXiv preprint arXiv:1502.03097}, 2015.

\bibitem[BBI01]{Burago-Course-2001}
Dmitri Burago, Yuri Burago, and Sergei Ivanov.
\newblock {\em A course in metric geometry}, volume~33 of {\em Graduate Studies
  in Mathematics}.
\newblock American Mathematical Society, Providence, RI, 2001.

\bibitem[BFL11]{Baez-Characterization-2011}
John~C. Baez, Tobias Fritz, and Tom Leinster.
\newblock A characterization of entropy in terms of information loss.
\newblock {\em Entropy}, 13(11):1945--1957, 2011.

\bibitem[BRO{\etalchar{+}}14]{Bertschinger-Quantifying-2014}
Nils Bertschinger, Johannes Rauh, Eckehard Olbrich, J\"urgen Jost, and Nihat
  Ay.
\newblock Quantifying unique information.
\newblock {\em Entropy}, 16(4):2161--2183, 2014.

\bibitem[Csi98]{Csiszar-Method-1998}
Imre Csisz\'ar.
\newblock The method of types.
\newblock {\em IEEE Trans. Inform. Theory}, 44(6):2505--2523, 1998.
\newblock Information theory: 1948--1998.

\bibitem[CT91]{Cover-Elements-1991}
Thomas~M. Cover and Joy~A. Thomas.
\newblock {\em Elements of information theory}.
\newblock Wiley Series in Telecommunications. John Wiley \& Sons, Inc., New
  York, 1991.
\newblock A Wiley-Interscience Publication.

\bibitem[DFZ11]{Dougherty-Non-Shannon-2011}
Randall Dougherty, Chris Freiling, and Kenneth Zeger.
\newblock Non-shannon information inequalities in four random variables.
\newblock {\em arXiv preprint arXiv:1104.3602}, 2011.

\bibitem[Fri09]{Friston-Free-2009}
Karl Friston.
\newblock The free-energy principle: a rough guide to the brain?
\newblock {\em Trends in cognitive sciences}, 13(7):293--301, 2009.

\bibitem[Gro12]{Gromov-Search-2012}
Misha Gromov.
\newblock In a search for a structure, part 1: On entropy.
\newblock Preprint available at \url{http://www. ihes. fr/gromov}, 2012.

\bibitem[Mat07]{Matus-Infinitely-2007}
Frantisek Matus.
\newblock Infinitely many information inequalities.
\newblock In {\em Information Theory, 2007. ISIT 2007. IEEE International
  Symposium on}, pages 41--44. IEEE, 2007.

\bibitem[SA15]{Steudel-Information-2015}
Bastian Steudel and Nihat Ay.
\newblock Information-theoretic inference of common ancestors.
\newblock {\em Entropy}, 17(4):2304--2327, 2015.

\bibitem[Sin76]{Sinai-ergodic-1976}
Ya.~G. Sinai.
\newblock {\em Introduction to ergodic theory}.
\newblock Princeton University Press, Princeton, N.J., 1976.
\newblock Translated by V. Scheffer, Mathematical Notes, 18.

\bibitem[VDP13]{Dijk-Informational-2013}
Sander~G Van~Dijk and Daniel Polani.
\newblock Informational constraints-driven organization in goal-directed
  behavior.
\newblock {\em Advances in Complex Systems}, 16(02n03):1350016, 2013.

\bibitem[Yeu12]{Yeung-First-2012}
Raymond~W Yeung.
\newblock {\em A first course in information theory}.
\newblock Springer Science \& Business Media, 2012.

\end{thebibliography}
	
\end{document}